\documentclass[12pt]{article}

\def\calx{{\mathcal{X}}}
\def\calz{{\mathcal{Z}}}

\def\calu{{\mathcal{U}}}
\def\calv{{\mathcal{V}}}

\def\({\left(}
\def\){\right)}

\def\vsp{\vspace*{1,5mm}\\ }

\def\bk{\bigskip }

\def\sk{\smallskip }

\def\n{\noindent }
\def\dd{\displaystyle}

\def\D{{\Delta}}

\def\barr{\begin{array}}
\def\earr{\end{array}}

\def\bit{\begin{itemize}}

\def\eit{\end{itemize}}

\def\Re{{\rm Re}}

\def\D{{\Delta}}

\usepackage{amsmath, amsthm, amscd, amsfonts, amssymb}
\usepackage[usenames]{color}
\usepackage{mathrsfs}

\numberwithin{equation}{section} %公式随章节编号%

\newtheorem{theorem}{Theorem}[section]
\newtheorem{proposition}[theorem]{Proposition}

\newtheorem{corollary}[theorem]{Corollary}

\newtheorem{lemma}[theorem]{Lemma}
\theoremstyle{definition}
\newtheorem{definition}[theorem]{Definition}

\newtheorem{remark}[theorem]{Remark}

\def\S{Schr\"o\-din\-ger}

\def\calf{{\mathcal{F}}}
\def\caln{\mathcal{N}}

\def\calx{\mathcal{X}}
\def\calz{\mathcal{Z}}

\def\calv{{\mathcal{V}}}
\def\caly{{\mathcal{Y}}}

\def\bbc{{\mathbb{C}}}
\def\bbe{{\mathbb{E}}}
\def\bbn{{\mathbb{N}}}
\def\bbr{{\mathbb{R}}}
\def\bbp{{\mathbb{P}}}
\def\bbx{{\mathbb{X}}}

\def\1{^{-1}}
\def\vsp{\vspace*{2mm}\\ }

\def\calf{{\mathcal{F}}}

\def\calx{{\mathcal{X}}}
\def\E{{\mathbb{E}}}
\def\rr{{\mathbb{R}}}
\def\nn{{\mathbb{N}}}
\def\9{{\infty}}

\def\lbb{{\lambda}}
\def\a{{\alpha}}
\def\b{{\beta}}
\def\na{{\nabla}}
\def\g{{\gamma}}
\def\wt{\widetilde}

\def\vf{{\varphi}}
\def\oo{{\omega}}
\def\ooo{{\Omega}}
\def\pp{{\partial}}
\def\p{{\partial}}
\def\D{{\Delta}}
\def\vp{{\varepsilon}}
\def\barr{\begin{array}}
\def\earr{\end{array}}

\def\dd{\displaystyle}
\def\bk{\bigskip }
\def\sk{\smallskip}
\def\n{\noindent }

\def\pas{\mathbb{P}\mbox{-a.s.}}

\def\vsp{\vspace*{2mm}\\ }
\def\ff{\forall }
\def\({\left(}
\def\){\right)}
\def\<{\left<}
\def\>{\right>}

\def\wt{\widetilde}
\def\wh{\widehat}
\def\ol{\overline}
\def\ve{{\varepsilon}}

\begin{document}

\begin{center}
{\Large{\bf Optimal bilinear control of stochastic nonlinear Schr\"odinger equations: \\
mass-(sub)critical case}}
\bigskip\bk

{
{\large{\bf Deng Zhang}}\footnote{Department of Mathematics,
Shanghai Jiao Tong University, 200240 Shanghai, China.
Email address: dzhang@sjtu.edu.cn  }}
\end{center}

\bk\bk\bk

\begin{quote}
\n{\small{\bf Abstract.}
We study optimal bilinear control problems for stochastic nonlinear Schr\"odinger equations
in both the mass subcritical and critical case.
For general initial data of the minimal $L^2$ regularity,
we prove the existence and first order Lagrange condition of an open loop control.
Furthermore,
we obtain uniform estimates of (backward) stochastic solutions
in new spaces of type $U^2$ and $V^2$,
adapted to evolution operators  related to linear Schr\"odinger equations
with lower order perturbations.
In particular,
we obtain a new temporal regularity of rescaled (backward) stochastic solutions,
which is the key ingredient in the proof of tightness of approximating controls
induced by Ekeland's variational principle.
}


\end{quote}

\vfill

\section{Introduction} \label{Sec-Intro}

We are concerned with the controlled stochastic system governed by the nonlinear \S\ equation
\begin{equation}\label{equa-x}
\barr{rcll}
idX(t,x)&\!\!=\!\!&\D X(t,x)dt+\lbb|X(t,x)|^{\a-1}X(t,x)dt -i\mu(x)X(t,x)dt \\
&&    +V_0(x)X(t,x)dt+\dd\sum\limits^m_{j=1} u_j(t)V_j(x)X(t,x)dt  \\
&&\dd+iX(t,x)dW(t,x),\ \ \ t\in(0,T), \ x\in\rr^d,\vsp
X(0)&\!\!=\!\!&X_0\ \mbox{ in }\rr^d. \earr
\end{equation}
Here,
$\lbb =-1$ (resp. $\lbb =1$) corresponds to the defocusing (resp. focusing) case,
$\a>1$, $V_j\in L^{\9}(\rr^d)$ are real valued functions, $0\leq j\leq m$,
$u(t) = (u_1(t),...,u_m(t)) \in \bbr^m$ represents the control,
$W$ is the colored Wiener process,
\begin{align} \label{W}
    W(t,x)=&\sum^N_{j=1}\mu_j e_j(x)\b_j(t),\ \ t\ge0,\ x\in\rr^d,
\end{align}
and
\begin{align*}
   \mu(x) =& \dd\frac12\sum^N_{j=1}|\mu_j|^2e^2_j(x),\ \  x\in\rr^d,\ d\ge1,
\end{align*}
where $\mu_j$  are purely imaginary numbers (i.e. $Re \mu_j =0$),
$e_j$ are real-valued functions,
and $\b_j$  are independent real Brownian motions on a probability space
$(\Omega,\calf,\mathbb{P})$ with natural filtration
$(\calf_t)_{t\ge0}$, $1\le j\le N.$
For simplicity, we assume $N<\9$,
but the arguments in this paper can be extended to the case where $N= \9$
under appropriate summable conditions.

In a quantum mechanical interpretation,
$X=X(t,x,\oo),$ $x\in\rr^d$, $t\ge0,$ $\oo\in\ooo$,
represents the quantum state at time $t$,
while the stochastic perturbation $iX\,dW$
represents a stochastic continuous measurement via the pointwise quantum observables
$R_j(X)=\mu_je_jX$
(see \cite[Section 2]{BG09}).
The functions $V_j$, $0\leq j\leq m$, describe an external potential.
In most situations
the control $u$ represents an external applied force
due to the interaction of the quantum system with an electric field or a laser pulse applied to a quantum system.

Here,
we are mainly concerned with {\it the conservative case}
where $\Re \mu_j =0$,
$1\leq j\leq N$.
In this case,
$-i\mu X dt + i X dW$ is indeed the Stratonovitch differential,
and so, via It\^o's formula,
one has the pathwise conservation $|X(t)|^2_{L^2}=|X_0|^2_{L^2}$, $\ff t\ge0$.
Hence,
with the normalized initial state $|X_0|_{L^2}=1$,
the quantum system evolves on the unit ball of $L^2$ and verifies the conservation of probability.
In particular,
in the focusing case $\lbb =1$,
stochastic Schr\"odinger equations with cubic nonlinearity
in dimension two arises in molecular aggregates with thermal fluctuations
(\cite{BCIR94,BCIRG95}).
See also \cite{RGBC95} for
stochastic Schr\"odinger equations with quartic nonlinearity
in dimension one.
Actually,
in both cases
the nonlinearities are critical when the initial
data are of only $L^2$ regularity.

We would like also to mention that
{\it the non-conservative case} (i.e., $\Re \mu_j \not =0$ for some $1\leq j\leq N$)
plays an important role
in the theory of open quantum systems.
See \cite{BG09, BH95,BPZ98} and references therein.

For the global well-posedness of \eqref{equa-x} in the mass-subcritical case
(i.e., $1<\a<1+\frac 4d$),
see
\cite{BRZ14,BRZ16,BD99,BD03,HRZ18,H18}.
See also \cite{BM13} for the compact manifold setting,
and \cite{BHW17,BHW18,BHM18.2,BH17,BH18}
for the study of martingale solutions.

Optimal control problem and exact controllability of
Schr\"odinger equations
have been extensively studied in the deterministic case.
See, e.g., \cite{BCMR07,HMMS13,IK07,MRT05}.
In the stochastic situation,
there are many results
on optimal control problems of dissipative equations,
see, e.g., \cite{FO16,FT02}.

However,
less results of optimal control problems are known for stochastic Schr\"odinger equations,
which is of dispersive type.
One of main difficulties lies in the weak regularization effect of free Schr\"odinger group $\{{e^{-it\Delta}}\}$
which,
unlike the heat semigroup,
cannot raise global Sobolev regularity.
Another difficulty arises from the defect of compact embeddings
in the probability space.
More precisely,
even if a space $\caly$ is compactly
imbedded into another space $\mathcal{Z}$, one generally does not have
the compact imbedding from $L^p(\Omega; \caly)$
to $L^p(\Omega; \mathcal{Z})$, $1\leq p\leq \9$.

For optimal control of
stochastic Schr\"odinger equations with Lipschitz nonlinearity and additive noise,
we refer to \cite{K13,K15}.
For the stochastic nonlinear Schr\"odinger equation \eqref{equa-x},
the existence  of an open-loop optimal control
and first-order Lagrange optimality conditions
were proved in the recent work \cite{BRZ18}
for a part of mass subcritical exponents
where $2\leq \a<1+\frac 4d$ with $1\leq d\leq 3$.
The primary idea of \cite{BRZ18} is
based on Ekeland's variational principle (\cite{E74,E79}) and
the subdifferential in the sense of Rockafellar (\cite{R79}),
which enable one
to obtain geometric characterizations of approximating controls
and so the tightness of associated distributions
by using appropriate estimates
of controlled solutions
and related dual backward stochastic solutions.

In the present work,
we are mainly interested in the
optimal control problems for \eqref{equa-x}
in the defocusing mass-critical case 
($\lbb =-1$, $\a=1+\frac 4d$, $d\geq 1$).

It should be mentioned that,
the global well-posedness and scattering of deterministic nonlinear Schr\"odinger equations
in the defocusing  case was
one of main conjectures in the field of dispersive equations.
The mass-critical case with general initial data has been proved in the recent work by Dodson \cite{D12,D16.1,D16.2},
where one main ingredient is the long-time Strichartz estimate.
We would also like to mention \cite{CKSTT08,RV07,V07}
for the proof of this conjecture in the defocusing energy-critical case.
In the stochastic case,
the global well-posedness
was proved in \cite{FX18.1,FX18.2}
for the defocusing mass-critical case with dimension $d=1$.
In the work \cite{Z18},
the global well-posedness of \eqref{equa-x}
has been proved in the defoucsing mass-critical case for all dimensions $d\geq 1$.
Moreover,
the defocusing energy-critical case,
scattering as well as
Stroock-Varadhan type support theorem
have been also studied in \cite{Z18}.

Here we give a unified treatment of optimal bilinear control problems for \eqref{equa-x}
in both the mass subcritical and critical case.
In particular,
the results apply to the cubic and quartic
nonlinear Schr\"odinger equations in dimensions $d=2$ and $d=1$, respectively.
Moreover,
the results are also applicable to
the whole subcritical range of exponents of nonlinearity
including the cases $\a\in (1,2)$ with $1\leq d\leq 3$ and  $\a\in (1,1+\frac 4d)$ with $d\geq 4$,
which were previously excluded in  \cite{BRZ18}.

Another novelty of the present work is to treat the situation where
the initial data of \eqref{equa-x} are of only $L^2$ regularity,
which is actually the minimal regularity required to define the objective functional
and is less regular than the $H^1$ regularity assumed in \cite{BRZ18}.
The key role here is played by the spaces $U^p$ and $V^p$,
which were first introduced by Koch and Tataru \cite{KT05,KT07}
and have been very effective in the study of dispersive equations.
See e.g. \cite{BOP15,HHK09,HTT11,KTV14} and references therein.

We obtain uniform estimates of solutions to \eqref{equa-x}
and dual backward stocastic equation (see \eqref{equa-back} below)
in new spaces of type $U^2$ and $V^2$,
adapted to evolution operators related to
linear Schr\"odinger equations with lower order perturbations.
These estimates are sharper than those in the Strichartz spaces  obtained in \cite{BRZ18}
and, more importantly,
reveal a new temporal regularity of rescaled (backward) stochastic solutions,
which is the key ingredient to obtain the tightness of approximating controls
induced by Ekeland's variational principle.

We would also like to mention that
the rescaling approach, developed in \cite{BRZ14,BRZ16, BRZ18},
also plays an important role
in  constructing appropriate spaces
for this new temporal regularity of  (backward) stochastic solutions.
This, in spirit, has similarity with the work \cite{BR14},
where a new temporal regularity of stochastic solutions has been obtained
for a large class of stochastic dissipative partial differential equations.

\paragraph{Notations.}
We denote by $|\cdot|_m$ the Euclidean norm in $\rr^m$ and by $v\cdot w$ the scalar product of vectors $v,w\in\bbr^m$.
For any $x=(x_1,\cdots,x_d) \in \bbr^d$
and any multi-index $\a=(\a_1,\cdots, \a_d)$,
$|\g|= \sum_{j=1}^d \g_j$,
$\<x\>=(1+|x|^2)^{1/2}$,
$\partial_x^\g=\partial_{x_1}^{\g_1}\cdots \partial_{x_d}^{\g_d}$,
and
$\<\na\>=(I-\Delta)^{1/2}$.

For $1\le p\le\9$,
$L^p = L^p(\rr^d)$ is the space of $p$-integrable (complex-valued) functions,
endowed with the norm $|\cdot|_{L^p}$.
In particular, the Hilbert space $L^2(\rr^d)$ is endowed with the scalar product
$\<v,w\>_2=\int_{\rr^d} v(x)\bar w(x)dx$.
As usual,
$L^q(0,T;L^p)$ means the space of all $L^q(0,T)$-integrable $L^p$-valued functions  with the norm
$\|\cdot\|_{L^q(0,T;L^p)}$,
and $C([0,T];L^p)$ denotes the space of all $L^p$-valued continuous functions on $[0,T]$ with the sup norm in $t$.

For any $2\leq q<\9$,
we use the $U^p$ and $V^p$ spaces as in \cite{BOP15,HHK09,HTT11,KTV14}.
We also use the local smoothing space defined by, for $\a, \beta\in \bbr$,
$L^2(I; H^\a_\beta)
   =\{u\in \mathcal{D}'; \int_I\int \<x\>^{2\beta} |\<\na\>^\a u(t,x)|^2 dx dt<\9\}$,
where $\mathcal{D}'$ means distributions.

For any two Banach spaces $\mathcal{X}$ and $\mathcal{Y}$,
$\mathcal{L}(\calx, \caly)$ is the space of linear continuous operators from $\calx$ to $\caly$.
Throughout this paper, we use $C$ for various constants that may change from line to line.

\section{Formulation of main results} \label{Sec-Main}

To begin with,
we recall the definition of a strong solution to equation \eqref{equa-x}.

\begin{definition}\label{d2.1}\rm
Let $X_0\in L^2$,  $0<T<\9$. Let $\a$ satisfy $1<\a \leq 1+\frac 4 d$, $d\geq 1$.
A strong $L^2$-solution to \eqref{equa-x} on $[0,T]$
is an $L^2$-valued continuous $(\calf_t)_{t\ge0}$-adapted process $X$ such that $|X|^{\a-1}X\in L^1(0,T;H^{-1})$
and $\pas$
\begin{align}\label{e2.1}
X(t) =& X_0 -\int^t_0 \bigg( i\Delta X(s)+\mu X(s)+\lbb i|X(s)|^{\a-1} X(s) +iV_0X(s) \nonumber \\
      & \qquad +i\dd\sum^m_{j=1}u_j(s)V_j X(s) \bigg)ds
          + \int^t_0X(s)dW(s),\ t\in [0,T],
\end{align}
as an It\^o equation in $H^{-2}$ (resp. $H^{-1}$).
Here,
the last integral in \eqref{e2.1} above is taken in the sense of It\^o.
See, e.g., \cite{DZ12} and \cite{LR15}.
\end{definition}

We assume the asymptotically flat condition below as in \cite{BRZ18}
(see also \cite{BRZ14,BRZ16,Z17,Z18}),
mainly for the global well-posedness of \eqref{equa-x}.
\begin{enumerate}
   \item[(H0)]
  For each $1\leq j\leq N$ and each  $1\leq k\leq m$, $e_j, V_k\in C^\9_b(\rr^d)$ satisfy that
   for any multi-index $\g \not = 0$,
\begin{align} \label{decay}
   \lim_{|x|\to \9} \<x\>^2 |\pp_x^\g e_j(x)| + |\pp_x^\g V_k(x)|=0.
\end{align}
\end{enumerate}

A pair $(p,q)$ is called a Strichartz pair,
if $\frac 2q = d(\frac 12 - \frac 1p)$,
$(p,q) \in [2, \9] \times (2,\9]$.
\footnote{The endpoint case where $q=2$ is not considered here. }
For any interval $I \subseteq \bbr^+$,
define the Strichartz spaces by
\begin{align*}
   S^0(I) = \bigcap\limits_{(p,q):Strichartz\ pair} L^q(I; L^p), \ \
   N^0(I) = \bigcup\limits_{(p,q):Strichartz\ pair} L^{q'}(I; L^{p'}).
\end{align*}

The following hypothesis is assumed  for
the integrability of solutions to \eqref{equa-x}
in the defocusing mass-critical case.
\begin{enumerate}
   \item[$(H0)^*$]
  $\{e_j\}_{j=1}^N$ are constants with $d\geq 1$.
  Or,
  $\{e_j\}_{j=1}^N \subseteq L^{\frac{2p}{p-2}}(\bbr^d)$ for some
  $p$ satisfying $\frac 1p \in (\max\{\frac{1}{2\a}, \frac 12 - \frac{1}{2d}\}, \frac{1}{\a}(\frac 12 + \frac 1d))$
  with $\a=1+\frac 4d$ and
  $1\leq d\leq 3$.
\end{enumerate}
(Note that,
the interval $(\max\{\frac{1}{2\a}, \frac 12 - \frac{1}{2d}\}, \frac{1}{\a}(\frac 12 + \frac 1d))$ is nonempty
in dimensions $1\leq d\leq 3$.
The condition on $p$ ensures that
there exists another Strichartz pair $(\wt{p}, \wt{q})$
such that $(\frac{1}{\wt{p}'}, \frac{1}{\wt{q}'}) = (\frac{\a}{p}, \frac{\a}{q})$,
where $q\in (2,\9)$ is such that $(p,q)$ is a Strichartz pair,
and $\wt{p}', \wt{q}'$ denote the conjugate numbers of $\wt{p}, \wt{q}$ respectively.
The restriction $p<\frac{2d}{d-1}$ arises from the integrability \eqref{E-Xu-U2-Integ} below.)

Denote by $U(t,s)$, $t,s\in \bbr^+$,
the evolution operators related to $-i e^{-W}\Delta(e^W\cdot)$,
i.e.,
for any $v\in L^2$, $v(t):=U(t,s)v$
solves the linear Schr\"odinger equation with lower order perturbations below
\begin{align}
   & i\p_t v(t) = e^{-W}\Delta(e^W v(t)) (= \Delta v(t) + (b(t)\cdot \na + c(t))v(t)), \label{equa-v-homo}\\
   & v(s) = v,  \nonumber
\end{align}
where $b(t)=2\na W(t)$,
$c(t) = \Delta W(t) + \sum_{j=1}^d (\p_jW(t))^2$.
It is known (see \cite{D96}, see also Theorem \ref{Thm-Stri} below)  that
$\{U(t,s)\}$ are bounded operators in $L^2$.
In particular,
in the case where $\{e_j\}$ are constants,
$-i e^{-W}\Delta(e^W\cdot) = -i \Delta$,
and so $U(t,s)=e^{-i(t-s)\Delta}$.

In order to obtain the temporal regularity of controlled solutions to \eqref{equa-x},
we introduce a new space adapted to the evolution operators as follows
\begin{align*}
  \calu^2(0,T):= \{v\in \mathcal{D}'((0,T)\times \bbr^d):
                  t\mapsto U(0,t)v(t) \in U^2(0,T; L^2) \}
\end{align*}
endowed with the norm
$\|v\|_{\calu^2(0,T)} := \|U(0,\cdot) v\|_{U^2(0,T;L^2)}$,
where $U^2(0,T;L^2)$ is the space of type $U^2$
(see Section \ref{Sec-Pre} below for the precise definition).

The global well-posedness, uniform estimates and the temporal regularity
of controlled solutions to \eqref{equa-x} are summarized below.

\begin{theorem} \label{Thm-Equa-X}
Consider the mass-subcritical case $\lbb = \pm 1$, $1<\a<1+\frac 4d$,
or the defocusing mass-critical case  $\lbb =-1$, $\a= 1+\frac 4d$, $d\geq 1$.
Assume $(H0)$.
For each $X_0\in L^2$, $u\in \calu_{ad}$ (see \eqref{u-ad} below)
and $0<T<\9$,
there exists $\bbp$-a.s. a unique strong $L^2$-solution $X^u$ to  \eqref{equa-x},
satisfying that $|X(t)|_{L^2} = |X_0|_{L^2}$, $t\in [0,T]$,
and
\begin{align} \label{E-Xu-U2}
     \sup\limits_{u\in \calu_{ad}}\|e^{-W} X^u\|_{\calu^2(0,T) \cap L^2(0,T;H^{\frac 12}_{-1})} <\9, \ \ \bbp-a.s..
\end{align}
In particular,
\begin{align} \label{E-Xu-S0}
      \sup\limits_{u\in \calu_{ad}}\|X^u\|_{S^0(0,T) \cap L^2(0,T;H^{\frac 12}_{-1})}  < \9, \ \ \bbp-a.s.,
\end{align}
and we have $\bbp$-a.s. the temporal regularity
\begin{align} \label{E-Xu-timereg}
    \sup\limits_{u\in \calu_{ad}}
    \sup\limits_{h\geq 0}
    h^{-\frac 1 2}\bigg\|U(0,\cdot +h)e^{-W(\cdot+h)} X^u(\cdot + h) - U(0,\cdot)e^{-W}X^u\bigg\|_{L^2(0,T; L^2)} <\9.
\end{align}
Moreover,
assume additionally that $(H0)^*$ holds
in the defoucsing mass-critical case.
Then, for any $1\leq \rho <\9$,
\begin{align} \label{E-Xu-U2-Integ}
     \bbe \sup\limits_{u\in \calu_{ad}}\|e^{-W} X^u\|^\rho_{\calu^2(0,T) \cap L^2(0,T;H^{\frac 12}_{-1})} <\9.
\end{align}
\end{theorem}

\begin{remark}
The global well-posedness  in the mass-subcritical case
was proved in \cite{BRZ14,BRZ18}.
However,
the  energy method used there is not applicable
to the mass-critical case.
The proof of global well-posedness
in the critical case
in Theorem \ref{Thm-Equa-X} benefits from the idea of recent work \cite{Z18},
based on the work of Dodson \cite{D12,D16.1,D16.2},
new rescaling transformations
and stability results for nonlinear Schr\"odinger equations with lower order perturbations.
\end{remark}

\begin{remark} \label{Rem-U2-Sharper-Xu}
The estimate \eqref{E-Xu-U2} is sharper than \eqref{E-Xu-S0}
in the Strichartz space
(see also Corollary \ref{Cor-S0N0-U2} below)
and, more importantly,
reveals the new temporal regularity \eqref{E-Xu-timereg} of controlled solutions.
\end{remark}

\begin{remark}
Hypothesis $(H0)$ is sufficient to
yield the pathwise global well-posedness of \eqref{equa-x}
in both the subcritical and critical case,
and it also suffices to get the integrability of solutions
in the subcritical case.
We also expect it to imply the integrability in the critical case,
however,
we will not treat this technical problem in the present paper.
\end{remark}

Similarly to \eqref{E-Xu-timereg},
we also have the temporal regularity of $VX$ as specified in Theorem \ref{Thm-VX} below,
which is crucial to yield the
tightness in Section \ref{Sec-Proof} below
when the initial data are of the minimal $L^2$ regularity.

We shall use another new space
adapted to the evolution operators $\{U(t,s)\}$ as follows
\begin{align*}
  \calv^2(0,T):=\{v\in \mathcal{D}'(0,T):
                  t \mapsto U(0,t)v(t) \in V^2(0,T;L^2)\},
\end{align*}
equipped with the norm
$\|v\|_{\calv^2(0,T)} =\|U(0,\cdot) v\|_{V^2(0,T)}$,
where $V^2(0,T;L^2)$ is the space of type $V^2$
(see Section \ref{Sec-Pre} below).

\begin{theorem} \label{Thm-VX}
Consider the situations of Theorem \ref{Thm-Equa-X}.
Assume $(H0)$.
Assume additionally $(H0)^*$ in the defocusing mass-critical case.
We have for each $1\leq k\leq m$ and for any $1\leq \rho <\9$,
\begin{align} \label{E-VXu-V2}
     \bbe \sup\limits_{u\in \calu_{ad}} \|e^{-W} V_k X^u\|^\rho_{\calv^2(0,T)} <\9.
\end{align}
In particular,
we have the temporal regularity
\begin{align} \label{E-VXu-timereg}
    \bbe \sup\limits_{u\in \calu_{ad}}
    \sup\limits_{h\geq 0}
    h^{-\frac \rho 2}
    \bigg\|U(0,\cdot +h)e^{-W(\cdot+h)}V_k X^u(\cdot + h) - U(0,\cdot)e^{-W} V_k X^u\bigg\|^\rho_{L^2(0,T; L^2)} <\9.
\end{align}
\end{theorem}

Now, let us  introduce the optimal control problem  as in \cite{BRZ18}.

Let
$L^2_{ad}(0,T;\rr^m)$ denote the space of all
$(\calf_t)_{t\ge0}$-adapted $\rr^m$-valued processes
$u:[0,T]\to\rr^m$ such that $u\in L^2((0,T)\times\ooo;\rr^m)$.
Similarly,  let
$L^2_{ad}(0,T;L^2(\ooo;L^2))$ denote the space of
$L^2$-valued $(\calf_t)_{t\ge0}$-adapted processes $v$ such that  $\E\int^T_0|v(t)|^2_{L^2}dt<\9.$

Given
$\bbx_T\in L^2(\ooo,\calf_T, \mathbb{P} ;L^2)$
and $\bbx\in L^2_{ad}(0,T;L^2(\ooo;L^2))$,
we define the objective functional
$\Phi:L^2_{ad}(0,T;\rr^m)\to\rr$ by
\begin{align} \label{def-Phi}
    \Phi(u)
    :=&  \bbe |X^u(T)-\bbx_T|_{L^2}^2
      + \g_1 \bbe \int_0^T |X^u(t)-\bbx(t)|_{L^2}^2 dt \nonumber \\
      &   + \g_2 \bbe \int_0^T |u(t)|_m^2dt + \g_3 \bbe \int_0^T |u'(t)|_m^2dt.
\end{align}
Here, the coefficients $\g_1, \g_3 \geq 0, \g_2 >0$,
$|\cdot|_m$ denotes the Euclidean norm in $\bbr^m$,
and  $u'$ denotes the time derivative of $u$ if it exists.
In most situations, $\bbx$ is a given trajectory of the uncontrolled system
or,  in particular, a steady state solution.
The control $u$ belongs to the
admissible set $\calu_{ad}$ defined by
\begin{align}\label{u-ad}
\calu_{ad} :=&\bigg\{u\in L^2_{ad}(0,T;\rr^m);\ u\in K,\ \ a.e.\ on\ (0,T)\times\ooo. \bigg\},
\end{align}
where $K$ is a compact  convex subset of $\rr^m$.
Note that,
$\sup_{u\in \calu_{ad}}\|u\|_{L^\9(0,T; \bbr^m)} \leq D_K <\9$,
where $D_K$ denotes the diameter of $K$.

The optimal control problem considered in this paper is formulated below
\begin{enumerate}
\item[(P)]
{\it Minimize}
\begin{align*}
\Phi(u) =&  \bbe \bigg(|X(T)-\mathbb{X}_T|^2_{L^2}+\g_1 \int^T_0 |X(t)-\bbx(t)|^2_{L^2}dt \\
         & \qquad \qquad    + \g_2  \int^T_0 |u(t)|^2_m   dt + \g_3  \int^T_0 |u'(t)|^2_m   dt \bigg)
\end{align*}
on all $(X,u)\in L^2_{ad}(0,T; L^2(\Omega; L^2)) \times \calu_{ad}$ subject to \eqref{equa-x}.
\end{enumerate}
Heuristically,
the objective of  control process is to steer the quantum system from an initial state $X_0$
to a target state $\bbx_T$ and also in the neighborhood of a given trajectory $\bbx$.
The last two terms in the cost functional represent the energy cost to obtain the desired objective.

As mentioned in the Introduction,
because of the defect of compact imbedding in the stochastic situation,
the existence of a solution in Problem ${\rm (P)}$ does not follow from
standard compactness techniques used in deterministic optimization problems.
Below
we consider the relaxed version of Problem ${\rm (P)}$ as in \cite{BRZ18}
which, actually, resembles standard weak solutions to stochastic equations.

\begin{definition}\label{Def-Contr}
Let $\caly:= L^2(\bbr^d) \times L^2((0,T)\times \bbr^d) \times C([0,T]; \bbr^N) \times L^2(0,T; \bbr^m)$ $\times L^2( (0,T) \times \bbr^d)$
and $(\Omega^*, \mathcal{F}^*, (\mathcal{F}_t^*)_{t\geq 0})$ be a new filtered probability space,
carrying $(\bbx^*_T, \bbx^*, \beta^*, u^*, X^*)$ in $\mathcal{Y}$.
Define
$L^2_{ad^*}(0,T;\bbr^m)$,
$L^2_{ad^*}(0,T; L^2(\Omega; L^2))$, $\calu_{ad^*}$ and $\Phi^*(u^*)$  similarly as above on this new filtered probability space.

The system
$(\ooo^*,\calf^*,\mathbb{P}^*,(\calf^*_t)_{t\ge0}, \bbx_T^*, \bbx^*, \beta^*, u^*,X^*)$ is said to be {\it admissible}, if
$\bbx^*_T \in L^2(\Omega, \mathcal{F}^*_T, \bbp^*; L^2)$, $\bbx^* \in L^2_{ad^*}(0,T; L^2(\Omega; L^2))$,
$\beta^* = (\beta^*_1, \ldots, \beta^*_N)$ is an $(\calf^*_t)_{t\ge0}$-adapted $\bbr^N$-valued Wiener process,
the joint distributions of $(\bbx^*_T, \bbx^*_1, \beta^*)$ and $(\bbx_T, \bbx, \beta)$ coincide,
$u^*\in \mathcal{U}_{ad^*}$, and $X^*$ is an $L^2$-valued continuous $(\calf^*_t)_{t\ge0}$-adapted process
that satisfies equation \eqref{equa-x}
corresponding to $(\beta^*, u^*)$.

The admissible system $(\ooo^*,\calf^*,\mathbb{P}^*,(\calf^*_t)_{t\ge0},\bbx_T^*, \bbx^*,\beta^*, u^*,X^*)$
is said to be a {\it relaxed solution} to the optimal control problem (P),
if
\begin{equation}\label{Def-relaxsol}
\Phi^*(u^*) \leq  \inf\{\Phi(u);\ u\in\calu_{ad},\ X^u\ satisfies\ \eqref{equa-x}\}.
\end{equation}
\end{definition}

The first result concerning the optimal control problem of \eqref{equa-x}
in the case where $\g_3>0$
is formulated below.

\begin{theorem} \label{Thm-Control-deriv}
Consider the mass-subcritical
or the defocusing mass-critical case.
Assume $(H0)$.
Then, for each $X_0\in L^2$,
$\bbx_T\in L^2(\Omega, \mathcal{F}_T, \bbp; L^2)$
and
$\bbx \in L^2_{ad}(0,T; L^2(\Omega; L^2))$,
$0<T<\9$,
there exists a relaxed solution in the sense of Definition \ref{Def-Contr}
to the optimal control problem {\rm (P)} in the case $\g_3>0$.
\end{theorem}

In the case where $\g_3=0$,
we are able to obtain a relaxed solution satisfying \eqref{Def-relaxsol}
with ``='' replacing ``$\leq$''.
For this purpose,
we introduce
the dual linearized backward stochastic equation  below
\begin{align}\label{equa-back}
 &d Y= -i\Delta Y\,dt - \lbb i h_1(X^u)Y dt +\lbb i  h_2(X^u)  \ol{Y} dt + \mu Y dt - iV_0Y dt  - i u\cdot V Y dt\nonumber  \\
 &\qquad \quad        + \g_1 (X^u - \bbx) dt  - \sum\limits_{k=1}^N \ol{\mu_k} e_k Z_k dt +  \sum\limits_{k=1}^N Z_k d\beta_k(t),    \\
 & Y(T) = -(X^u(T)-\bbx_T),  \nonumber
\end{align}
where
\begin{equation}\label{h1-h2}
h_1(X^u):= \frac{\a+1}{2} |X^u|^{\a-1},\ \ h_2(X^u):=  \frac{\a-1}{2} |X^u|^{\a-3} (X^u)^2 .
\end{equation}
The functions $h_j$, $j=1,2$, are the usual complex derivatives of the function $z\mapsto |z|^{\a-1} z$,
i.e., $h_1(z) = \partial_z (|z|^{\a-1}z)$ and $h_2(z)= \partial_{\ol{z}} (|z|^{\a-1}z)$, $z\in \mathbb{C}$.

It should be mentioned that,
although equation \eqref{equa-back} is linear with respect to $Y$,
the presence of singular coefficient $h_2(X^u)$ destroys standard
Lipschitz property or monotonicity  required in the literature of stochastic backward equations
(see, e.g., \cite{FT02, HP91,T96,YZ99}).
Hence, the standard energy method,
based on the It\^o formula of $|Y(t)|_2^2$,
is not applicable to \eqref{equa-back}.
Moreover,
unlike the heat semigroup,
the free Schr\"odinger group $\{e^{-i t \Delta}\}$
is a unitary evolution in $L^2$ and so
has no global regularization effect to raise the Sobolev regularity of solutions,
which also makes the analysis of \eqref{equa-back} more difficult.
Actually,
a great effect of \cite{BRZ18} is dedicated to this issue.
The key idea is to use the duality arguments to
reduce the analysis of \eqref{equa-back} to that of
an associated variational equation (see \eqref{equa-var-Psi} below).
See also \cite{FO16,LZ14}
for  duality arguments in analyzing backward stochastic equations
arising form control problems.

The global well-posedness
and temporal regularity of
stochastic backward solutions are formulated below.

\begin{theorem} \label{Thm-BSPDE}
Consider the mass subcritical or the
defocusing mass-critical case.
Assume Hypothesis $(H0)$.
Assume additionally that $e_j$ are constants,
$1\leq j\leq N$.

Then,
given any $\bbx_T \in L^{2 + \nu}(\Omega, \mathcal{F}_T; L^2)$ and $\bbx\in L^{2+\nu}(\Omega; L^2(0,T; L^2))$
for some $\nu \in (0,1)$,
there exists a unique $(\mathcal{F}_t)$-adapted solution $(Y^u,Z^u)$
to \eqref{equa-back} corresponding to $u\in \calu_{ad}$,
satisfying that
$Y^u \in C([0,T]; L^2)$, a.s.,
for any $1\leq \rho <2+\nu$,
\begin{align} \label{E-Yu-V2}
     \sup\limits_{u\in \calu_{ad}} \bbe \|e^{-W}Y^u\|^\rho_{\calv^{2}(0,T)}  < \9,
\end{align}
and for each $1\leq k\leq N$,
\begin{align} \label{E-Zu-L2}
     \sup\limits_{u\in \calu_{ad}} \bbe  \|Z_k^u\|^\rho_{L^2(0,T; L^2)} < \9.
\end{align}
In particular, for any $1\leq \rho <2+\nu$,
\begin{align} \label{E-Yu-S0}
     \sup\limits_{u\in \calu_{ad}}  \bbe \|Y^u\|^\rho_{S^0(0,T)}  < \9,
\end{align}
and we  have the temporal regularity
\begin{align} \label{E-Yu-timereg}
   \sup\limits_{u\in \calu_{ad}} \bbe
    \sup\limits_{h\geq 0}
    h^{-\frac \rho 2}
    \bigg\|U(0,\cdot +h)e^{-W(\cdot+h)} Y^u(\cdot + h) - U(0,\cdot)e^{-W}Y^u \bigg\|^\rho_{L^2(0,T; L^2)} <\9.
\end{align}
\end{theorem}

\begin{remark} \label{Rem-V2-Sharper-Yu}
As in Remark  \ref{Rem-U2-Sharper-Xu},
the estimate \eqref{E-Yu-V2} is sharper than \eqref{E-Yu-S0}
in the Strichartz space
and yields the temporal regularity \eqref{E-Yu-timereg} of backward solutions.
The latter fact enables us to obtain the tightness of approximating controls
with initial data of minimal $L^2$ regularity.
See Section \ref{Sec-Proof} below.
\end{remark}

The main result of this paper is formulated below.
\begin{theorem}\label{Thm-Control}
Consider the mass-subcritical case
or the defocusing mass-critical case.
Assume Hypothesis $(H0)$.
Assume additionally that $e_j$ are constants, $1\leq j\leq N$.

Then, for each $X_0\in L^2$, $\bbx_T\in L^{2+\nu}(\Omega, \mathcal{F}_T;L^2) $ and
$ \bbx \in L^{2+\nu}(\Omega; L^2(0,T;$ $L^2)) $
for some $\nu\in (0,1)$, $0<T<\9$, there exists a relaxed solution
$(\Omega^*, \mathscr{F}^*, \bbp^*,$ $(\mathscr{F}^*_t)_{t\geq 0}, \bbx^*_T, \bbx^*, \beta^*, u^*, X^*)$
in the sense of Definition \ref{Def-Contr} to Problem {\rm (P)} in the case  $\g_3=0$,
satisfying
\begin{equation}\label{e2.10}
 \Phi^*(u^*)=\inf\{\Phi(u);\ u\in\calu_{ad},\ X^u\ satisfies\ \eqref{equa-x} \}.
\end{equation}
Moreover, we have the explicit characterization of $u^*$ below
\begin{equation}\label{e5.1}
u^* (t)=P_K \(\frac 1{\g_2}\,{\rm Im}\int_{\rr^d}V(x) X^*(t,x) \ol{Y^*}(t,x)dx\),\ \ff t\in[0,T],\ \bbp^*-a.s.
\end{equation}
where $P_K$ is the projection on $K$,
and
$Y^*$ is the first component of backward solution $(Y^*,Z^*)$ to \eqref{equa-back}
with
$\bbx^*_T, \bbx^*, \beta^*, u^*, X^*$ replacing
$\bbx_T, \bbx, \beta, u, X^u$, respectively.
\end{theorem}

\begin{remark}
Theorem \ref{Thm-Control} is applicable to the whole mass-subcritical case
where $1<\a <1+\frac 4d$, $d\geq 1$,
which includes the ranges $\a\in (1,2)$ with $1\leq d\leq 3$
and $\a\in (1,1+\frac d4)$ with $d\geq 4$
that were previously excluded in \cite{BRZ18}.
More importantly,
Theorem \ref{Thm-Control} applies also to the defocusing mass-critical case,
which includes the cubic and quartic stochastic nonlinear Schr\"odinger equations
in dimension $d=2$ and $d=1$, respectively.
\end{remark}

\begin{remark}
The additional condition that $e_j$ are constants
is imposed in Theorem \ref{Thm-Control} above
mainly for the integrability of variational solutions
(see Proposition \ref{Prop-Var-Psi} below).
This condition can be reduced to that
$\bbe \|X^u\|^\rho_{S^0(0,T)} + \bbe \|\psi^u\|_{C([0,T];L^2)}^\rho <\9$, $\forall 1\leq \rho<\9$,
where $X^u$ and $\psi^u$ solve \eqref{equa-x} and \eqref{equa-var-Psi} below respectively.
See Remark \ref{Rem-var-Integ} below for detail explanations.
As a matter of fact,
Hypothesis $(H0)$ suffices to give the pathwise gloabl well-posedness
and estimates of solutions to both \eqref{equa-x} and \eqref{equa-var-Psi},
and so it is also expected  to yield the integrability of solutions.
\end{remark}

\begin{remark}
It is also interesting to study the case where the controls
also enter into the diffusion coefficients.
In that case,
one usually needs to deal with a second adjoint equation
(see e.g. \cite{FO16,LZ14,P90} and references therein).
In the present work,
we will not treat this issue.
\end{remark}

As a byproduct,
we have the following result in the deterministic case.

\begin{corollary} \label{Thm-deter}
In the deterministic case (i.e. $\mu_k=0$, $1\leq k\leq N$),
consider
the mass-subcritical case
or the defocusing mass-critical case as in Theorem \ref{Thm-Control}.

Then, for each $X_0\in L^2$, $\bbx_T\in L^2 $ and
$ \bbx \in  L^2(0,T; L^2) $,  $0<T<\9$, there exists an optimal control $u$ to Problem ${\rm (P)}$
in the case $\g_3=0$, such that
\begin{equation}\label{e2.10-deter}
\Phi(u)=\inf\{\Phi(v);\ v\in\calu_{ad},\ X^v\ satisfies\ \eqref{e2.1}\}.
\end{equation}
Moreover,
\begin{equation}\label{e5.1-deter}
u(t)=P_K\(\frac 1{\g_2}\,{\rm Im}\int_{\rr^d}V(x) X(t,x) \ol{Y}(t,x)dx\),\ \ff t\in[0,T],
\end{equation}
where $P_K$ is the projection on $K$, and $Y$ is the solution to backward
equation \eqref{equa-back} with $Z_k=0$,
$1\leq k\leq N$.
\end{corollary}

The remaining part of this paper is organized as follows.
In Section \ref{Sec-Pre}
we present  preliminaries used in this paper,
including the Strichartz and local smoothing estimates,
$U^p$-$V^p$ spaces
and stability results for nonlinear Schr\"odinger equations
with lower order perturbations.
Sections \ref{Sec-Control-equa} and \ref{Sec-Backequa}
constitute the technical parts of this paper.
In Section \ref{Sec-Control-equa},
we treat the controlled stochastic equation \eqref{equa-x}.
We first prove Theorems \ref{Thm-Equa-X} and \ref{Thm-VX} in Subsection \ref{Subsec-GWP-ContrEqua},
and then in Subsection \ref{Subsec-Sta-ContrEqua}
we prove stability results for controlled solutions and
objective functionals,
which consequently lead to the proof of Theorem \ref{Thm-Control-deriv}.
In Section \ref{Sec-Backequa},
we mainly analyze the variational equations
and the dual backward stochastic equations.
As a consequence,
we obtain the directional derivative of objective functional.
Section \ref{Sec-Proof} is mainly devoted to the proof of main result  Theorem \ref{Thm-Control}.
For simplicity of the exposition, some technical proofs are postponed to the Appendix.

\section{Preliminaries} \label{Sec-Pre}

\subsection{Strichartz and local smoothing estimates}

We start with the Strichartz and local smoothing estimates
related to  operators of the form $e^{-\Phi}\Delta(e^{\Phi} \cdot)$.
In the sequel, we assume that $\Re \Phi =0$.

\begin{theorem}   \label{Thm-Stri}
Let $ I=[t_0,T]\subseteq \bbr^+$.
Consider the equation
\begin{align} \label{equa-stri}
   & i\partial_t v = e^{-\Phi} \Delta(e^\Phi v) + f.
\end{align}
Here,
the function
$\Phi =\Phi(t,x)$ is continuous on $t$ for each $x\in \bbr^d$, $d\geq 1$,
$\Re \Phi =0$,
and
for each multi-index $\g$,
\begin{align} \label{psi-Stri}
     \sup\limits_{t\in I} |\p_x^\g \Phi(t,x)| \leq C(\g) \sup\limits_{t\in I}g(t)  \<x\>^{-2}
\end{align}
for some positive and continuous function $g$.
Then, for any  $v(t_0)\in L^2$
and $f\in N^0(I) + L^2{(I; H^{-\frac 12}_1)}$,
the solution $v$ to \eqref{equa-stri} satisfies
\begin{align}  \label{L2-Stri}
     \|v\|_{S^0(I) \cap L^2(I;H^\frac 12_{-1})}
    \leq& C_T (|v({t_0})|_{L^2}+
           \|f\|_{N^0(I)+L^2(I;H^{-\frac 12}_{1})})
\end{align}
Moreover,
if $\bbe \sup_{0\leq t\leq T} (g(t))^\rho \leq C(T, \rho)$, $\forall 1\leq \rho <\9$.
then so is $C_T$, i.e.,
\begin{align} \label{E-CT}
    \bbe (C_T)^\rho \leq C(\rho, T)<\9,\ \ \forall 1\leq \rho<\9.
\end{align}
\end{theorem}

{\bf Proof.} The proof of \eqref{L2-Stri}
is quite similar to that of \cite[Theorem 3.3]{Z18},
based on the pseudo-differential calculus.
Actually, the asymptotical flatness condition \eqref{psi-Stri}
ensures that the
lower order perturbations arising in the operator $e^{-\Phi}\Delta( e^{\Phi} \cdot)$
can be controlled, via the G{\aa}rding inequality,
by the Poisson bracket $i[\Psi_h,\Delta]$
for some appropriate symbol $h\in S^0$.
Concerning the constant $C_T$,
similar arguments as in the proof of \cite[Theorem 2.6]{Z17}
show that $C_T$
depends polynomially on $\sup_{0\leq t\leq T} g(t)$,
which implies \eqref{E-CT}.
\hfill $\square$

Let $\{V(t,s)\}$ be the evolution operators related to $-ie^{-\Phi}\Delta(e^\Phi \cdot)$,
i.e.,
for any $v\in L^2$,
$v(t):=V(t,s)v$ satisfies the equation
\begin{align*}
  & i\p_t v(t) = e^{-\Phi}\Delta(e^\Phi v(t)), \ \ t\in \bbr^+,
\end{align*}
with $v(s) =v$.
Denote by $V^*(t,s)$ the dual operator of $V(t,s)$, $t,s\in \bbr^+$.
We have the following result concerning the evolution operators.
\begin{proposition} \label{Prop-V-l2}
Let $\Phi$ be as in Theorem \ref{Thm-Stri}.
Then, for any $v\in L^2$, $t\in \bbr^+$,
\begin{align}
   & V^*(t,0) v= V(0,t) v,   \label{V*-V} \\
   & V^*(0,t) v = V(t,0) v.  \label{V*-V.2}
\end{align}
In particular, for any $v_1, v_2\in L^2$, $t\in \bbr^+$,
\begin{align} \label{V*-V-L2}
   \<v_1, v_2\>_2 = \<V(0,t)e^{-\Phi(t)} v_1, V(0,t)e^{-\Phi(t)} v_2\>_2.
\end{align}
\end{proposition}

The proof is postponed to the Appendix for simplicity.

We end this subsection with Theorem \ref{Thm-Rescale-sigma} below,
which relates the stochastic equation \eqref{equa-x} to
a random Schr\"odinger equation with lower order perturbations.
Note that,
the local well-posedness of \eqref{equa-x}
in the mass-critical case can be proved by using similar arguments as in \cite{BRZ14}.

\begin{theorem} \label{Thm-Rescale-sigma}
Consider the defocusing mass-critical case,
i.e., $\lbb =-1$, $\a=1+\frac 4d$, $d\geq 1$.
Let $X$ be the $L^2$-solution to \eqref{equa-x} on
$[0,\tau^*)$
with $X(0)=X_0 \in L^2$,
where $\tau^*$ is the maximal existing time.
Given any $(\mathscr{F}_t)$-stopping time $\sigma$
satisfying  $0\leq \sigma <\tau^*$,
we define the rescaling transformation
\begin{align}  \label{vsigma}
     v_\sigma(t) := e^{-W_\sigma(t)} X(\sigma+t),\ \ t\in[0,\tau],
\end{align}
where
$W_\sigma(t)
   := W(\sigma+t) - W(\sigma)$.
Set $f(u_\sigma) := V_0 + u_\sigma \cdot V$
with $u_\sigma(t) = u(\sigma+t)$.
Then,
$v_\sigma$ satisfies $\bbp$-a.s. the random equation
\begin{align} \label{equa-RNLS}
   i\p_t v_\sigma
  =& e^{-W_\sigma} \Delta (e^{W_\sigma} v_\sigma) -   F(v_\sigma) + f(u_\sigma)v_\sigma, \\
  v_\sigma(0)=& X(\sigma), \nonumber
\end{align}
on $[0,\tau^*-\sigma)$ in the space $H^{-2}$.
\end{theorem}

The proof is similar to that of \cite[Thoerem 2.17]{Z18}
and is omitted here.
Note that,
since
$\Re \mu_j=0$, $1\leq j\leq N$,
the term $\wh{\mu}$ defined in \cite[(2.22)]{Z18}
equals to zero and so does not appear in \eqref{equa-RNLS}.

\subsection{$U^p$-$V^p$ spaces}

In this subsection,
we first recall some basic facts of the spaces $U^p$ and $V^p$,
following \cite{BOP15,HHK09,HTT11,KTV14}.
Then, we enhance the Strichartz estimates in Theorem \ref{Thm-Stri}
to the spaces of type $U^2$ and $V^2$,
adapted to the evolution operators $\{U(t,s)\}$
as in Theorem \ref{Thm-Equa-X},
which will play an important role in the proof of tightness
in Section \ref{Sec-Proof} below.

Let $H$ be a separable Hilbert space over $\bbc$,
it will be chosen as either $L^2(\bbr^d)$ or $\bbc$.
Let $\calz$ be the set of finite partitions $\{t_k\}_{k=0}^K$ of $\bbr$.
If $t_K=\9$, we use the convention that $u(t_K)=0$
for all functions $u:\bbr \mapsto H$.
We also use the notation $\calx_I$ for the  characteristic function of a set $I \subseteq \bbr$.

\begin{definition}
Let $1\leq p<\9$.
\begin{enumerate}
  \item[(i)] A $U^p$-atom is a right continuous step function $u: \bbr \mapsto H$ of the form
  \begin{align*}
       a = \sum\limits_{k=1}^K \calx_{[t_{k-1}, t_k)}\phi_{k-1},
  \end{align*}
  for some $\{t_k\}_{k=0}^K \in \calz$ and $\{\phi_k\}_{k=0}^K \subseteq H$
  with $\sum_{k=0}^{K-1}\|\phi_k\|_{H}^p =1$.
  The $U^p$-space $U^p(:=U^p(\bbr; H))$ is defined by
  \begin{align*}
      U^p:=\{u=\sum\limits_{j=1}^\9 \lbb_j a_j; a_j\ are\ U^p-atoms\ and\ \{\lbb_j\}_{j\in \bbn} \in l^1(\bbn; \bbc)\}
  \end{align*}
  with the norm
  \begin{align*}
      \|u\|_{U^p}:= \inf\{\sum\limits_{j=1}^\9 |\lbb_j|; u=\sum\limits_{j=1}^\9 \lbb_j a_j, \lbb_j\in \bbc, a_j\ are\ U^p-atoms\}.
  \end{align*}
  Likewise, $U^p_c$ denotes the closed subspace in $U^p$ of all continuous functions $u:\bbr \mapsto H$.

  \item[(ii)]
  The $V^p$-space $V^p(:=V^p(\bbr; H))$ is the space of all functions of bounded $p$-variation
  $v: \bbr \mapsto H$ satisfying $v(\9)=0$,
  endowed with the norm
  \begin{align*}
     \|v\|_{V^p(\bbr; H)}
     := \sup\limits_{\{t_k\}_{k=0}^K \in \calz} (\sum\limits_{k=1}^K\|v(t_k)-v(t_{k-1})\|_{H}^p)^{\frac 1p}.
  \end{align*}
  Likewise, let $V^p_{rc}$  denote the closed subspace of all right continuous functions
  $v: \bbr \mapsto H$ such that $\lim_{t\to -\9} v(t)=0$.
\end{enumerate}
\end{definition}

\begin{lemma} (\cite{HHK09}) \label{Lem-Embed-UpVp}
The spaces $U^p$, $U^p_c$, $V^p$, $V^p_{rc}$ are Banach spaces.
For any $1\leq p<q<\9$, we have the embeddings
\begin{align} \label{embed-Up-Vp}
    U^p \hookrightarrow V_{rc}^p \hookrightarrow U^q  \hookrightarrow L^\9 .
\end{align}
Moreover, for $1<p<\9$,
$V^p \hookrightarrow \dot{B}^{\frac 1p}_{p,\9}$, i.e.,
for $v\in V^p$,
\begin{align} \label{embed-Vp-Bp}
    \sup\limits_{h>0} h^{-\frac 1p} \|v(\cdot+h)- v \|_{L^p(\bbr;H)} \leq C \|v\|_{V^p(\bbr;H)}
\end{align}
\end{lemma}

The duality between  spaces $U^p$ and $V^p$ are related to a bilinear form
$B(u,v): U^p \times V^{p'} \mapsto B(u,v)$
defined as follows:
for $u\in U^p$, $v\in V^{p'}$
and any $\mathfrak{t}=\{t_k\}_{k=0}^K \in \calz$,
$B_{\mathfrak{t}}(u,v) := \sum_{k=1}^K\<u(t_{k-1}), v(t_k)-v(t_{k-1})\>_2$.
Then, there exists a unique number $B(u,v)$ such that
for all $\ve >0$, there exists $\mathfrak{t} \in \calz$
such that for any $\mathfrak{t}' \supseteq \mathfrak{t}$,
$|B_{\mathfrak{t}'}(u,v) - B_{\mathfrak{t}}(u,v)| <\ve$.
(See \cite[Proposition 2.7]{HHK09}.)

\begin{lemma} (\cite[Propositions 2.7, 2.10]{HHK09}) \label{Lem-Buv}
Let $1\leq p<\9$. For $u\in U^p$, $v\in V^{p'}$,
$|B(u,v)| \leq \|u\|_{U^p} \|v\|_{V^{p'}}$.
Moreover, if in addition $u \in V^1$ and $u$ is absolutely continuous on compact intervals,
then,
\begin{align*}
   B(u,v) = - \int\limits_{\bbr} \<u'(t), v(t)\>_2 dt.
\end{align*}
\end{lemma}

Given any time interval $I= [a,b) \subseteq \bbr$,
we also consider the restriction space $U^p(I)$ with the norm
$\|u\|_{U^p(I)} := \inf \{\|v\|_{U^p(\bbr)}; v|_I =u \}$.
Actually, the infimum is attained by $v= \calx_{I} u$
(See \cite[Remark A.2]{BOP15}).

\begin{lemma} (\cite[Lemmas A.3, A.6]{BOP15}) \label{Lem-UI}
Given any interval $I= [a,b) = \cup_{j=1}^\9 I_j$, we have
$\|u\|_{U^p(I)} \leq \sum_{j=1}^\9 \|u\|_{U^p(I_j)}$.
Moreover, for $u\in U_c^p(I)$,
the mapping $t \mapsto \|u\|_{U^p([a,t))}$ is continuous.
\end{lemma}

We have the following result from  \cite[Corollary 2.4]{KTV14}
for $v\in C([0,T];H)$.

\begin{lemma} (\cite[Corollary 2.4]{KTV14}) \label{Lem-Dual-UpVp}
Let $1 < p,p'<\9$, $\frac 1 p + \frac{1}{p'} =1$,
and $v\in V^p(0,T; H) \cap C([0,T]; H)$.
We have
\begin{align*}
    \|v\|_{V^p(0,T; H)} = \sup\{B(u,v): u\in C_0^\9(0,T; H), \|u\|_{U^{p'}(0,T;H)} =1\}.
\end{align*}
\end{lemma}

Similarly to \eqref{embed-Vp-Bp},
we have
\begin{lemma} \label{Lem-Embed-UpVp*}
Let $v\in V^p(0,T;H)\cap C([0,T];H)$, $1\leq p<\9$,
$0<T<\9$. Then,
\begin{align} \label{embed-Vp-Bp*}
    \|v(\cdot +h)-v\|_{L^p(0,T-h;H)}
    \leq 2^{1+\frac 1p} h^{\frac 1p} \|v\|_{V^p(0,T)}.
\end{align}
\end{lemma}
(See the Appendix for the proof.)

The following estimates relates the $U^p-V^p$ spaces and Strichartz spaces.

\begin{proposition}  \label{Pro-LpLq-UqVq}
Consider the situations in Theorem \ref{Thm-Stri}
and let $V(t,s)$, $t,s\in \bbr^+$, be the evolution operators
related to $-ie^{-\Phi}\Delta(e^{\Phi}\cdot)$.
Then, for any $I=[t_0, T) \subseteq \bbr^+$
and any Strichartz pair $(p,q)$,
\begin{align}
    & \|u\|_{L^q(I; L^p)} \leq C(T) \|V(t_0,\cdot)u\|_{U^q(I)}, \label{Hom-LpLq-Uq} \\
    & \bigg\|\int_{t_0}^\cdot V(t_0,s) f(s) ds \bigg\|_{V^{q'}(I)}
     \leq C(T) \|f\|_{L^{q'}(I; L^{p'})}, \label{Inhom-Vq-LpLq}
\end{align}
and
\begin{align}     \label{Inhom-V2-LpLq-LS}
   \bigg\|\int_{t_0}^\cdot V(t_0, s) f(s) ds \bigg\|_{V^2(I)}
   \leq C(T)(1+ \|f\|_{L^{q'}(I;L^{p'}) + L^2(I; H^{-\frac 12}_1)}).
\end{align}
where $C(T)$ is independent of $p,q$ and $C(T) \in L^\rho (\Omega)$ for any $1\leq \rho <\9$.
\end{proposition}

We postpone the proof to the Appendix.
As a consequence, we have

\begin{corollary} \label{Cor-S0N0-U2}
Consider the situations in Theorem \ref{Thm-Stri}.
Then,
\begin{align}
   \|u\|_{S^0(I)} \leq C(T) \|V(t_0, \cdot)u\|_{V^2(I)}
                  \leq C(T) \|V(t_0,\cdot)u\|_{U^2(I)}, \label{Hom-S0-U2}
\end{align}
and
\begin{align}
   \bigg\|\int_{t_0}^\cdot V(t_0,s)f(s) ds \bigg\|_{U^2(I)} \leq C(T) \|f\|_{N^0(I)},   \label{Inhom-U2-N0}
\end{align}
where $C(T) \in L^\rho(\Omega)$, $\forall 1\leq \rho<\9$.
\end{corollary}

{\bf Proof.}
The proof follows from Proposition \ref{Pro-LpLq-UqVq} and
the embedding  $V^{q'} \hookrightarrow U^2(I) \hookrightarrow V^2(I) \hookrightarrow U^q(I)$
for $2<q<\9$.
\hfill $\square$

\subsection{Mass-critical stability results}

We first present the recent the global well-posedness results obtained by Dodson \cite{D12,D16.1,D16.2}
in the defocusing mass-critical case.

\begin{theorem} (\cite{D12,D16.1, D16.2}) \label{Thm-L2GWP-Det}
For any $v_0\in L^2$,
there exists a unique global $L^2$-solution $v$ to the equation
\begin{align} \label{equa-u-L2}
   i\p_t v =&   \Delta v -  |v|^{\frac{4}{d}} v,  \\
   v(0)=& v_0  \nonumber
\end{align}
with $d\geq 1$. Moreover,
\begin{align} \label{globdd-u-L2-Det}
   \|v\|_{L^{2+\frac 4d}(\bbr\times\bbr^d) \cap L^2(\bbr; H^\frac 12_{-1})} \leq B_0(|v_0|_{L^2}) <\9,
\end{align}
where $B_0(|v_0|_{L^2})$ depends continuously on $v_0$ in $L^2$.
\end{theorem}

\begin{remark}
The proof of global well-posedness (and also scattering)
in \cite{D12,D16.1, D16.2} is highly nontrivial and is
based on the concentration-compact arguments.
One key ingredient there is the long-time Strichartz estimate.
The bound of $\|v\|_{L^2(\bbr; H^{\frac 12}_{-1})}$ in \eqref{globdd-u-L2-Det}
follows from Strichartz estimates
and that of $\|v\|_{L^{2+\frac 4d}(\bbr\times\bbr^d)}$,
and the continuity of $B_0(\cdot)$
follows from the mass-critical stability result Lemma 3.6 of \cite{TVZ07}.
We also refer   to \cite{CKSTT08,RV07,V07} for the defocusing energy-critical case.
\end{remark}

Below we state the mass-critical stability result
for Schr\"odinger equations with lower order perturbations
which is crucial in the proof of Theorem \ref{Thm-Equa-X}.

\begin{theorem}  \label{Thm-Sta-L2} ({\it Mass-Critical Stability Result}).
Fix $I=[t_0,T]\subseteq \bbr^+$.
Take the Strichartz pair $(p,q)$,
$p=2+\frac 4d$ with $d\geq 1$,
or $\frac 1p \in (\max\{\frac{1}{2\a}, \frac 12 - \frac{1}{2d}\}, \frac{1}{\a}(\frac 12 + \frac 1d))$
with $\a=1+\frac 4d$ and $1\leq d\leq 3$.
Let $V(t,s)$, $t,s \in I$, be the evolution operators related to
the operator $-i e^{-\Phi} \Delta (e^{\Phi} \cdot)$,
where $\Phi$ is as in Theorem \ref{Thm-Stri}
with the time function $g$.
Let $v$ be the solution to
\begin{align} \label{equa-v-p}
     v(t) = V(t,t_0) v(t_0) + \int_{t_0}^t V(t,s) (i |v|^{\frac 4d} v
               + G v) ds + R(t)
\end{align}
for some functions $R$ with $R(t_0) =0$,
$v(t_0)\in L^2$,
and $G\in L^\9(I\times \bbr^d)$.
Let $\wt{v}$ solve the equation
\begin{align} \label{equa-wtv-p}
    \wt{v}(t) = V(t,t_0) \wt{v}(t_0) + \int_{t_0}^t V(t,s) (i |\wt{v}|^{\frac 4d} \wt{v}
               + G \wt{v}+  e) ds
\end{align}
for some functions $e$ and $\wt{v}(t_0) \in L^2$.
Assume that
$ \|\wt{v}\|_{L^{q}(I; L^p)} \leq L $
for some positive constant $L$.
Assume also the smallness conditions
\begin{align} \label{Sta-L2-ve.1}
    \| V(\cdot, t_0)(v(t_0)-\wt{v}(t_0)) \|_{L^{q}(I; L^p)} \leq \ve,
\end{align}
\begin{align} \label{Sta-L2-ve.2}
    \|R\|_{L^{q}(I; L^p) \cap L^1(I; L^2)} \leq \ve,\
     \| e \|_{N^0+L^2(I; H^{-\frac 12}_1)}\leq \ve
\end{align}
for some $0<\ve\leq \ve_*$,
$\ve_* = \ve_*(C_T,  L)>0$ is a small constant,
and $C_T$ is the Strichartz constant in Theorem \ref{Thm-Stri}.
Then,
\begin{align}
   & \|v-\wt{v} -R\|_{L^{q}(I; L^p) \cap C(I; L^2)} \leq C_*(C_T, L) \ve , \label{Sta-L2.1}\\
   & \|v\|_{L^{q}(I; L^p)} \leq C_*(C_T, L).  \label{Sta-L2.3}
\end{align}
where $(\ve_*(C_T, L))^{-1}$, $C_*(C_T, L)$
can be taken to be nondecreasing with respect to arguments.
\end{theorem}

The proof of Theorem \ref{Thm-Sta-L2} is quite similar to that of
\cite[Theorem 4.1]{Z18}
and \cite[Proposition 4.6]{FX18.1}.
Actually, the linear terms $Gv$, $G\wt{v}$
and the error terms $R$ and $e$
act as perturbations of the nonlinearity.
So, the two solutions stay close to each other
as long as these perturbations are small enough.
For simplicity, the proof is postponed to the Appendix.

\section{The controlled equations} \label{Sec-Control-equa}

\subsection{Global well-posedness and integrability} \label{Subsec-GWP-ContrEqua}

In this subsection,
we  prove Theorems \ref{Thm-Equa-X} and \ref{Thm-VX}.
Let us first treat the subcritical case in Theorem \ref{Thm-Equa-X}.

{\bf Proof of Theorem \ref{Thm-Equa-X}.}
({\bf Mass-subcritical case}.)
Using $v:= e^{-W} X$ we have
\begin{align} \label{equa-v-contr}
   i\p_t v = e^{-W}\Delta(e^W v) + \lbb F(v)+ f(u)v,
\end{align}
where
$F(v) := |v|^{\a-1}v$ and
$f(u):= V_0 + u\cdot V$.
It has been proved in \cite[(3.19)]{BRZ18}
that for the Strichartz pair $(p,q) = (\a+1, \frac{4(\a+1)}{d(\a-1)})$,
$d\geq 1$,
\begin{align} \label{esti-xu-sub}
   \sup\limits_{u\in \calu_{ad}}
   \|X^u\|_{L^q(0,T; L^p)}
   \leq ([\frac T t]+1)^{\frac 1q} \frac{\a}{\a-1} C_T (1+D_*T) |X_0|_2,
\end{align}
where $\theta = 1-\frac{d(\a-1)}{4}>0$,
$t= \a^{-\frac{\a}{\theta}}(\a-1)^{\frac{\a-1}{\theta}}(|X_0|_2+1)^{-\frac{\a-1}{\theta}}
 C_T^{-\frac{\a}{\theta}} (1+ D_*T)^{-\frac{\a-1}{\theta}}$,
$D_* =  |V_0|_{L^\9}+D_K \|V\|_{L^\9(\bbr^d; \bbr^m)} $,
and $C_T$ is the constant  in Theorem \ref{Thm-Stri}.

Then, using \eqref{E-CT} we infer that
\begin{align} \label{bdd-Xu-LpLq-sub}
  \bbe \sup\limits_{u\in \calu_{ad}}
       \|X^u\|^\rho_{L^q(0,T; L^p)}
  \leq C(\rho, T)<\9,\ \ \forall 1\leq \rho <\9,
\end{align}
which along with Strichartz estimates \eqref{L2-Stri}
yields \eqref{E-Xu-S0}.

In order to obtain \eqref{E-Xu-U2},
we reformulate \eqref{equa-v-contr} as follows
\begin{align*}
   v(t) = U(t,0)X_0
        + (-i) \int_0^t U(t,s) (\lbb F(v) + f(u) v) ds,
\end{align*}
where $U(t,s)$, $t,s\geq 0$, are the evolution operators related to $-ie^{-W}\Delta(e^W\cdot)$.
Since $U(0,t)U(t,s) = U(0,s)$ for all $0\leq s\leq t$,
we have
\begin{align*}
   U(0,t)v(t) = X_0
        +(-i)\int_0^t U(0,s) (\lbb F(v) + f(u) v) ds.
\end{align*}
Then, by Corollary \ref{Cor-S0N0-U2},
H\"older's inequality and
$\|v\|_{C([0,T]; L^2)} = |X_0|_{L^2}$,
\begin{align*}
  \|U(0,\cdot) v\|_{U^2}
  \leq& |X_0|_{L^2}
       + C_T(\|F(v)\|_{L^{q'}(0,T;L^{p'})}+ \|f(u)v\|_{L^1(0,T; L^2)}) \\
  \leq& |X_0|_{L^2}
       + C_T(T^\theta\|X^u\|^\a_{L^{q}(0,T;L^{p})}+D_*T|X_0|_{L^2}).
\end{align*}
This along with \eqref{bdd-Xu-LpLq-sub} yields
\eqref{E-Xu-U2-Integ},
so \eqref{E-Xu-U2} follows.
Using \eqref{embed-Vp-Bp*} we get \eqref{E-Xu-timereg}.

Therefore,
the proof in the mass-subcritical case is complete.
\hfill $\square$

Next we consider \eqref{equa-x} in the defocusing mass-critical case.
We shall first prove the global well-posedness of \eqref{equa-x}
by using the idea of \cite{Z18}, based on stability result Theorem \ref{Thm-Sta-L2}
and a series of rescaling transformations.
Then,
we adapt the arguments in \cite{FX18.1}
to get the integrability of global solutions.

{\bf Proof of Theorem \ref{Thm-Equa-X}.} ({\bf Mass-critical case}.)
First,
by using Theorem \ref{Thm-Stri}
and  similar arguments as in \cite{BRZ14},
we have a unique $L^2$-solution $X^u$ to \eqref{equa-x}
on a maximal existing time interval $[0,\tau^*)$,
where $\tau^* (\leq T)$ is an $(\mathscr{F}_t)$ stopping time.
Moreover, $\tau^* = T$ a.s. if
\begin{align} \label{bdd-Xu-V}
     \|X^u\|_{L^{2+\frac 4d}((0,\tau^*)\times \bbr^d)} <\9,\ \ a.s..
\end{align}
Below we prove the global existence and $L^\rho(\Omega)$-integrability of $X^u$ separately.

{\bf Global existence.}
We shall prove the global estimate \eqref{bdd-Xu-V} above.
For this purpose,
we reformulate \eqref{equa-u-L2} as follows
\begin{align}
   &i\p_t \wt{v}
   = e^{-W}\Delta(e^W  \wt{v}) + F(\wt{v}) + f(u)\wt{v}  + e, \\
   & \wt{v}(0) = X_0, \nonumber
\end{align}
with the error term
\begin{align*}
   e= -(b(t)\cdot \na + c(t)) \wt{v} - f(u) \wt{v},
\end{align*}
where $b(t)$ and $c(t)$ are as in \eqref{equa-v-homo}.

In order to estimate the error term,
we see that
\begin{align*}
   \<x\> \<\na\>^{-\frac 12} (b(t)\cdot \na + c(t))
   = \Psi_{p(t)} \<x\>^{-1} \<\na\>^{\frac 12},
\end{align*}
where $\Psi_{p(t)}:= \<x\> \<\na\>^{-\frac 12} (b(t)\cdot \na + c(t))  \<\na\>^{-\frac 12} \<x\>$
is a pseudo-differential operator of order $0$
satisfying that for any $l\geq 1$,
the semi-norms $|p(t)|_{S^0}^{(l)} \leq C(l) \sup_{0\leq s\leq t}|\beta(s)|$,
where $|\beta(s)| =\max_{1\leq j\leq N} |\beta_j(s)|$.
(Here we may take $t$ small enough such that
$|\beta_j(t)|\leq 1$
and so $|\beta_j(t)|^2 \leq |\beta_j(t)|$.)

Then, using the $L^2$-boundedness of $\Psi_{p(t)}$ (see \cite[Lemma 2.4]{Z17})
and Theorem \ref{Thm-L2GWP-Det} we get
\begin{align} \label{esti-wtv-LS}
   \|(b(t)\cdot \na + c(t)) \wt{v}\|_{L^2(0,t; H^{-\frac 12}_1)}
   \leq& C_1 \|\wt{v}\|_{L^2(0,t; H^{\frac 12}_{-1})}  \sup\limits_{0\leq s\leq t} |\beta(s)|  \nonumber \\
   \leq& C_1 B_0(|X_0|_{L^2}) \sup\limits_{0\leq s\leq t} |\beta(s)| .
\end{align}

Moreover, since $|\wt{v}(t)|_{L^2} = |\wt{v}(0)|_{L^2} = |X_0|_{L^2}$,
we have,
if $D_*$ is as in \eqref{esti-xu-sub},
\begin{align} \label{esti-wtu-uV}
   \| f(u) \wt{v}\|_{L^1(0,t; L^2)}
   \leq D_*  |X_0|_{L^2} t.
\end{align}

Hence, letting $h(t):= |\beta(t)|+t$, $D_0 := (C_1+D_*)(B_0(|X_0|_{L^2})+|X_0|_{L^2})$,
we obtain
\begin{align} \label{esti-e-Xu}
   \|e\|_{N^0(0,t) + L^2(0,t; H^{-\frac 12}_1)}
   \leq& \|(b(t)\cdot \na + c(t)) \wt{v}\|_{L^2(0,t; H^{-\frac 12}_1)}
        + \| f(u)\wt{v}\|_{L^1(0,t; L^2)}  \nonumber \\
   \leq& D_0  \sup\limits_{s\in [0,t]} h(s).
\end{align}

Thus,
set $\tau_1:= \inf\{0\leq t\leq \tau^*: D_0 \sup_{s\in [0,t]}h(s) \geq \ve_1(t)\} \wedge \tau^*$,
where the small constant $\ve_1(t):= \ve_*(C_t, B_0(|X_0|_{L^2}))$
is as in Theorem \ref{Thm-Sta-L2}
with  $B_0(|X_0|_{L^2})$
replacing $L$.
We deduce from \eqref{esti-e-Xu} that
\begin{align*}
   \|e\|_{N^0(0,\tau_1) + L^2(0,\tau_1; H^{-\frac 12}_1)}
   \leq \ve_*(C_t, B_0(|X_0|_{L^2})).
\end{align*}
This, via Theorem \ref{Thm-Sta-L2} with $G=-if(u)$ and $R=0$, yields that
\begin{align*}
   \|X^u\|_{L^{2+\frac 4d}((0,\tau_1)\times \bbr^d)}
   =\|z\|_{L^{2+\frac 4d}((0,\tau_1)\times \bbr^d)}
   \leq C_*(C_{\tau_1}, B_0(|X_0|_{L^2}))
   =: C_*(C_{\tau_1}).
\end{align*}

Next, we define the random times $\tau_{j+1}$ and $\sigma_{j+1}$
inductively as follows:
$\sigma_0 := 0$, $\sigma_1 := \tau_1$,
and for $j\geq 0$,
\begin{align*}
   & \tau_{j+1}:= \inf\{0\leq t\leq (\tau^* - \sigma_j):
       D_0  \sup\limits_{s\in [0,t]}h_{\sigma_j}(s) \geq \ve_j(t)\} \wedge (\tau^* -\sigma_j), \\
   & \sigma_{j+1} := \sigma_j +\tau_{j+1} (\leq \tau^*),\ \  L :=  \inf\{j\geq 1: \sigma_j = \tau^*\},
\end{align*}
where $h_{\sigma_j}(t) := |\beta(\sigma_j+t) - \beta(\sigma_j)| + t$,
and $\ve_j(t) := \ve_*(C_{\sigma_j+t}, B_0(|X_0|_{L^2}))$ is as Theorem \ref{Thm-Sta-L2}.

We have that $L<\9$, a.s.,
by using similar arguments as in the proof of \cite[Theorem 2.4]{Z18}
involving the $(\frac 12-\kappa)$-H\"older continuity of $t\mapsto h_{\sigma_j}(t)$,
$0<\kappa<\frac 12$.

Then, for each $1\leq j\leq L-1$,
let
$v_{\sigma_j} (t) := e^{-W_{\sigma_j}(t)} X(\sigma_j +t)$
with
$W_{\sigma_j}(t) := W(\sigma_j +t) - W(\sigma_j)$
and set
$u_{\sigma_j}(t) := u(\sigma_j +t)$,
$t\geq 0$.
In view of Theorem \ref{Thm-Rescale-sigma},
we get
\begin{align*}
    & i\p_t v_{\sigma_j}
     = e^{-W_{\sigma_j}} \Delta(e^{W_{\sigma_j}} v_{\sigma_j})
        - F(v_{\sigma_j}) +f(u_{\sigma_j}) v_{\sigma_j}, \\
    & v_{\sigma_j}(0) = X(\sigma_j).
\end{align*}

Similarly as above,
we   compare $v_{\sigma_j}$ with  $\wt{v}_j$
which solves the equation
\begin{align*}
   i\p_t \wt{v}_j =& \Delta \wt{v}_j - F(\wt{v}_j) \\
   \wt{v}_j(0) =& v_{\sigma_j}(0) = X(\sigma_j),
\end{align*}
or, equivalently,
\begin{align*}
    i\p_t \wt{v}_j
    =  e^{-W_{\sigma_j}} \Delta(e^{W_{\sigma_j}} \wt{v}_j)
                   - F(\wt{v}_j) + f(u_{\sigma_j})\wt{v}_j + e_j
\end{align*}
with the error term
\begin{align*}
   e_j = -(b_{\sigma_j} \cdot \na + c_{\sigma_j}) \wt{v}_j
         -  f(u_{\sigma_j}) \wt{v}_j,
\end{align*}
and $b_{\sigma_j}(t) := 2 \na W_{\sigma_j}(t)$,
$c_{\sigma_j}(t) := \Delta W_{\sigma_j}(t) + \sum_{j=1}^d (\p_j W_{\sigma_j}(t))^2$.

Since
$|\wt{v}_j(t)|_{L^2}= |X(\sigma_j)|_{L^2} = |X_0|_{L^2}$,
applying Theorem \ref{Thm-L2GWP-Det} we have
\begin{align*}
   \|\wt{v}_j\|_{L^{2+\frac 4d}(\bbr \times \bbr^d)}
   \leq B_0(|X(\sigma_j)|_{L^2})
   = B_0(|X_0|_{L^2}).
\end{align*}

Then, estimating as in \eqref{esti-e-Xu}, we have
\begin{align*}
    \|e_j\|_{N^0(0,\tau_{j+1}) + L^2(0,\tau_{j+1}; H^{-\frac 12}_1)}
    \leq&  \|(b_{\sigma_j}\cdot \na + c_{\sigma_j}) \wt{v}_j\|_{L^2(0,\tau_{j+1}; H^{-\frac 12}_1)} \\
        &  + \| f(u_{\sigma_j})\wt{v}_j\|_{L^1(0,\tau_{j+1}; L^2)} \\
    \leq& D_0  \sup\limits_{t\in [0,\tau_{j+1}]} h_{\sigma_j}(t)
    \leq \ve_*(C_{\sigma_{j+1}}, B_0(|X_0|_{L^2}) ).
\end{align*}
This, via Theorem \ref{Thm-Sta-L2}, implies that $\bbp$-a.s.
for each $1\leq j\leq L-1$,
\begin{align*}
   &\|X^u\|_{L^{2+\frac 4d}((\sigma_j, \sigma_{j+1})\times \bbr^d)}
   = \|v_{\sigma_j}\|_{L^{2+\frac 4d}((0,\tau_{j+1})\times \bbr^d)} \nonumber \\
   \leq& C_*(C_{\sigma_{j+1}}, B_0(|X_0|_{L^2}))
   =: C_*(C_{\sigma_{j+1}}).
\end{align*}

Thus, taking into account $L<\9$, a.s., we conclude that
\begin{align} \label{esti-Xu-T}
   \sup\limits_{u\in \calu_{ad}} \|X^u\|_{L^{2+\frac 4d}((0,\sigma_L)\times \bbr^d)}
   \leq \sum\limits_{j=0}^{L-1} C_*(C_{\sigma_{j+1}})
   \leq L C_*(C_T) <\9, \ \ a.s.,
\end{align}
which, in particular,
yields \eqref{bdd-Xu-V}
and so the global existence of $X^u$ on $[0,T]$.

Therefore,
applying Corollary \ref{Cor-S0N0-U2} to \eqref{equa-v-contr}
and using \eqref{esti-Xu-T}
we obtain \eqref{E-Xu-U2},
which along with \eqref{embed-Vp-Bp} and \eqref{Hom-S0-U2}
implies \eqref{E-Xu-S0} and \eqref{E-Xu-timereg}.

{\bf $L^\rho(\Omega)$-integrability.}
In the case that $\{e_j\}$ are constants,
equation \eqref{equa-v-contr} reduces to
the deterministic equation
\begin{align} \label{equa-v-det}
   i\p_t v = \Delta v  + \lbb F(v)+ f(u)v,
\end{align}
which implies that $\|v\|_{S^0(0,T)}\in L^\9(\Omega)$
and so \eqref{E-Xu-U2-Integ} follows.

Below we consider the case that $e_j \in L^{\frac{2p}{p-2}}$
for some $p$ satisfying the condition in Hypothesis $(H0)^*$.
Take $q\in (2,\9)$ such that $(p,q)$ is a Strichartz pair.
Set
$M_1^*(t):= \sup_{0\leq r_1<r_2\leq t}$ $|\int_{r_1}^{r_2} e^{-i(t-s)\Delta}X(s)dW(s)|_{L^2}$,
$M_2^*(t):= \sup_{0\leq r_1<r_2\leq t}$ $|\int_{r_1}^{r_2} e^{-i(t-s)\Delta}X(s)dW(s)|_{L^p}$,
$t\in(0,T)$.
Using similar arguments as in the proof of \cite[Proposition 2.7]{FX18.1}
we have that  for any $1\leq  \rho<\9$,
\begin{align} \label{Integ-M12}
   \|M_1^*\|_{L^\rho(\Omega; L^1(0,T))}
   + \|M_2^*\|_{L^\rho(\Omega; L^q(0,T))} <\9.
\end{align}
(See the Appendix for the proof.
Note that, the restrictions $p<\frac{2d}{d-1}$ and $1\leq d\leq 3$
are used here.)
In particular, \eqref{Integ-M12} implies that
$$ \|M_1^*\|_{L^1(0,T)} + \|M_2^*\|_{L^q(0,T)} <\9,\ \ a.s.. $$

Then, we take a partition $\{[t'_j, t'_{j+1}]\}_{j=0}^{l'}$ of $[0,T]$
such that
$t'_{j+1} = \inf\{t>t'_j: \|M_1^*\|_{L^1(t_j, t)}
          =  \ve_* /2\} \wedge T$,
where $\ve_*:= \ve_*(C_0,$ $B_0(|X_0|_{L^2}))$ is as in Theorem \ref{Thm-Sta-L2},
$C_0$ is the deterministic Strichartz constant related to $\{e^{-it\Delta}\}$.
Note that,
$ l' \leq  2\|M^*_1\|_{L^1(0,T)} /\ve_*$.

Similarly,
take another partition $\{[t''_j, t''_{j+1}]\}_{j=0}^{l''}$ of $[0,T]$,
such that
$t''_{j+1} = \inf\{t>t''_j: \|M_2^*\|_{L^q(t_j, t)}
          =  \ve_* /2\} \wedge T$.
Then,
$ l'' \leq  (2\|M^*_2\|_{L^1(0,T)} /\ve_*)^q$.

Moreover,
taking $\delta :=  \ve_* /((D_*+|\mu|_{L^\9})|X_0|_{L^2})$
and $\{t_j\}_{j=0}^{L'+1} =\{t'_j\}_{k=0}^{l'+1} \cup \{j\delta \}_{j=0}^{[T/\delta]}$,
we get that for any $1\leq \rho <\9$,
\begin{align} \label{integ-L}
     L' \leq l'  + [\frac{T}{\delta}]
       \leq& \frac{2}{\ve_*} \|M^*_1\|_{L^1(0,T)} + (\frac{2}{\ve_*} \|M_2^*\|_{L^q(0,T)})^q \nonumber  \\
       & + \frac{1}{\ve^*}T(D_*+|\mu|_{L^\9})|X_0|_{L^2}
       \in L^\rho(\Omega).
\end{align}

Now, for each $t\in  I_j:= [t_j, t_{j+1}]$, $0\leq j\leq L'$,
we infer from \eqref{equa-x} that
\begin{align} \label{equa-x-M}
   X(t) = e^{-i(t-t_j) \Delta} X(t_j)
           + \int_{t_j}^t e^{-i(t-t_j) \Delta} (iF(X) -\mu X -i f(u) X)ds + M_j(t),
\end{align}
where $F(X)=|X|^{\frac 4d}X$, $f(u)=V_0+u\cdot V$,
and $M_j(t)= \int_{t_j}^t e^{-i(t-s)\Delta}X(s)dW(s)$, $t\in[0,T]$.
In order to apply Theorem \ref{Thm-Sta-L2} to obtain
the boundedness of $\|X\|_{L^q(I_j;L^p)}$,
we compare \eqref{equa-x-M} with the equation
\begin{align*}
     i\p_t \wt{v}_j = \Delta \wt{v}_j - F(\wt{v}_j),
\end{align*}
with $\wt{v}_j(t_j) = X(t_j)$, or equivalently,
\begin{align*}
    \wt{v}_j(t) = e^{-i(t-t_j)\Delta}X(t_j)
                  + \int_{t_j}^t e^{-i(t-s)\Delta}(iF(\wt{v}_j) - \mu\wt{v}_j -i f(u)\wt{v}_j +e )ds
\end{align*}
with the error term $e=\mu\wt{v}_j +i f(u)\wt{v}_j$.

Note that, Theorem \ref{Thm-L2GWP-Det} implies
\begin{align*}
    \|\wt{v}_j\|_{L^q(\bbr;L^p)} \leq B_0(|\wt{v}_j(0)|_{L^2}) = B_0(|X(t_j)|_{L^2}) = B_0(|X_0|_{L^2}).
\end{align*}
Moreover, by the construction of $I_j$,
\begin{align*}
   \|M_j\|_{L^q(I_j;L^p)}
    + \|M_j\|_{L^1(I_j;L^2)}
   \leq \|M_1^*\|_{L^1(I_j)} + \|M_2^*\|_{L^q(I_j)}
      \leq  \ve_*,
\end{align*}
and
\begin{align*}
  \|e\|_{N^0(I_j)}
  \leq \|\mu\wt{v}_j +i f(u)\wt{v}_j\|_{L^1(I_j; L^2)}
  \leq& (D_*+|\mu|_{L^\9}) \|\wt{v}_j\|_{C(I_j; L^2)} |I_j| \\
  =& (D_*+|\mu|_{L^\9}) |X_0|_{L^2} \delta
  \leq \ve^*,
\end{align*}
which implies the conditions \eqref{Sta-L2-ve.1} and \eqref{Sta-L2-ve.2} of Theorem \ref{Thm-Sta-L2}.

Thus, applying Theorem \ref{Thm-Sta-L2}
with $\Phi \equiv 0$, $R=M_j$,
$G= -\mu -i f(u)$
and $V(t,s) = e^{-i(t-s)\Delta}$,
we obtain
\begin{align*}
   \|X\|_{L^q(I_j;L^p)} \leq C_*(C_0, B_0(|X_0|_{L^2})) \in L^\9(\Omega).
\end{align*}
Summing over $j$ and using \eqref{integ-L} we get
\begin{align*}
    \|X\|_{L^q(0,T;L^p)} \leq L' C_*(C_0,  B_0(|X_0|_{L^2})) \in L^\rho(\Omega), \  \ \forall 1\leq \rho <\9.
\end{align*}
Then, taking  the Strichartz pair $(\wt{p},\wt{q})$
such that $(\frac{1}{\wt{p}'}, \frac{1}{\wt{q}'}) = (\frac{\a}{p}, \frac{\a}{q})$,
applying Theorem \ref{Thm-Stri} to \eqref{equa-v-contr}, 
using H\"older's inequality and \eqref{Inhom-U2-N0}
we finally obtain 
\begin{align*}
       & \|X\|_{L^2(0,T; H^{\frac 12}_{-1})}  + \|e^{-W}X\|_{\calu^2(0,T)} \\
    \leq& C(T) |X_0|_{L^2}  + C(T) \|X\|^{1+\frac 4d}_{L^q(0,T;L^p)} +  C(T)D_*T |X_0|_{L^2} \in L^\rho(\Omega), \ \ 1\leq \rho <\9.
\end{align*}

Therefore, the proof  in the defocusing mass-critical case is complete.
\hfill $\square$

Below
we prove  Theorem \ref{Thm-VX},
which is important to derive the tightness  in Section \ref{Sec-Proof} below.

{\bf Proof of Theorem \ref{Thm-VX}.}
First,
choose the Strichartz pair $(p,q)=(\a+1, \frac{4(\a+1)}{d(\a-1)})$
and divide $[0,T]$ into $\{I_{j}\}_{j=0}^L:=\{[t_{j},t_{j+1}]\}_{j=0}^L$ such that
$t_{j+1} = \inf\{t>t_j: \|X\|^{\a-1}_{L^q(t_j, t; L^p)} = (2C^*(T))^{-1}\} \wedge T$,
where $C^*(T) \in L^\rho(\Omega)$ is to be specified in \eqref{esti.3} below.
Then,
$L\leq  (2C^*(T))^{\frac{q}{\a-1}}\|X\|^q_{L^q(0,T; L^p)} \in L^\rho (\Omega)$,
$\forall 1\leq \rho <\9$.

Letting $z_k := e^{-W} V_kX$, we have
\begin{align}
   &i\p_t z_k = e^{-W}\Delta (e^{W} z_k) + \lbb |X|^{\a-1} z_k  + f(u) z_k  - 2 e^{-W} \na V_k \cdot \na X - e^{-W}\Delta V_k X, \nonumber \\
   & z_k(0) = V_k X_0,
\end{align}
where $f(u)=V_0+u\cdot V$.
This yields that for every $0 \leq j\leq L$,
\begin{align*}
   z_k(t)
   =& U(t, t_j) z_k(t_j) + (-i)\int_{t_j}^{t} U(t, s)
          \bigg[\lbb |X|^{\a-1} z_k  +f(u) z_k  \\
    & \qquad \qquad \qquad \qquad \qquad \qquad - 2e^{-W} \na V_k \cdot \na X - e^{-W}\Delta V_k X \bigg] ds.
\end{align*}
Then, using
$U(t_j, t) U(t,s) = U(t_j, s)$ we get
\begin{align*}
   U(t_j, t) z_k(t)
   =&  z_k(t_j) + (-i) \int_{t_j}^{t} U(t_j, s)
          \bigg[\lbb |X|^{\a-1} z_k  + f(u) z_k  \\
    & \qquad \qquad \qquad \qquad \qquad - 2e^{-W} \na V_k \cdot \na X - e^{-W}\Delta V_k X \bigg] ds.
\end{align*}
Applying Proposition \ref{Pro-LpLq-UqVq} we get
\begin{align} \label{esti-zk}
  \|U(t_j,\cdot)z_k\|_{V^2(I_j)}
  \leq& C(T)\bigg[1+ |z_k(t_j)|_{L^2}
        + \||X|^{\a-1}z_k\|_{L^{q'}(I_j; L^{p'})}
        + \|f(u) z_k\|_{L^1(I_j; L^2)}   \nonumber \\
      & + 2\|e^{-W} \na V_k \cdot \na X\|_{L^2(I_j; H^{-\frac 12}_1)}
       + \|e^{-W} \Delta V_k X\|_{L^1(I_j; L^2)}  \bigg].
\end{align}
Since $|z_k(t)|_{L^2} \leq |V_k|_{L^\9}|X_0|_{L^2}$,
we have
\begin{align} \label{esti.1}
    &|z_*(t_j)|_{L^2} + \|f(u)z_k\|_{L^1(I_j; L^2)} + \|e^{-W} \Delta V_k X\|_{L^1(I_j;L^2)} \nonumber  \\
     \leq& (|V_k|_{L^\9}+|V_k|_{L^\9}D_*T+|\Delta V_k|_{L^\9}T ) |X_0|_{L^2}
\end{align}
with $D_*$ as in \eqref{esti-xu-sub}.
Moreover,  for any multi-index $\g$, by $(H0)$,
\begin{align*}
   \sup\limits_{0\leq t\leq T} |\p_x^\g (e^{-W} \na V_k)|
   \leq C\sup\limits_{0\leq t\leq T} |\beta(t)|^{|\g|} \<x\>^{-2}.
\end{align*}
Then, estimating as in the proof of \eqref{esti-wtv-LS}
we have for some $l\geq 1$,
\begin{align} \label{esti.2}
    \|e^{-W}\na V_k \cdot \na X\|_{L^2(I_j; H^{-\frac 12}_{1})}
    \leq C_0 \sup\limits_{0\leq t\leq T} |\beta(t)|^{l}\| X\|_{L^2(I_j; H^{\frac 12}_{-1})}.
\end{align}
Thus, plugging \eqref{esti.1} and \eqref{esti.2} into \eqref{esti-zk}
and using
$\||X|^{\a-1}z_k\|_{L^{q'}(I_j; L^{p'})}
\leq T^\theta \|X\|^{\a-1}_{L^q(I_j; L^p)} \|z\|_{L^{q}(I_j; L^p)}$
where $\theta = 1-\frac{d(\a-1)}{4} \geq 0$,
we obtain
\begin{align} \label{esti.3}
  \|U(t_j,\cdot)z_k\|_{V^2(I_j)}
  \leq& C^*(T)  (1+ |X_0|_{L^2} + \|X\|_{L^2(0,T; H^{\frac 12}_{-1})} + \|X\|^{\a-1}_{L^{q}(I_j; L^{p})} \|z\|_{L^{q}(I_j; L^{p})})  \nonumber \\
  \leq&  C^*(T) (1+ |X_0|_{L^2} + \|X\|_{L^2(0,T; H^{\frac 12}_{-1})}) + \frac 12 \|z\|_{L^{q'}(I_j; L^{p'})}
\end{align}
with
$C^*(T):=C(T) (1+|V_k|_{L^\9}(1 + D_*T) +|\Delta V_k|_{L^\9} T + T^\theta + 2C_0 \sup_{0\leq t\leq T} |\beta(t)|^{l})$.
It follows that
\begin{align*}
    \|U(t_j,\cdot)z_k\|_{V^2(I_j)}
     \leq 2 C^*(T)  (1+ |X_0|_{L^2} + \|X\|_{L^2(0,T; H^{\frac 12}_{-1})}).
\end{align*}

Moreover, taking into account
$U(0,\cdot) = U(0,t_j) U(t_j, \cdot)$,
$\|U(0,t_j)\|_{\mathcal{L}(L^2, L^2)}$ $\leq C(T) \in L^\rho(\Omega)$
for any $1\leq \rho <\9$,
we have
\begin{align*}
  \|U(0,\cdot)z_k\|_{V^2(I_j)}
  \leq&  C(T) \|U(t_j,\cdot)z_k\|_{V^2(I_j)} \\
  \leq&  2 C(T)C^*(T)  (1+ |X_0|_{L^2} + \|X\|_{L^2(0,T; H^{\frac 12}_{-1})}).
\end{align*}

Thus, summing over $0\leq j\leq L$
and using \eqref{E-Xu-U2-Integ}
we obtain
\begin{align*}
     \|z_k\|_{\calv^2(0,T)}
     \leq& \sum\limits_{j=0}^L \|U(0,\cdot)z_k\|_{V(I_j)}  \\
     \leq& 2 L C(T) C^*(T) (1+ |X_0|_{L^2}+ \|X\|_{ L^2(0,T; H^{\frac 12}_{-1})}) \in L^\rho (\Omega)
\end{align*}
for any $1\leq \rho<\9$, which implies \eqref{E-VXu-V2}.

The temporal regularity \eqref{E-VXu-timereg} now follows from \eqref{E-VXu-V2}
and the embedding \eqref{embed-Vp-Bp*}.
Therefore, the proof of Theorem \ref{Thm-VX} is complete.
\hfill $\square$

\subsection{Stability with respect to controls} \label{Subsec-Sta-ContrEqua}

In this Subsection,
we prove the continuous dependence of controlled solutions and objective functionals
on control and noise.

Consider the sequence $((\bbx_T)_n, \bbx_n, \beta_n, u_n)$,
$(\bbx_T, \bbx,  \beta, u) \in \caly$, $n\geq 1$.
where
$\caly$ is the space as in Definition \ref{Def-Contr},
$u_n, u\in \calu_{ad}$, and
$\beta_n$, $\beta$ are $\bbr^N$-dimensional Wiener processes.
Assume $\bbp$-a.s.,
\begin{align}
   &u_n \to u,\ \ in\ L^2(0,T; \bbr^m), \ \
   \beta_n \to \beta, \ \ in\ C([0,T]; \bbr^N), \label{conv-un-betan.2} \\
   &\bbx_n \to \bbx,\ \ in\ L^2((0,T)\times \bbr^d), \ \ (\bbx_T)_n \to \bbx_T,\ \ in\ L^2(\bbr^d). \label{conv-bbxn}
\end{align}
Let $X_n$ (resp. $X$) be the controlled solution to \eqref{equa-x}
corresponding to $(u_n, \beta_n)$ (resp. $(u,\beta)$), $n\geq 1$.
We have
\begin{lemma} \label{Lem-Xn-X}
Consider the mass-subcritical case
or the defocusing mass-critical case.
Assume $(H0)$.
Then, for each $X_0\in L^2$, $0<T<\9$,
\begin{align} \label{conv-Xn-X}
     \|X_n - X\|_{S^0(0,T)} \to 0,\ \ as\ n\to \9,\ a.s..
\end{align}
\end{lemma}

{\it \bf Proof. }
In the subcritical case,
\eqref{conv-Xn-X} follows from \cite[(3.11)]{BRZ18}
and Strichartz estimates \eqref{L2-Stri}.
Below we prove the defocusing mass-critical case
by using the stability result Theorem \ref{Thm-Sta-L2}.

Set $W_n(t,x) = \sum_{j=1}^N \mu_j e_j(x) \beta_{j,n}(t)$,
$W(t,x) = \sum_{j=1}^N \mu_j e_j(x) \beta_{j}(t)$,
where
$\beta_{j,n}$ and $\beta_j$
are the $j$-th component of $\beta_n$ and $\beta$ respectively,
$1\leq j\leq N$,
$t\geq 0$, $x\in \bbr^d$.
We may assume $T \geq 1$ without loss of generality.

Using the rescaling transformations
$v_n := e^{-W_n} X_n$ we have
\begin{align} \label{equa-yn*}
    i \p_t v_n
  =& e^{-W_n}\Delta(e^{W_n} v_n) - F(v_n) +  f(u_n) v_n, \\
  v_n(0) =& X_0, \nonumber
\end{align}
where $f(u_n):= V_0 + u_n\cdot V$.
Similarly $v:=e^{-W}X$ solves \eqref{equa-v-contr} with $\lbb =-1$.

We shall compare $v_n$ and $v$.
For this purpose,
we reformulate \eqref{equa-v-contr}
\begin{align} \label{equa-v-contr.2}
    i\p_t v
  =& e^{-W_n}\Delta(e^{W_n} v) - F(v) +  f(u_n) v+e_n,
\end{align}
where the error term
\begin{align*}
   e_n = ((b(t)-b_n(t)) \cdot \na + (c(t)-c_n(t))) v
           + (f(u)-f(u_n))v
\end{align*}
with
$b_n(t)=2\na W_n(t)$,
$c_n(t)= \Delta W_n(t) + \sum_{j=1}^d(\p_jW_n(t))^2$,
and $b(t), c(t)$ defined similarly as in \eqref{equa-v-homo}.

Note that,
by $(H0)$,
$\bbp$-a.s. for any multi-index $\g$,
\begin{align*}
   |\p_x^\g W_n(t,x)| \leq C\sup\limits_{n\geq 1} |\beta_n(t)| \<x\>^{-2} <\9,
\end{align*}
where $C$ is independent of $n$.
Hence, Theorem \ref{Thm-Stri} yields that
Strichartz and local smoothing  estimates hold for the
operator $e^{-W_n}\Delta(e^{W_n} \cdot)$
and the related Strichartz constants are uniformly bounded for all $n$.

Moreover,
in  order to estimate the error term,
we see that
\begin{align*}
   \<x\> \<\na\>^{-\frac 12} ((b(t)-b_n(t))\cdot \na + (c(t)-c_n(t)))
   = \Psi_{p_n(t)} \<x\>^{-1} \<\na\>^{\frac 12},
\end{align*}
where $\Psi_{p_n(t)}:= \<x\> \<\na\>^{-\frac 12} ((b(t)-b_n(t))\cdot \na + (c(t)-c_n(t)))  \<\na\>^{-\frac 12} \<x\>$
is a pseudo-differential operator of order $0$,
satisfying
$|p_n(t)|_{S^0}^{(l)}
\leq C(l) \|\beta_n\|_{C([0,t]; \bbr^N)}$
$\|\beta - \beta_n \|_{C([0,t]; \bbr^N)}$
for any $l\geq 1$.
Then,
estimating as in \eqref{esti-wtv-LS} we get
\begin{align} \label{esti-wtu-LS}
   &\|((b(t)-b_n(t))\cdot \na + (c(t)-c_n(t))) v\|_{L^2(0,T; H^{-\frac 12}_1)} \nonumber \\
   \leq& C \|\beta_n\|_{C([0,T]; \bbr^N)}
          \|\beta - \beta_n \|_{C([0,T]; \bbr^N)} \|v\|_{L^2(0,T; H^{\frac 12}_{-1})}.
\end{align}
This along with H\"older's inequality and \eqref{conv-un-betan.2} yields that $\bbp$-a.s. as $n\to \9$,
\begin{align*}
   \|e_n\|_{N^0(0,T)+ L^2(0,T;H^{-\frac 12}_{1})}
   \leq& \|((b(t)-b_n(t))\cdot \na + (c(t)-c_n(t))) v \|_{L^2(0,T; H^{-\frac 12}_1)} \\
       & + \|(f(u)-f(u_n))v\|_{L^1(0,T; L^2)} \\
   \leq&C \sup\limits_{n\geq 1}\|\beta_n\|_{C([0,T]; \bbr^N)}
          \|v\|_{L^2(0,T;H^{\frac 12}_{-1})} \|\beta - \beta_n \|_{C([0,T]; \bbr^N)} \\
       & + C T^\frac 12  \|v\|_{C([0,T];L^2)}\|u_n - u\|_{L^2(0,T; \bbr^m)}
   \to 0.
\end{align*}

Thus, by virtue of Theorem \ref{Thm-Sta-L2} with $R=0$,
we obtain
\begin{align} \label{conv-vn-v}
    \|v_n - v\|_{L^{2+\frac 4d}((0,T)\times\bbr^d)} \to 0,\ \ as\ n\to \9,\ a.s..
\end{align}
Note that,
for any $n$ large enough,
\begin{align*}
   & \|X_n-X\|_{L^{2+\frac 4d}((0,T)\times\bbr^d)} \\
   \leq&
   \|v_n- v\|_{L^{2+\frac 4d}((0,T)\times\bbr^d)}
   + \|e^{-W}X- e^{-W_n}X\|_{L^{2+\frac 4d}((0,T)\times\bbr^d)}  \\
   \leq&\|v_n- v\|_{L^{2+\frac 4d}((0,T)\times\bbr^d)}
         + e\|W_n-W\|_{C([0,T]; L^\9(\bbr^d; \bbr^N))} \|X\|_{L^{2+\frac 4d}((0,T)\times\bbr^d)} ,
\end{align*}
where in the last step we used the inequality $|e^x-1|\leq e|x|$
for $|x|\leq 1$.

Therefore,
using  \eqref{conv-un-betan.2}, \eqref{conv-vn-v}
we obtain \eqref{conv-Xn-X}
and finish the proof.
\hfill $\square$

\begin{lemma}  \label{Lem-conv}
Assume the conditions of Lemma \ref{Lem-Xn-X} to hold.
Assume additionally that
$\sup_n (\bbe|(\bbx_T)_n|^2 + \bbe \|X_n\|^2_{L^2(0,T; L^2)}) <\9$.
Then, as $n\to \9$,
\begin{align}
      \bbe Re\<X_n(T), (\bbx_T)_n\>_2
     &\to \bbe Re \<X(T), \bbx_T\>_2, \label{App-conv-Xtau} \\
      \bbe\int_0^{T} Re\<X_n(t), \bbx_{n}(t)\>_2 dt
     &\to \bbe \int_0^{T} Re\<X(t), \bbx(t)\>_2 dt. \label{App-conv-intX}
\end{align}
\end{lemma}

{\bf Proof.}
First, we infer from \eqref{conv-bbxn} and Lemma \ref{Lem-Xn-X} that
$\bbp$-a.s.,
\begin{align}
        Re\<X_n(T), (\bbx_T)_n\>_2
     &\to Re \<X(T), \bbx_T\>_2. \label{App-conv-Xtau*}
\end{align}
Moreover, since $|X_n(T)|_{L^2}= |X_0|_{L^2}$,
\begin{align*}
   \sup\limits_{n\geq 1} \bbe |Re\<X_n(T), (\bbx_T)_n\>_2|^2
   \leq |X_0|_{L^2}  \sup\limits_{n\geq 1} \bbe |(\bbx_T)_n|_{L^2}^2<\9,
\end{align*}
which implies the uniform integrability of
$\{Re\<X_n(T), (\bbx_T)_n\>_2\}$.
Hence,
taking into account \eqref{App-conv-Xtau*} we obtain \eqref{App-conv-Xtau}.
The proof of \eqref{App-conv-intX} is similar.
\hfill $\square$

Below we prove the continuous dependence of  objective functionals with respect to
control and noise when $\g_3=0$.
Similarly to \eqref{def-Phi},
define
\begin{align}  \label{Def-Phin}
\Phi_n(u)
    := & \bbe |X_n(T)-(\bbx_T)_n|_{L^2}^2
         + \g_1 \bbe \int_0^T |X_n(t)-\bbx_n(t)|_{L^2}^2 dt \nonumber \\
       &  + \g_2 \bbe \int_0^T |u_n(t)|_m^2dt.
\end{align}

\begin{proposition} \label{Prop-Phi-conti}
Let $\Phi$ be as in \eqref{def-Phi} with $\g_3=0$
and $\Phi_n$ be as in \eqref{Def-Phin} above.
Assume $(H0)$.
Assume additionally that
$\bbe |(\bbx_T)_n|_{L^2}^2 \to \bbe |\bbx_T|_{L^2}^2$
and
$\bbe \int_0^T |\bbx_n(t)|_{L^2}^2 dt \to \bbe \int_0^T |\bbx(t)|_{L^2}^2 dt$.
Then, for each $X_0\in L^2$, $0<T<\9$,
\begin{align} \label{conv-Phin-Phi}
    \Phi_n(u_n) - \Phi(u) \to 0,\ \ as\ n\to \9,\ a.s..
\end{align}
\end{proposition}

The proof
follows by expanding the objective functions $\Phi_n$ and $\Phi$
and then using Lemma \ref{Lem-conv}.

We are now ready to prove Theorem  \ref{Thm-Control-deriv}.

{\bf Proof of Theorem \ref{Thm-Control-deriv}.} The proof is similar to that of Theorem $2.5$ of \cite{BRZ18}.
Actually,
take any approximating controls $\{u_n\} \subseteq \calu_{ad}$ such that,
if $I:= \inf \{\Phi(u); u\in \calu_{ad}, X^u\ satisfies\ \eqref{equa-x}\}$,
\begin{align*}
   I \leq \Phi(u_n) \leq I+ \frac 1n,\ \ n\geq 1.
\end{align*}
When $\g_3>0$,
we have additionally the differentiability of $\{u_n\}$,
which yields
the tightness of associated distributions on $C([0,T]; \bbr^m)$.
Hence,
one can apply the Skorohod representation theorem and
the stability result Proposition \ref{Prop-Phi-conti}
to obtain a relaxed solution to Problem $(P)$ with $\g_3>0$.
We refer to Section 3 of \cite{BRZ18} for more details.
\hfill $\square$

\section{The dual backward stochastic equations} \label{Sec-Backequa}

This section is mainly devoted to the proof of Theorem \ref{Thm-BSPDE}.
We first analyze the variational equation in Subsection \ref{Subsec-Vraequa},
and then in Subsection \ref{Subsec-Backequa} we prove Theorem \ref{Thm-BSPDE}
for the dual backward  equation.
As a direct application,
we obtain the directional derivative of
objective functional.

\subsection{The variational equations}  \label{Subsec-Vraequa}

Consider the equation
\begin{align} \label{equa-var-Psi}
  &id \psi =   \Delta \psi dt  + \lbb   f_1(X^u) \psi dt  + \lbb  f_2(X^u) \ol{\psi} dt -i \mu \psi dt  \nonumber \\
  &\qquad \quad  + V_0 \psi dt + u\cdot V \psi dt  - i \Psi dt +i \psi d W(t),  \\
  &\psi(0)= 0. \nonumber
\end{align}
Here,
$f_j$ satisfies $|f_j(z)| \leq C_* |z|^{\a-1}$, $z\in \mathbb{C}$, $j=1,2$,
$X^u$ is the controlled solution to \eqref{equa-x} related to $u\in \calu_{ad}$,
and $\Psi$ satisfies $e^{-W}\Psi\in \caln^2(0,T)$,
i.e.,
$$\bigg\|\int_0^\cdot U(0,s)e^{-W(s)}\Psi(s)ds \bigg\|_{U^2(0,T)}<\9.$$

We mention that,
in the special case   where $\Psi = i \wt{u} \cdot V X^u$
with
$\wt{u} =v-u$,
$u,v\in \calu_{ad}$,
equation \eqref{equa-var-Psi} is indeed
the variational equation (see Proposition \ref{Prop-Var-0} below).
Here, we consider more general case,
in order to analyze  the dual linearized backward stochastic equation \eqref{equa-back}
in Subsection \ref{Subsec-Backequa} later.

\begin{proposition} \label{Prop-Var-Psi}
Consider the mass-subcritical case
or the defocusing mass-critical case.
Assume $(H0)$.
Let $u\in \calu_{ad}$ and
$e^{-W}\Psi\in \caln^2(0,T)$,
$0<T<\9$.
Then,
there exists a unique global $L^2$-solution $\psi^u$ to \eqref{equa-var-Psi},
satisfying that
\begin{align} \label{E-vf-X0}
       \sup\limits_{u\in \calu_{ad}}\|e^{-W}\psi^u\|_{\calu^2(0,T)}
      \leq C(T) \|e^{-W}\Psi\|_{\caln^2(0,T)},\ \ a.s..
\end{align}
If in addition $e_j$ are constants, $1\leq j\leq N$,
then $C(T) \in L^\9(\Omega)$.
\end{proposition}

\begin{remark}
The additional condition that $\{e_j\}$ are constants
can be reduced to the integrability of
$\|X^u\|_{S^0(0,T)}$ and $\|\psi^u\|_{C(0,T; L^2)}$
(see Remark \eqref{Rem-var-Integ} below)
of which the proof is, however, technically unclear due to
the singular coefficient $h_2(X^u)$.
\end{remark}

{\bf Proof of Proposition \ref{Prop-Var-Psi}.}
We use the idea of \cite{BRZ18}
and give a unified treatment in  the subcritical and critical cases.
For convenience,
we use  $C_T(\geq 1)$
for the constants in Theorem \ref{Thm-Stri}, Proposition \ref{Pro-LpLq-UqVq} and Corollary \ref{Cor-S0N0-U2}.

We first see that
$z:= e^{-W} \psi^u$ satisfies the equation
\begin{align} \label{equa-var-z}
   & i \p_t z = e^{-W}\Delta (e^W z) + \lbb h(X^u, z) + f(u)z -i e^{-W} \Psi, \\
   & z(0)= 0, \nonumber
\end{align}
where
$h(X^u, z) = f_1(X^u) z + f_2 (X^u) e^{-2W} \ol{z}$,
and $f(u):=  V_0 + u \cdot V $.
Note that,
$\sup_{u\in \calu_{ad}} \|f\|_{L^\9((0,T)\times \bbr^d)}
   \leq |V_0|_{L^\9(\bbr^d)}
           + D_K \|V_0\|_{L^\9(\bbr^d;\bbr^m)}
   =:D_*$.

Define the operator $F_1$ on $\calu^2(0,T)$ by
\begin{align*}
   F_1(\phi)(t):= (-i) \int_0^t U(t,s) [\lbb h(X^u, \phi)(s) + f(u(s))\phi(s) - i e^{-W(s)} \Psi(s)] ds,
\end{align*}
where $\phi\in \calu^2(0,T)$, $0\leq t\leq T$.

Take $(p,q)=(\a+1, \frac{4(\a+1)}{d(\a-1)})$
and set
$R_1(t):=  2C_* t^\theta \|X^u\|_{L^q(0,t; L^p)}^{\a} + D_* t$,
where $\theta = 1-\frac{d(\a-1)}{4} \in [0,1)$.
Define
$\tau_1:= \inf \{t\in[0,T]:C^2_t R_1(t) \geq \frac 14 \} \wedge T$,
and
$\mathcal{Z}^{\tau_1}_{M_1}
    := \{\phi\in \calu^{2}(0,T):  \|\phi\|_{\calu^{2}(0,\tau_1)}
    \leq   M_1 := 2 C_{\tau_1} \|e^{-W}\Psi\|_{\caln^2(0,\tau_1)}\}$.
Note that,
$R(t)\to 0$ as $t\to 0$, even in the critical case where $\theta=0$.
Moreover,
$\tau_1$ is an $\{\calf_t\}$-stopping time
and $\mathcal{Z}^{\tau_1}_{M_1}$ is a closed ball of Banach space.

Since $U(0,t)U(t,s)= U(0,s)$,
we deduce that
\begin{align*}
    U(0,t) (F_1(\phi) - F_1(\phi_2)) =& (-i)\int_0^t U(0,s) \big[\lbb h(X^u, \phi_1 - \phi_2)(s) \\
         &\qquad \qquad \qquad + f(u(s))(\phi_1- \phi_2)(s) \big] ds.
\end{align*}
Then,
by Corollary \ref{Cor-S0N0-U2} and the H\"older inequality,
for any $\phi_1, \phi_2  \in \mathcal{Z}^{\tau_1}_{M_1}$,
\begin{align} \label{esti-F1.1}
    &\|F_1(\phi_1)- F_1(\phi_2)\|_{\calu^{2}(0,\tau_1)}    \nonumber  \\
    \leq& C_{\tau_1} \|h(X^u, \phi_1- \phi_2) \|_{L^{q'}(0,\tau_1; L^{p'})}
          + C_{\tau_1}  \|f(u)(\phi_1- \phi_2)\|_{L^{1}(0,\tau_1; L^{2})}  \nonumber \\
    \leq& C_{\tau_1} R_1(\tau_1) (\|\phi_1- \phi_2\|_{L^q(0,\tau_1; L^p)} + \|\phi_1- \phi_2\|_{L^\9(0,\tau_1; L^2)}) \nonumber \\
    \leq& \frac 12 \|\phi_1 - \phi_2\|_{\calu^2(0,\tau_1)}.
\end{align}
Similarly, we have for any $\phi \in \mathcal{Z}^{\tau_1}_{M_1}$,
\begin{align}  \label{esti-F1.2}
     \|F_1(\phi)\|_{\calu^{2}(0,\tau_1)}
   \leq& C^2_{\tau_1} R_1(\tau_1)\|\phi\|_{\calu^{2}(0,\tau_1)}
         + C_{\tau_1} \|e^{-W}\Psi\|_{\caln^2(0,\tau_1)}
   \leq M_1.
\end{align}

Thus,
we infer from \eqref{esti-F1.1} and \eqref{esti-F1.2} that
$F_1: \mathcal{Z}^{\tau_1}_{M_1} \mapsto \mathcal{Z}^{\tau_1}_{M_1}$ is a contraction map,
and so
there exists a solution $z_1 \in \mathcal{Z}^{\tau_1}_{M_1}$
such that $z_1 = F_1(z_1)$,
which implies that
$z_1$ solves \eqref{equa-var-Psi} on $[0,\tau_1]$.
Moreover,
\eqref{esti-F1.2} yields that
\begin{align*}
    \|z_1\|_{\calu^{2}(0,\tau_1)} \leq 2 C_{\tau_1}  \|e^{-W} \Psi\|_{\caln^2(0,\tau_1)}.
\end{align*}

Next we use inductive arguments to extend $z_1$ to
the whole time regime $[0,T]$.
For this purpose,
we set $\sigma_1 :=\tau_1$,
and for each $j\geq 1$,
\begin{align*}
   & \tau_{j+1} := \{t\in [0,T-\sigma_j]: C^2_{\sigma_j+t} R_{j+1}(t) \geq \frac 14 \} \wedge (T-\sigma_j), \\
   & \sigma_{j+1} := \sigma_j + \tau_{j+1} (\leq T),\ \
     L := \inf\{j\geq 1: \sigma_j = T\},
\end{align*}
where $R_{j+1}(t) = 2 C_* t^\theta \|X^u\|^\a_{L^q(\sigma_j, \sigma_j+t; L^p)} + D_*t$.

Suppose that for some $1\leq j\leq L-1$
there exists a solution $z_j$ to \eqref{equa-var-Psi} on $[0,\sigma_j]$,
satisfying $z_j(\cdot) = z_j(\cdot \wedge \sigma_j)$,
and
\begin{align} \label{induc-esti-zj}
   \|z_j\|_{\calu^2(0,\sigma_j)}
   \leq M_j:= j(2 C_{\sigma_j})^j \|e^{-W}\Psi\|_{\caln^2(0,\sigma_j)}.
\end{align}
We shall construct a solution $z_{j+1}$ to \eqref{equa-var-Psi} on $[\sigma_j,\sigma_{j+1}]$ below.

For this purpose, we define the operator $F_{j}$ on $\calu^2(0,T)$ by
\begin{align*}
   F_{j}(\phi)(t):=&  U(\sigma_j+t,\sigma_j) z_j(\sigma_j)
               + (-i) \int_0^t U(\sigma_j+t, \sigma_j+s)
               \bigg[ f(u(\sigma_j+s))\phi(s) \\
              &\qquad  + \lbb  h(X^u(\sigma_j+s), \phi(s))
                 - i e^{-W(\sigma_j+s)} \Psi(\sigma_j+s) \bigg] ds,
\end{align*}
where $ \phi\in \calu^2(0,T)$.
Set $\|\phi\|_{\calu^2_j(\sigma_j, \sigma_{j+1})}
:= \|U(\sigma_j, \sigma_j+\cdot) \phi\|_{U^2(0,\tau_{j+1})}$
and
$\mathcal{Z}^{\tau_{j+1}}_{M_{j+1}}
   :=  \{\phi\in \calu^{2}(0,T):
       \|\phi\|_{\calu^2_j(\sigma_j, \sigma_{j+1})} \leq M_{j+1}  \}$,
where $M_{j+1}$ is as in \eqref{induc-esti-zj}
with $j+1$ replacing $j$.

Then,
estimating as in \eqref{esti-F1.1} and \eqref{esti-F1.2},
we get for any $\phi_1, \phi_2, \phi \in \mathcal{Z}^{\tau_{j+1}}_{M_{j+1}}$,
\begin{align} \label{esti-Fn1.1}
    \|F_{j}(\phi_1) - F_{j}(\phi_2)\|_{\calu^2_j(\sigma_j, \sigma_{j+1})}
   \leq& 2C^2_{\sigma_{j+1}} R_{j+1}(\tau_{j+1})
          \| \phi_1 - \phi_2 \|_{\calu^2_j(\sigma_j, \sigma_{j+1})}  \nonumber  \\
   \leq& \frac 12 \| \phi_1 - \phi_2 \|_{\calu^2_j(\sigma_j, \sigma_{j+1})},
\end{align}
and
\begin{align} \label{esti-Fn1.2}
         \|F_{j}(\phi)\|_{\calu^2_j(\sigma_j, \sigma_{n+1})}
   \leq&  |z_j(\sigma_j)|_{L^2}
        + C^2_{\sigma_{j+1}} R_{j+1}(\tau_{j+1})
              \|\phi\|_{\calu^2_j(\sigma_j, \sigma_{j+1})}  \\
        & + \bigg\|\int_0^\cdot U(\sigma_j, \sigma_j+s) e^{-W(\sigma_j+s)}\Psi(\sigma_j+s) ds \bigg\|_{U^2(0,\tau_{j+1})}.  \nonumber
\end{align}
Since $U(\sigma_j, \sigma_j+s)= U(\sigma_j, 0) U(0,\sigma_j+s)$
and $\|U(\sigma_j, 0) \|_{\mathcal{L}(L^2, L^2)} \leq C_{\sigma_j}$,
\begin{align*}
  &\bigg\|\int_0^\cdot U(\sigma_j, \sigma_j+s) e^{-W(\sigma_j+s)}  \Psi(\sigma_j+s) ds \bigg\|_{U^2(0,\tau_{j+1})}  \\
  \leq& C_{\sigma_j} \bigg\| \int_{\sigma_j}^\cdot U(0, s) e^{-W(s)}  \Psi(s) ds \bigg\|_{U^2(\sigma_j,\sigma_{j+1})}.
\end{align*}
Splitting the integral $\int_{\sigma_j}^\cdot = \int_0^\cdot - \int_0^{\sigma_j}$
and using  $\|v\|_{U^2(\sigma_j,\sigma_{j+1})} \leq \|v \|_{U^2(0,\sigma_{j+1})} $
for any $v\in U^2$,
we obtain
\begin{align} \label{esti-Psin-U2}
    \bigg\|\int_0^\cdot U(\sigma_j, \sigma_j+s) e^{-W(\sigma_j+s)}  \Psi(\sigma_j+s) ds \bigg\|_{U^2(0,\tau_{j+1})}
   \leq 2 C_{\sigma_j}  \|e^{-W} \Psi\|_{\caln^2(0,\sigma_{j+1})}.
\end{align}
Plugging this into \eqref{esti-Fn1.2}
and using \eqref{induc-esti-zj} and \eqref{Hom-LpLq-Uq} we obtain
\begin{align} \label{esti-Fn1.3}
        & \|F_{j}(\phi)\|_{\calu_j^2(\sigma_j, \sigma_{j+1})}  \nonumber  \\
    \leq&  C_{\sigma_{j}} \|z_j\|_{\calu^2_j(0,\sigma_j)}
          + C^2_{\sigma_{j+1}} R_{j+1}(\tau_{j+1}) \|\phi\|_{\calu^2_j(0,\sigma_j)}
          + 2 C_{\sigma_{j+1}} \|e^{-W}\Psi\|_{\caln^2(0,\sigma_{j+1})} \nonumber  \\
    \leq& j  C_{\sigma_{j}}  (2C_{\sigma_{j}})^j \|e^{-W}\Psi\|_{\caln^2(0,\sigma_{j})}
          + \frac 12 M_{j+1}
          + 2 C_{\sigma_{j+1}} \|e^{-W}\Psi\|_{\caln^2(0,\sigma_{j+1})}  \nonumber  \\
    \leq&  M_{j+1}.
\end{align}

Hence,
we deduce from \eqref{esti-Fn1.1} and \eqref{esti-Fn1.3} that
$F_{j}: \mathcal{Z}^{\tau_{j+1}}_{M_{j+1}} \mapsto \mathcal{Z}^{\tau_{j+1}}_{M_{j+1}}$
is a contraction map
and so,
for some $\wt{z}_{j+1} \in \mathcal{Z}^{\tau_{j+1}}_{M_{j+1}}$,
$F_{j}(\wt{z}_{j+1}) = \wt{z}_{j+1}$ on $[0,\tau_{j+1}]$.
Then, letting
\begin{align*}
    z_{j+1}(t) :=
    \left\{
      \begin{array}{ll}
        z_j(t), & \hbox{$t\in [0,\sigma_j]$;} \\
        \wt{z}_{j+1}(t-\sigma_j), & \hbox{$t\in [\sigma_j, T]$,}
      \end{array}
    \right.
\end{align*}
we obtain that $z_{j+1}$ is an $L^2$-solution to \eqref{equa-var-Psi} on $[0,\sigma_{j+1}]$.

Thus, using inductive arguments we obtain a solution $z$ to \eqref{equa-var-Psi} on
the maximal existing time interval $[0,\sigma_L)$.

Below we prove that
$L<\9$, a.s..

For this purpose,
we set
$\delta:=   (16 D_* C^2_T)^{-1}$
and divide $[0,T]$ into finite subintervals $\{[t'_j,t'_{j+1}]_{j=0}^{l'}\}$,
such that
$t_{j+1} = \inf\{t>t_j:
C_*C^2_T T^\theta \|X^u\|^\a_{L^q(t'_j, t; L^p)}$ $=  \frac{1}{32}\} \wedge T$.
Then,
$l' \leq  (32 C_* C_T^2 T^\theta)^{\frac q\a} \|X^u\|^q_{L^q(0,T;L^p)}$.

Hence, letting
$\{t_j\}_{j=0}^{L^*+1} = \{j \delta; 0\leq j\leq [\frac{T}{\delta}] \} \cup
 \{t_j'; 0\leq j\leq l'+1\}$,
we obtain a partition $\{[t_j, t_{j+1}]\}_{j=0}^{L^*}$ of $[0,T]$
satisfying that
\begin{align} \label{esti-l}
   L^* \leq [\frac {T}{\delta}] + l'
      \leq 16 D_*C_T^2 T+ (32 C_* C^2_TT^\theta)^{\frac q \a}\|X^u\|^q_{L^q(0,T;L^p)} <\9,\  a.s.,
\end{align}
and
\begin{align} \label{L*}
   C_T^2 (2C_* T^\theta \|X^u\|_{L^q(t_j, t_{j+1}; L^p)}^{\a}
   + D_* (t_{j+1}-t_j)) \leq \frac 18.
\end{align}

We claim that
\begin{align} \label{L-L*}
   L\leq L^* +1, \ \ a.s..
\end{align}
To this end,
suppose that $L>L^*+1$.
Then, for some $0\leq k\leq L^*$
and some $1\leq j\leq L$,
$[t_k, t_{k+1}]$ contains both $\sigma_j$ and $\sigma_{j+1}$.
If $0\leq k < L^*$,
we have $t_k\leq \sigma_j<\sigma_{j+1} \leq t_{k+1}<T$,
and so, via the definition of $\tau_{j+1}$
and \eqref{L*},
\begin{align*}
   \frac 14
   =  C^2_{\sigma_{j+1}} R_{j+1}(\tau_{j+1})
   \leq  C_T^2 (2C_* T^\theta \|X^u\|_{L^q(t_k, t_{k+1}; L^p)}^{\a}
   + D_* (t_{k+1}-t_k)) \leq \frac 18,
\end{align*}
yielding a contradiction.
If $k=L^*$,
since $L>L^*+1$,
we have $j<L-1$.
But this yields that
$t_{L^*} \leq \sigma_j<\sigma_{j+1} <t_{L^*+1}=T$,
which also leads to a contradiction by reasoning as above.
Thus, we prove \eqref{L-L*}, as claimed.

Therefore, $L<\9$ and $\sigma_L = T$, a.s..
So, the solution exists globally.
Moreover,
we infer from \eqref{induc-esti-zj} that
\begin{align} \label{induc-esti-zj*}
   \|e^{-W}\psi^u\|_{\calu^2(0,T)}
   \leq L (2C_T)^L \|e^{-W}\Psi\|_{\caln^2(0,T)},\ \ a.s.,
\end{align}
which implies \eqref{E-vf-X0} with $C(T)=L (2C_T)^L$.

The uniqueness of solutions to \eqref{equa-var-Psi}
can be proved similarly as in \cite{BRZ18} by using Strichartz estimates.

In the case $e_j$ are constants,
equation \eqref{equa-v-contr} reduces to \eqref{equa-v-det},
which implies that $\|v\|_{S^0(0,T)} \in L^\9(\Omega)$
and so is $\|X^u\|_{S^0(0,T)}$.
Moreover, we have $U(t,s) = e^{-i(t-s)\Delta}$,
implying that
the related Strichartz constant
$C_T$ is deterministic and is independent of $T$.
Hence,
taking into account \eqref{esti-l}, \eqref{L-L*}
and \eqref{induc-esti-zj*}
we obtain $C(T) \in L^\9(\Omega)$.

Therefore, the proof of Proposition \ref{Prop-Var-Psi} is complete.
\hfill $\square$

\begin{remark} \label{Rem-var-Integ}
The $L^\rho(\Omega)$-integrability of $\|e^{-W}\psi^u\|_{\calu^2(0,T)}$
can be reduced to that of $\|X^u\|_{S^0(0,T)}$ and $\|\psi^u\|_{C([0,T];L^2)}$.
Actually,
by
\eqref{esti-Fn1.2}
and \eqref{esti-Psin-U2},
\begin{align} \label{esti-z-calx0.1}
   \|z\|_{\calu^2_j(\sigma_j, \sigma_{j+1})}
   \leq 2  |z(\sigma_j)|_{L^2}
         + 4 C_{\sigma_{j+1}} \|e^{-W}\Psi\|_{\caln^2(0,\sigma_{j+1})},\ 0\leq j< L.
\end{align}
Since
$U(0,\cdot) = U(0,\sigma_j) U(\sigma_j, \cdot)$ and
$\|U(0,\sigma_j)\|_{\mathcal{L}(L^2, L^2)} \leq C_{\sigma_j}$,
by Lemma \ref{Lem-UI},
\begin{align} \label{esti-z-calx0.2}
   \|z\|_{\calu^2(0,T)}
   \leq  \sum\limits_{j=0}^{L-1} \|U(0,\cdot)z\|_{U^2(\sigma_j, \sigma_{j+1})}
    \leq  \sum\limits_{j=0}^{L-1} C_{\sigma_j}\|z\|_{\calu^2_j(0, \tau_{j+1})}.
\end{align}
Then,
plugging \eqref{esti-z-calx0.1} into \eqref{esti-z-calx0.2}
and using $z=e^{-W}\psi^u$
we obtain $\bbp$-a.s.,
\begin{align} \label{esti-z-S0.1}
    \|e^{-W} \psi^u\|_{\calu^{2}(0, T)}
    \leq& 4 C^2_T L \(\|\psi^u\|_{C([0,T]; L^2)}
             +   \|e^{-W}\Psi\|_{\caln^2(0,T)}\),
\end{align}
which along with \eqref{esti-l} and \eqref{L-L*}
yields the statement above.
\end{remark}

Below we concentrate on the variational equation \eqref{equa-var-Psi}
with $\Psi= i\wt{u}\cdot V X^u$, $\wt{u} = v- u$, $v,u \in \calu_{ad}$, i.e.,
\begin{align} \label{equa-var}
  &id \vf =   \Delta \vf dt  + \lbb   h_1(X^u) \vf dt  + \lbb  h_2(X^u) \ol{\vf} dt -i \mu \vf dt  \nonumber \\
  &\qquad \quad  + V_0 \vf dt + u\cdot V \vf dt  + \wt{u}\cdot V X^udt +i \vf d W(t), \nonumber \\
  &\vf(0)= 0,
\end{align}
where $h_1(X^u)$, $h_2(X^u)$ are defined in \eqref{h1-h2}.

In this case,
since $L^1(0,T; L^2) \hookrightarrow N^0(0,T)$,
by \eqref{Inhom-U2-N0},
\begin{align*}
   \|e^{-W}\Phi\|_{\caln^2(0,T)}
   \leq& C(T) \|\wt{u} \cdot V e^{-W}X^u\|_{L^1(0,T; L^2)} \\
   \leq& C(T) D_K \|V\|_{L^\9(\bbr^d; \bbr^m)} |X_0|_{L^2} T<\9,\ \ a.s..
\end{align*}
We deduce from Proposition \ref{Prop-Var-Psi}
that there exists a unique global $L^2$-solution $\vf^{u,\wt{u}}$ to \eqref{equa-var},
satisfying that
\begin{align} \label{vfu-X0}
   \|e^{-W} \vf^{u,\wt{u}}\|_{\calu^2(0,T)}
   \leq C(T) |X_0|_{L^2},\ \ a.s..
\end{align}

Proposition \ref{Prop-Var-0} below shows that
$\vf^{u,\wt{u}}$ is indeed the first order approximation
of the controlled solution to \eqref{equa-x}.

\begin{proposition} \label{Prop-Var-0}
Consider the situations in Proposition \ref{Prop-Var-Psi}.
Assume $(H0)$
and that $e_j$ are constants, $1\leq j\leq N$.
Given any $u, v \in \calu_{ad}$,
let $\wt{u} := v-u$ and
$\vf^{u,\wt{u}}$ be the solution to \eqref{equa-var}, $\ve \in (0,1)$.
Set $u_\ve:= u+ \ve\wt{u}$
and let $X^{u_\ve}$ (resp. $X^u$)
be the controlled solution to \eqref{equa-x}
related to $u_\ve$ (resp. $u$).
Then,
\begin{align} \label{Asym-Xuve-Xu}
    \lim\limits_{\ve \to 0^+}
    \bbe \sup\limits_{t\in [0,T]}
    |\ve^{-1}(X^{u_\ve}(t) - X^u(t)) - \vf^{u,\wt{u}}(t)|_{L^2}^2 =0.
\end{align}
\end{proposition}

{\bf Proof.}
Set $\wt{X}^{u, \wt{u}}_\ve:= \ve^{-1}(X^{u_\ve} - X^u) -\vf^{u,\wt{u}}$,
$\wt{y}_\ve^{u, \wt{u}}:= e^{-W} \wt{X}_\ve^{u, \wt{u}}$.
We have
\begin{align} \label{equa-wtyve}
    \p_t \wt{y}_\ve^{u, \wt{u}}
    = -i\Delta \wt{y}_\ve^{u, \wt{u}} - \lbb i \int_0^1 h(X_{u,r,\ve}, \wt{y}_\ve^{u, \wt{u}})dr
      -i f(u_\ve) \wt{y}_\ve^{u, \wt{u}} + e^{-W} R(\ve, \vf^{u, \wt{u}}),
\end{align}
where
$f(u_\ve)= V_0+u_\ve \cdot V$,
$h(X_{u,r,\ve}) = h_1(X_{u,r,\ve})\wt{y}_\ve^{u, \wt{u}} + h_2(X_{u,r,\ve})e^{-2W} \ol{\wt{y}}_\ve^{u, \wt{u}} $,
$X_{u,r,\ve} := X^u + r(X^{u_\ve} - X^u)$,
$r\in (0,1)$,
and
$R(\ve, \vf^{u,\wt{u}}) := -i (\lbb R_1(\ve) \vf^{u,\wt{u}} + \lbb R_2(\ve) \ol{\vf}^{u,\wt{u}} + \ve \wt{u}\cdot V \vf^{u,\wt{u}})$
with
$R_j(\ve) = \int_0^1 (h_j(X_{u,r,\ve}) - h_j(X^u)) dr$,
$j=1,2$.

In order to prove \eqref{Asym-Xuve-Xu},
it is equivalent to prove that
\begin{align} \label{y-uuve-0}
   \lim\limits_{\ve \to 0} \bbe \|\wt{y}^{u, \wt{u}}_\ve \|^2_{C([0,T];L^2)} =0.
\end{align}

For this purpose,
by virtue of Proposition \ref{Prop-Var-Psi},
we have
\begin{align*}
   \|\wt{y}^{u, \wt{u}}_\ve\|_{\calu^2(0,T)}
    \leq C(T) \|e^{-W}R(\ve, \vf^{u, \wt{u}})\|_{\caln^2(0,T)},
\end{align*}
where $C(T) \in L^\rho(\Omega)$, $\forall 1\leq \rho<\9$.
This along with Corollary \ref{Cor-S0N0-U2} yields
\begin{align*}
  \|y^{u,\wt{u}}_\ve\|_{S^0(0,T)}
  \leq& C'(T) \|e^{-W}R(\ve, \vf^{u, \wt{u}})\|_{N^0(0,T)} \\
  \leq& C'(T) \bigg(\|R_1(\ve) \vf^{u,\wt{u}}\|_{L^{q'}(0,T; L^{p'})}
               + \|R_2(\ve) \ol{\vf}^{u,\wt{u}}\|_{L^{q'}(0,T; L^{p'})}  \\
      &\qquad \qquad   + \|\ve \wt{u}\cdot V \vf^{u,\wt{u}}\|_{_{L^{1}(0,T; L^{2})}} \bigg),
\end{align*}
where
$C'(T)\in L^\rho(\Omega)$,
$(p,q)=(\a+1, \frac{4(\a+1)}{d(\a-1)})$,
and $p', q'$ are the conjugate numbers of $p$ and $q$, respectively.
Thus, by Cauchy's inequality,
\begin{align*}
   \bbe \|y^{u,\wt{u}}_\ve\|^2_{S^0(0,T)}
   \leq& C
         \bigg(\bbe \|R_1(\ve) \vf^{u,\wt{u}}\|^4_{L^{q'}(0,T; L^{p'})}
               + \bbe \|R_2(\ve) \ol{\vf}^{u,\wt{u}}\|^4_{L^{q'}(0,T; L^{p'})}  \\
      &\qquad \qquad   + \bbe \|\ve \wt{u}\cdot V \vf^{u,\wt{u}}\|^4_{_{L^{1}(0,T; L^{2})}} \bigg)^\frac 12.
\end{align*}

Note that, by Corollary \ref{Cor-S0N0-U2} and \eqref{vfu-X0},
\begin{align*}
       \(\bbe \|\ve \wt{u}\cdot V \vf^{u,\wt{u}}\|^4_{L^1(0,T;L^2)} \)^\frac 12
   \leq   C \ve^2
   \to 0.
\end{align*}

Thus, it remains to prove that
as $\ve \to 0$,
\begin{align} \label{E-Rj-0}
    \bbe \|R_1(\ve) \vf^{u,\wt{u}}\|^4_{L^{q'}(0,T; L^{p'})}
               + \bbe \|R_2(\ve) \ol{\vf}^{u,\wt{u}}\|^4_{L^{q'}(0,T; L^{p'})}
    \to 0.
\end{align}

For this purpose,
we take $R_1(\ve) \vf^{u,\wt{u}}$ for an example below,
the argument  for $R_2 (\ve) \ol{\vf}^{u,\wt{u}}$ is similar.
We shall prove that
\begin{align} \label{conv-h1-R1}
   \bbe \bigg\|\int_0^1 (h_1(X_{u,r,\ve}) - h_1(X^u)) \vf^{u,\wt{u}} dr \bigg\|^4_{L^{q'}(0,T; L^{p'})} \to 0, \ \ as\ \ve\to 0.
\end{align}

To this end,
we first note that,  Lemma \ref{Lem-Xn-X} implies that
there exists a null set $N'$ such that
$\bbp (N')=0$
and for each $\omega \notin N'$,
\begin{align} \label{Xue-Xu-0}
   \|X^{u_\ve}(\omega) - X^u(\omega) \|_{L^q(0,T; L^p)} \to 0,\ as\ \ve\to 0.
\end{align}
Below we omit the argument $\omega$ for simplicity.
This implies that
for any subsequence $\{\ve_n\}$,
we can extract a further subsequence (still denoted by $\{\ve_n\}$) such that
\begin{align} \label{Xue-Xu-Lp-0}
     |X^{u_{\ve_n}}(t) - X^u(t)|_{L^p} \to 0,\  \ as\ n\to \9,\  dt-a.e..
\end{align}
In particular, for $dt$-a.e. $t \in (0,T)$,
\begin{align*}
    X^{u_{\ve_n}}(t) \to X^u(t),\  \  in\ measure\ dx,
\end{align*}
which yields that  for $dt$-a.e. $t \in (0,T)$,
\begin{align} \label{h1Xue-h1Xu-0.1*}
   h_1(X_{u,r,\ve_n})(t) \to h_1(X^u)(t),\  \  in\ measure\ dx.
\end{align}
Moreover, by \eqref{Xue-Xu-Lp-0},   for $dt$-a.e. $t \in (0,T)$,  as $n\to \9$,
\begin{align} \label{h1Xue-h1Xu-0.2}
  |h_1(X_{u,r,\ve_n})(t)|_{L^{\frac{\a+1}{\a-1}}}
  =& \frac{\a+1}{2} ||X_{u,r,\ve_n}(t)|^{\a-1}|_{L^{\frac{\a+1}{\a-1}}}
  = \frac{\a+1}{2} |X_{u,r,\ve_n}(t)|^{\a-1}_{L^{p}}  \nonumber \\
  \to & \frac{\a+1}{2} |X^u(t)|^{\a-1}_{L^{p}}
  = |h_1(X^u)(t)|_{L^{\frac{\a+1}{\a-1}}}.
\end{align}
We infer from \eqref{h1Xue-h1Xu-0.1*} and \eqref{h1Xue-h1Xu-0.2} that
$\forall r\in (0,1)$, $dt$-a.e. $t \in (0,T)$,
\begin{align} \label{h1Xue-h1Xu-0.1}
   | h_1(X_{u,r,\ve_n})(t) - h_1(X^u)(t)|_{L^{\frac{\a+1}{\a-1}}} \to 0,\ \ as\ n\to \9.
\end{align}
This, via H\"older's inequality, yields that
\begin{align} \label{h1Xuevf-h1Xuvf-0}
   |(h_1(X_{u,r,\ve_n})(t)  - h_1(X^u)(t) ) \vf^{u,\wt{u}}(t)|_{L^{p'}} \to 0,\ \  as\ n\to \9.
\end{align}

Next, we claim that  $\forall r\in (0,1)$,
$\{| (h_1(X_{u,r,\ve_n}) - h_1(X^u)) \vf^{u,\wt{u}} |^{q'}_{L^{p'}}\}$
is uniformly integrable.

Actually, for any Borel set $A \subseteq (0,T)$,
H\"older's inequality implies that
\begin{align*}
  &\int\limits_A | (h_1(X_{u,r,\ve_n})  - h_1(X^u)) \vf^{u,\wt{u}} |^{q'}_{L^{p'}} dt  \\
  \leq&  C T^{\theta q'} \|\vf^{u,\wt{u}}\|_{L^q(0,T;L^p)}^{q'}
           (\int\limits_A |X^{u_{\ve_n}}|^q_{L^p} dt + \int\limits_A |X^{u}|^q_{L^p} dt)^{\frac{(\a-1)q'}{q}}  \\
  \leq& C(T) \|X^{u_{\ve_n}} - X^u\|^{q'}_{L^{(\a-1)q'}(0,T; L^{p'})}
        + C(T) (\int\limits_A |X^{u}|^p_{L^p} dt)^{\frac{(\a-1)q'}{q}},
\end{align*}
where $\theta = 1-\frac{d(\a-1)}{4} \in [0,1)$.
Taking into account \eqref{Xue-Xu-0} we get
\begin{align*}
   & \limsup\limits_{n\to \9}
   \int\limits_{A} | (h_1(X_{u,r,\ve_n}) - h_1(X^u)(t)) \vf^{u,\wt{u}} |^{q'}_{L^{p'}} dt  \\
   \leq&  C(T) (\int\limits_A |X^{u}|^p_{L^p} dt)^{\frac{(\a-1)q'}{q}}
   \to 0,\ \ as\ |A|\to 0,
\end{align*}
which implies the statement,
as claimed.

Hence, taking into account \eqref{h1Xuevf-h1Xuvf-0}, we get that
for $\forall r\in (0,1)$,
\begin{align} \label{h1Xuevf-h1Xuvf-LqLp-0}
   \|(h_1(X_{u,r,\ve_n})  - h_1(X^u) ) \vf^{u,\wt{u}}\|_{L^{q'}(0,T; L^{p'})}  \to 0,\ \ as\ n\to \9.
\end{align}

We also see that
\begin{align} \label{esti-h1}
    \sup\limits_{n \geq 1}
      \|h_1(X_{u,r,\ve_n})\vf^{u,\wt{u}}& \|_{L^{q'}(0,T; L^{p'})}
   \leq C(\a) T^\theta (\|X^u\|^{\a-1}_{L^{q}(0,T; L^{p})} \|\vf^{u,\wt{u}}\|_{L^q(0,T; L^p)}   \nonumber \\
       &  +  \sup\limits_{n\geq 1}\|X^{u_{\ve_n}}\|^{\a-1}_{L^{q}(0,T; L^{p})} \|\vf^{u,\wt{u}}\|_{L^q(0,T; L^p)} ) <\9.
\end{align}
This, via the bounded  convergence theorem, yields that as $n\to \9$,
\begin{align} \label{h1-vf-0}
   &\bigg\|\int_0^1 (h_1(X_{u,r,\ve_n}) - h_1(X^u))\vf^{u,\wt{u}} dr \bigg\|_{L^{q'}(0,T; L^{p'})}  \nonumber \\
   \leq& \int_0^1 \|(h_1(X_{u,r,\ve_n}) - h_1(X^u))\vf^{u,\wt{u}}  \|_{L^{q'}(0,T; L^{p'})} dr \to 0
\end{align}
for any subsequence $\{\ve_n\}$.
Since $\{\ve_n\}$ is arbitrary,
we conclude that
$\bbp$-a.s. \eqref{h1-vf-0} holds with $\ve$ replacing $\ve_n$.

Finally,  in view of \eqref{E-Xu-S0}, \eqref{vfu-X0} and \eqref{esti-h1},
we infer that
the left-hand side of \eqref{h1-vf-0} with $\ve$ replacing $\ve_n$ is uniformly integrable,
and so we obtain \eqref{conv-h1-R1}.
The proof for $R_2(\ve) \ol{\vf}$ is similar.

Therefore,
we obtain \eqref{E-Rj-0} and finish
the proof of Proposition \ref{Prop-Var-0}.
\hfill $\square$

\subsection{Proof of Theorem \ref{Thm-BSPDE}} \label{Subsec-Backequa}

We prove Theorem \ref{Thm-BSPDE} for the stochastic backward equation \eqref{equa-back}.

{\it \bf Proof of Theorem \ref{Thm-BSPDE}. }
We use the duality arguments as in \cite{BRZ18}
to reduce the analysis of  \eqref{equa-back}
to that of the forward equation \eqref{equa-var-Psi}.

First, we construct approximating solutions to \eqref{equa-back}.
Precisely,
we use the truncation
$h_{j,n}(X^u):= g(\frac{|X^u|}{n}) h_j(X^u)$,
$j=1,2$,
where $g$ is a radial smooth cut-off function  such that $g=1$ on $B_1(\bbr)$, and $g=0$ on $B^c_2(\bbr)$.
Note that $ |h_{1,n}(X^u)| +  |h_{2,n}(X^u)| \leq \min\{\a 2^{\a-1} |g|_{L^\9} n^{\a-1}, \a|X^u|^{\a-1}\}$.

Consider the approximating backward stochastic equation
\begin{align} \label{appro-back-equa}
 &d Y_n= -i\Delta Y_n\,dt - \lbb i h_{1,n}(X^u)Y_n dt +\lbb i  h_{2,n}(X^u)  \ol{Y_n} dt + \mu Y_n dt - iV_0Y_n dt \nonumber  \\
 &\qquad \quad        - i u\cdot V Y_n dt + \g_1 (X^u - \bbx_1) dt  - \sum\limits_{k=1}^N \ol{\mu_k} e_k Z_{k,n} dt +  \sum\limits_{k=1}^N Z_{k,n} d\beta_k(t),\nonumber    \\
 & Y_n(T) = -(X^u(T)-\bbx_T).
\end{align}
For $n\geq 1$,
there exists a
unique $(\mathcal{F}_t)$-adapted solution
$(Y_n,Z_n) (:= (Y_n^u, Z^u_n)) \in  L^2(\Omega; C([0,T]; L^2))  \times (L^2_{ad}(0,T; L^2(\Omega; L^2)))^N$
to \eqref{appro-back-equa}
(see, e.g., \cite{FT02}, \cite{HP91}).

In order to pass to the limit $n\to \9$,
we prove the uniform estimate blow
\begin{align}  \label{bdd-E-Yn-V2}
       \sup\limits_{n\geq 1} \sup\limits_{u\in \calu_{ad}} \|e^{-W}Y_n\|_{L^\rho(\Omega; \calv^2(0,T))} <\9, \ \ 1\leq \rho <\rho_\nu:=2+\nu.
\end{align}
In particular,
by \eqref{Hom-S0-U2}
and $|e^{-W}|=1$,
\begin{align}  \label{bdd-E-Yn-S0}
       \sup\limits_{n\geq 1}\sup\limits_{u\in \calu_{ad}} \|Y_n\|_{L^\rho(\Omega; S^0(0,T))} <\9, \ \ 1\leq \rho <\rho_\nu.
\end{align}

To this end, we define the functional $\Lambda_n$ on the space
$L^\9(\Omega\times (0,T) \times \bbr^d)$,
\begin{align} \label{Oper-Phi}
   \Lambda_n(\Psi) := \bbe& Re \<X^u(T)- \bbx_T, \psi_n(T)\>_2
                      +\g_1 \bbe \int_0^{T} Re\<X^u(t)-\bbx(t), \psi_n(t)\>_2dt,
\end{align}
where $\Psi \in L^\9(\Omega\times (0,T) \times \bbr^d)$,
$X^u$ is the controlled solution to \eqref{equa-x},
and $\psi_n(:= \psi_n(u,\Psi))$  satisfies \eqref{equa-var-Psi}
with $h_{j,n}(X^u)$ replacing $f_j(X^u)$, $j=1,2$.

By It\^o's formula,
we  have for any $\Psi\in L^\9(\Omega\times (0,T) \times \bbr^d)$,
\begin{align} \label{dual-Psi-Y}
    \Lambda_n(\Psi)= \bbe \int_0^{T} Re \<\Psi, Y_n\>_2 dt.
\end{align}
Moreover, by virtue of Proposition \ref{Prop-Var-Psi},
we have
\begin{align} \label{psi-Psi*}
     \sup\limits_{n} \sup\limits_{u\in \calu_{ad}}
     \|e^{-W}\psi_n\|_{\calu^2(0,T)}
   \leq  C(T)  \|e^{-W}\Psi\|_{\caln^2(0,T)},\ \ a.s.,
\end{align}
where $C(T) \in L^\rho(\Omega)$,  $\forall 1\leq \rho <\9$,
and is independent of $\Psi$.
This, via H\"older's inequality,
yields that  for  any $1\leq \rho_1<\rho_2<\9$,
\begin{align} \label{psi-Psi}
   \sup\limits_{n\geq 1} \sup\limits_{u\in \calu_{ad}}\|e^{-W}\psi_n\|_{L^{\rho_1}(\Omega; \calu^2(0,T))}
   \leq C(\rho_1, \rho_2, T) \|e^{-W}\Psi\|_{L^{\rho_2}(\Omega; \caln^2(0,T))},
\end{align}
where $C(\rho_1, \rho_2,T)$ is independent of  $\Psi$.
Then, combining together \eqref{E-Xu-U2-Integ}, \eqref{Oper-Phi} and \eqref{psi-Psi},
we obtain that  for any $\ve\in (0,1)$,
\begin{align} \label{dual-esti}
   |\Lambda_n(\Psi)|
   \leq&  \g_1 \|X^u-\bbx_1\|_{L^{\rho_\nu}(\Omega; L^2(0,T; L^2))} \|\psi_n\|_{L^{\rho_\nu'}(\Omega; L^2(0,T; L^2))} \nonumber  \\
       & + \|X^u(T) - \bbx_T\|_{L^{\rho_\nu}(\Omega; L^2)}   \|\psi_n(T)\|_{L^{\rho_\nu'}(\Omega; L^2)}  \nonumber  \\
   \leq&  C(\rho_\nu, \ve, T)  \|e^{-W}\Psi\|_{L^{(\rho_\nu-\ve)'}(\Omega; \caln^2(0,T))}.
\end{align}
Hence, we conclude from \eqref{dual-Psi-Y} and \eqref{dual-esti} that
\begin{align} \label{esti-Lam}
  \sup\limits_{n\geq 1} \sup\limits_{u\in \calu_{ad}}
  \bigg|\bbe \int_0^T \Re \<\Psi, Y_n\>_2 dt \bigg|
  \leq C(\rho_\nu, \ve, T) \|e^{-W} \Psi\|_{L^{(\rho_\nu-\ve)'}(\Omega; \caln^2(0,T))}.
\end{align}

Now,
setting
$\calz := \{v\in L^\9(\Omega; C_0^\9(0,T; L^2)); \|v\|_{L^{(\rho_\nu -\ve)'}(\Omega; U^2(0,T))} \leq 1\}$
and using Lemmas \ref{Lem-Buv} and \ref{Lem-Dual-UpVp} and \eqref{V*-V.2} we get
\begin{align*}
     \| e^{-W} Y_n\|_{L^{\rho_\nu -\ve}(\Omega; \calv^2(0,T))}
   =& \sup\limits_{v\in \calz}
     |\bbe \Re  B(v, U(0, \cdot) e^{-W} Y_n)| \\
   =& \sup\limits_{v\in \calz}
      \bigg|\bbe \int_0^T \Re \<  v'(t), U(0,t) e^{-W(t)} Y_n(t)\>_2 dt \bigg| \\
   =& \sup\limits_{v\in \calz}
     \bigg|\bbe \int_0^T \Re \<e^{W(t)} U(t,0)v'(t), Y_n(t)\>_2 dt\bigg|,
\end{align*}
Then, by \eqref{esti-Lam},
the right-hand side above is bounded by
\begin{align*}
        C(\rho_\nu,\ve,T) \sup\limits_{v\in \calz}
     \bigg\|\int_0^\cdot v'(t) dt \bigg\|_{L^{(\rho_\nu -\ve)'}(\Omega; U^2(0,T))}
   =& C (\rho_\nu,\ve,T) \sup\limits_{v\in \calz}
     \|v \|_{L^{(\rho_\nu -\ve)'}(\Omega; U^2(0,T))} \\
   \leq& C(\rho,\ve,T),
\end{align*}
where $C(\rho,\ve,T)$ is independent of $n$ and $u\in \calu_{ad}$.
This yields \eqref{bdd-E-Yn-V2}
and so \eqref{bdd-E-Yn-S0}, as claimed.

Below
we choose the Strichartz pair $(p,q)= (\a+1, \frac{4(\a+1)}{d(\a-1)})$
and fix  $2 \leq \rho <\rho_\nu$.
Estimate \eqref{bdd-E-Yn-S0} implies that
for some  $\wt{Y} \in L^{\rho}(\Omega; S^0(0,T))$,
along a subsequence of
$\{n\}\to \9$ (still denoted by $\{n\}$),
\begin{align} \label{conv-Yn-Lp}
    Y_n  \overset{\omega^*}{\rightharpoonup} \wt{Y},\ \ in\  L^{\rho}(\Omega; S^0)),
\end{align}
where ``$\overset{\omega^*}{\rightharpoonup}$''
denotes the weak-star convergence.

Moreover,
since for each $j=1,2$, $h_{j,n}(X^u) \to h_j(X^u)$,  $d\bbp \otimes dt \otimes dx$-a.e.,
and $\sup_{n\geq 1}|h_{j,n}(X^u)| \leq C |X^u|^{\a-1} \in L^{2\rho'}(\Omega; L^{\frac{q}{\a-1}}(0,T;L^{\frac{p}{\a-1}}))$,
by the dominated convergence theorem,
\begin{align} \label{conv-hin}
    h_{j,n}(X^u) \to h_j(X^u), \ in\ L^{2\rho'}(\Omega; L^{\frac{q}{\a-1}}(0,T;L^{\frac{p}{\a-1}})), \ \ as\ n\to \9.
\end{align}
Hence, using H\"older's inequality, \eqref{conv-Yn-Lp} and \eqref{conv-hin}
we obtain
\begin{align} \label{conv-hin-Yn}
    h_{1,n} (X^u) Y_n \overset{\omega}{\rightharpoonup}  h_1(X^u) \wt{Y},\  h_{2,n} (X^u) \ol{Y_n} \overset{\omega}{\rightharpoonup}  h_2(X^u) \ol{\wt{Y}},
    \ \ in\ L^{(2\rho')'}(\Omega; L^{\frac{q}{\a}}(0,T;L^{p'})),
\end{align}
where ``$\overset{\omega}{\rightharpoonup}$'' denotes the weak convergence.

Regarding $Z_{k,n}$, $1\leq k\leq N$,
apply It\^o's formula to $e^{\eta t}|Y_n(t)|_{L^2}^2$ for $\eta$ large enough and
using similar arguments as in the proof of \cite[$(4.44)$]{BRZ18},
involving the Burkholder-Davis-Gundy inequality,
we obtain
\begin{align*}
  &\sup\limits_{n\geq 1}
  \frac{1}{2} \sum\limits_{k=1}^N \bbe (\int_0^T e^{\eta s} |Z_{k,n}|_{L^2}^2 ds)^{\frac{\rho}{2}} \\
  \leq& C(\rho) \sup\limits_{n\geq 1} \bbe V_{T,n}^{\frac \rho 2}
        + C(\rho, N) e^{\frac12\rho \eta T}
          \sup\limits_{n\geq 1} \bbe \|Y_n\|^\rho_{C([0,T]; L^2)},
\end{align*}
where
\begin{align*}
   V_{T,n}:=
   & e^{\eta T} |X^u(T) - \bbx_T|_{L^2}^2 + 2 \g_1 \int_0^T e^{\eta T} |X^u-\bbx|_{L^2}^2 ds \\
   & + \a |g|_{L^\9} e^{\eta T} T^\theta \|X^u\|^{\a-1}_{L^q(0,T; L^p)} \|Y_n\|^2_{L^q(0,T; L^p)}.
\end{align*}
Hence,
using  \eqref{E-Xu-U2-Integ} and \eqref{bdd-E-Yn-S0}
we get
\begin{align} \label{bdd-Zn-L2}
   \sup\limits_{n\geq 1} \sup\limits_{u\in \calu_{ad}} \|Z_{k,n}\|_{L^{\rho}(\Omega; L^2(0,T; L^2))} \leq C <\9,\ \ 1\leq k\leq N.
\end{align}
In particular,
there exists $Z^u_k\in L^{\rho}(\Omega; L^2(0,T; L^2))$  such that
\begin{align} \label{conv-Zn-L2}
    Z_{k,n} \overset{\omega}{\rightharpoonup}  Z^u_k,\ \ in\ L^2(\Omega; L^2(0,T; L^2))
\end{align}
for a further subsequence if necessary,
which implies that
\begin{align} \label{conv-Zn-L2*}
   \int_\cdot^T Z_{k,n} d\beta_k(s) \overset{\omega}{\rightharpoonup}   \int_\cdot^T Z^u_{k} d\beta_k(s),\ \  in\ L^2(\Omega; L^2(0,T; L^2)).
\end{align}

Now,
since equation \eqref{appro-back-equa} is linear with respect to $Y_n$, $Z_{k,n}$,
$1\leq k\leq N$,
we can use \eqref{conv-Yn-Lp}-\eqref{conv-Zn-L2*}
to pass to the limit in \eqref{appro-back-equa}.
Moreover, using  similar arguments
as those below $(4.44)$ of \cite{BRZ18},
we can  obtain a continuous version $Y^u$ of $\wt{Y}$,
such that $(Y^u, Z^u)$ solves \eqref{equa-back} in $H^{-2}$
for all $t\in [0,T]$, $\bbp$-a.s..
This yields the global existence of solutions to \eqref{equa-back}.
The uniqueness of solutions can be proved via the duality \eqref{dual-Psi-Y},
as in the proof of \cite[Proposition $7.2$]{BRZ18}.

Moreover,
we obtain \eqref{E-Yu-V2}, \eqref{E-Zu-L2} and \eqref{E-Yu-S0}
from  \eqref{bdd-E-Yn-V2}, \eqref{bdd-E-Yn-S0} and \eqref{bdd-Zn-L2}.
Then,
in view of \eqref{embed-Vp-Bp*},
we also get \eqref{E-Yu-timereg}.

Therefore, the proof of Theorem \ref{Thm-BSPDE} is complete. \hfill $\square$

As a consequence of Proposition \ref{Prop-Var-Psi},
Theorem \ref{Thm-BSPDE} and the duality \eqref{dual-Psi-Y}
with $\Psi= i\wt{u}\cdot V X^u$,
we obtain the directional derivative of objective functional
$\Phi$
in the case where $\g_3=0$.

\begin{proposition} \label{Prop-Deriv-Phi}
Assume  the conditions of Theorem \ref{Thm-Control} to hold.
Then, for each $X_0\in L^2$ and any $u, v\in \calu_{ad}$, we have
\begin{equation}\label{e4.1}
\lim_{\vp\to0}\frac1\vp\,(\Phi(u+\vp\wt u)-\Phi(u))=\E\int^{T}_0\eta(u)(t)\cdot{\wt u}(t)dt,
\end{equation}
where $\wt{u} = v-u$, and
\begin{equation}\label{eta-dPhi}
\eta(u)=2\(\g_2 u-{\rm Im}\int_{\rr^d}V(x) X^u(x) \ol{Y^u}(x)dx \).
\end{equation}
Here, $(Y^u, Z^u)$ is the solution to the dual backward stochastic equation \eqref{equa-back}.
\end{proposition}

\section{Proof of Theorem \ref{Thm-Control}} \label{Sec-Proof}

The key idea here is to apply Ekeland's principle (see \cite[Theorem 1]{E79} or \cite{E74})
to obtain  approximating controls,
which are minimizers of perturbed objective functionals
and so can be characterized explicitly
by using the calculus of subdifferential in the sense of Rockafellar (\cite{R79}).
Then,
by virtue of the new temporal regularities \eqref{E-VXu-timereg} and \eqref{E-Yu-timereg}
of controlled and backward solutions respectively,
we can obtain the tightness of the associated distributions
of approximating controls
which, combined with Skorohod's representation theorem,
enables us to obtain
a relaxed optimal control for Problem $(P)$.

To be precise, let us first deduce from
Proposition \ref{Prop-Phi-conti}  that,
when $\g_3=0$,
$\Phi$ is  continuous on the metric space $\calu_{ad}$ endowed with the distance
$d(u,v)= \|u-v\|=(\bbe \int_0^T |u(t)-v(t)|_m^2 dt)^{1/2}$.
Then, by virtue of Ekeland's variational principle,
for every $n\in\nn$,
we have $u_n\in \calu_{ad}$ such that
\begin{align}
 &\Phi(u_n)\le\Phi(u)+\dd\frac{1}{n}\,d(u_n,u),\ \ \forall  u\in\calu_{ad},   \label{inf-phi-n}
\end{align}
which implies that
\begin{equation}\label{un-mini}
u_n={\rm arg\,min}\left\{\Phi(u)+\frac{1}{n}\,\|u_n-u\|;\ u\in\calu_{ad}\right\}.
\end{equation}

We have the geometric characterization of $\{u_n\}$ as specified below.

\begin{lemma} \label{Lem-chara-contr}
Assume the conditions of Theorem \ref{Thm-Control} to hold.
Let $X_n$ (resp. $(Y_n, Z_n)$)
be the solution to \eqref{equa-x} (resp. \eqref{equa-back}),
corresponding to $u_n$ above, $n\geq 1$.
We have
\begin{equation}\label{un-char}
u_n(t)=P_K\(\dd\frac 1{\g_2}\,{\rm Im}\int V(x) X_n(t,x) \ol{Y_n}(t,x)dx-\frac1{2\g_2 n}\,\eta_n(t)\),
\end{equation}
where $\eta_n \in \p(\|u_n-u\|)|_{u=u_n}$ satisfying
$\E\int^{T}_0|\eta_n(t)|^2_mdt=1$,
``$\p$'' means the subdifferential,
and $P_K$ is the projection operator on  $K$.
\end{lemma}

{\bf Proof.}
We reformulate \eqref{un-mini} as follows
\begin{align*}
   u_n ={\rm arg\,min}\left\{\wt{\Phi}(u) + \frac 1n \|u_n-u\|, u\in L^2_{ad}(0,T; \bbr^m)   \right\},
\end{align*}
where $\wt{\Phi}(u) = \Phi(u) + I_{\calu_{ad}}(u)$,
$I_{\calu_{ad}}(u) =0$ if $u\in \calu_{ad}$,
or $I_{\calu_{ad}}(u) = \9$ otherwise.
Then, applying the subdifferential in the sense of Rockafellar
(see \cite{R79}, see also \cite{BRZ18}) we obtain
\begin{align} \label{0-un.1}
    0 \in \p (\wt{\Phi}(u) + \frac 1n \|u_n-u\|)|_{u=u_n}.
\end{align}
Moreover, by virtue of Theorem $2$ of \cite{R79} we get
\begin{align} \label{0-un.2}
   \p (\wt{\Phi}(u) + \frac 1n \|u_n-u\|) |_{u=u_n}
   \subseteq& \p  \wt{\Phi}(u)|_{u=u_n} + \frac 1n  \p (\|u_n-u\|)|_{u=u_n}  \nonumber \\
   \subseteq& \eta(u_n) + \p I_{\calu_{ad}}(u_n) +  \frac 1n  \p (\|u_n-u\|)|_{u=u_n}  \nonumber \\
   =& \eta(u_n) + \caln_{\calu_{ad}}(u_n) + \frac 1n  \p (\|u_n-u\|)|_{u=u_n},
\end{align}
Here,
$\eta(u_n)$  is the directional derivative of $\Phi$ at $u_n$ obtained in \eqref{eta-dPhi},
and $\caln_{\calu_{ad}}(u_n)$ is the normal cone to $\calu_{ad}$ at $u_n$, i.e.,
\begin{align*}
   \caln_{\calu_{ad}}(u_n)
   = \{v\in L^2_{ad}(0,T; \bbr^m);\ \<v, u_n - \wt{v} \>\geq 0,\ \  \forall \wt{v} \in \calu_{ad} \},
\end{align*}
where $\<\ ,\  \>$ denotes the inner product of $L^2_{ad}(0,T; \bbr^m)$.
We also have the following characterization of $\caln_{\calu_{ad}}(u_n)$
(see \cite[(5.6)]{BRZ18})
\begin{align} \label{equa-cone}
\mathcal{N}_{\calu_{ad}}(u_n)=\{v\in L^2_{ad}(0,T;\rr^m);\  v\in N_K(u_n),\ \ a.e.\ on\ \Omega  \times  (0,T).\},
\end{align}
where $N_K(u_n)$ is the normal cone to $K\subset \rr^m$ at $u_n\in K$, i.e.,
$$N_K(u_n)=\{v\in\rr^m;\  v\cdot( u_n-{\wt v}) \ge0,\ \ \ff\wt v\in K\}.$$

Thus, plugging \eqref{0-un.2} into \eqref{0-un.1} we come to
\begin{align*}
    0 \in \eta(u_n) + \caln_{ad}(u_n) + \frac 1n \partial (\|u_n-u\|)|_{u=u_n},
\end{align*}
which, via \eqref{equa-cone}, yields that
for some $\zeta_n \in N_K(u_n)$, $d\bbp \times dt$-a.e.,
and $\eta_n \in\partial (\|u_n-u\|)|_{u=u_n}$,
\begin{align} \label{un-zetan-etan}
    \eta(u_n) + \zeta_n + \frac 1n \eta_n = 0.
\end{align}
Taking into account \eqref{eta-dPhi},
we arrive at
\begin{equation}\label{e4.12}
\barr{r}
\dd u_n(t)+\frac1{2\gamma_2}\,\zeta_n(t)
= \frac 1{\gamma_2}\,{\rm Im}\int_{\rr^d} V(x)  X_n(t,x) \ol{Y_n}(t,x)dx-\frac1{2\g_2 n}\,\eta_n(t),\earr\end{equation}
a.e. on $(0,T)\times \Omega$,
where  $X_n (:=X^{u_n})$ and $(Y_n,Z_n) (:= (Y^{u_n}, Z^{u_n}))$ are the solutions to
\eqref{equa-x} and
\eqref{equa-back} corresponding to $u_n$, respectively.

Therefore,
applying the projection operator $P_K$ to both sides of \eqref{e4.12}
we obtain  \eqref{un-char} and finish the proof.
\hfill $\square$

Below we prove the crucial tightness of
distributions of $\{u_n\}$ on $L^1(0,T; \bbr^m)$.
\begin{lemma} \label{Lem-tight-contr}
Consider the situations in Lemma \ref{Lem-chara-contr}.
Then, the induced probability measures $\bbp \circ u_n^{-1}$, $n\geq 1$,
are tight on $L^1(0,T; \bbr^m)$.
\end{lemma}

{\bf Proof.}
In view of \cite[Lemma $A.2$]{BRZ18},
we need to verify that
\begin{align} \label{tight-space}
  \lim\limits_{R\to \9} \limsup\limits_{n \to \9} \bbp \bigg\{ \int_0^T  |u_n(t) |_m dt > R \bigg\} =0,
\end{align}
and for any $\ve >0$,
\begin{align} \label{tight-time}
    \lim\limits_{\delta \to 0} \limsup\limits_{n\to \9} \bbp \bigg\{\sup\limits_{0<h\leq \delta}
       \int_0^{T-h} |u_n(t+h)-u_n(t)|_m dt  >\ve \bigg\} =0.
\end{align}

The estimate \eqref{tight-space} follows immediately from the uniform boundedness of $\{u_n\}$.
In order to prove \eqref{tight-time},
using Markov's inequality,
we only need to prove that for some positive exponent $b >0$ and for any $\delta\in(0,1)$,
\begin{align} \label{L1-h-un}
   \limsup\limits_{n\to \9} \bbe \sup\limits_{0<h\leq \delta}  \int_0^{T-h}  | u_n(t+h) - u_n(t)|_m  dt  \leq C \delta^b.
\end{align}

For this purpose,
since $P_K$ is Lipschitz,
using  \eqref{un-char}
we have
\begin{align} \label{esti-un-h.1}
    |u_n(t+h)&-u_n(t)|_m
  \leq \frac{1}{2\g_2 n} (|\eta_n(t+h)|_m + |\eta_n(t)|_m)  \nonumber \\
      & + \frac{1}{\g_2} \bigg|\int V X_n(t+h)\ol{Y}_n(t+h) dx - \int V X_n(t)\ol{Y}_n(t) dx \bigg|_m.
\end{align}

First note that, by H\"older's inequality,  as $n\to \9$,
\begin{align} \label{esti-etan}
   \bbe \sup\limits_{0\leq h\leq \delta}
  \int_0^{T-h} \frac{1}{2\g_2 n} (|\eta_n(t+h)|_m + |\eta_n(t)|_m) dt
  \leq &  \frac{1}{\g_2 n} T^\frac 12 (\bbe \int_0^{T}  |\eta_n(t)|_m^2 dt)^\frac 12 \nonumber   \\
   =&   \frac{1}{\g_2 n} T^\frac 12 \to 0.
\end{align}

Let us treat the second term on the right-hand side of \eqref{esti-un-h.1}.
For simplicity,
set
$\wt{X}_n := U(0,t)e^{-W(t)}(VX_n(t))$,
$\wt{Y}_n := U(0,t)e^{-W(t)}Y_n(t)$.
By \eqref{V*-V},
\begin{align} \label{esti-un-h.2}
    & \bigg|\int V X_n(t+h)\ol{Y}_n(t+h) dx - \int V X_n(t)\ol{Y}_n(t) dx \bigg|_m  \nonumber \\
  =&\bigg|\<\wt{X}_n(t+h), \wt{Y}_n(t+h)\>_2 - \<\wt{X}_n(t), \wt{Y}_n(t)\>_2\bigg|_m  \\
  \leq&\bigg|\<\wt{X}_n(t+h) - \wt{X}_n(t) , \wt{Y}_n(t+h)\>_2\bigg|_m
        + \bigg| \<\wt{X}_n(t), \wt{Y}_n(t+h) - \wt{Y}_n(t)\>_2\bigg|_m. \nonumber
\end{align}
Note that, by Cauchy's inequality,
\begin{align} \label{esti-wtXn}
   & \int_0^{T-h} \bigg|\<\wt{X}_n(t+h) - \wt{X}_n(t) , \wt{Y}_n(t+h)\>_2 \bigg|_m dt \nonumber \\
  \leq& T^\frac 12 \|\wt{Y}_n\|_{C([0,T]; L^2)} \(\int_0^{T-h} |\wt{X}_n(t+h) - \wt{X}_n(t)|^2_2 dt\)^\frac 12.
\end{align}
Taking into account \eqref{embed-Up-Vp} and \eqref{embed-Vp-Bp*}  we obtain
\begin{align} \label{esti-un-h.3}
    &\int_0^{T-h} \bigg|\<\wt{X}_n(t+h) - \wt{X}_n(t) , \wt{Y}_n(t+h)\>_2 \bigg|_m dt  \nonumber \\
   \leq&  C  T^{\frac 1 2}  h^{\frac{1}{2}} \|e^{-W}Y_n\|_{\calv^2(0,T)} \|e^{-W}VX_n\|_{\calv^2(0,T)}.
\end{align}
Similarly,
estimating as above
and using \eqref{embed-Up-Vp} and \eqref{embed-Vp-Bp*} we obtain
\begin{align} \label{esti-un-h.4}
   & \bigg| \<\wt{X}_n(t), \wt{Y}_n(t+h) - \wt{Y}_n(t)\>_2 \bigg|_m dt \nonumber \\
   \leq& C T^\frac 12 \|\wt{X}_n\|_{C([0,T]; L^2)} \(\int_0^{T-h} |\wt{Y}_n(t+h) - \wt{Y}_n(t)|^2_2 dt\)^\frac 12 \nonumber \\
   \leq& C T^\frac 12 h^\frac 12 \|e^{-W}Y_n\|_{\calv^2(0,T)} \|e^{-W}X_n\|_{\calu^2(0,T)}
\end{align}
Thus, combining \eqref{esti-un-h.2}, \eqref{esti-un-h.3} and \eqref{esti-un-h.4} together we obtain
\begin{align*}
        & \int_0^{T-h} \bigg|\int V X_n(t+h)\ol{Y}_n(t+h) dx - \int V X_n(t)\ol{Y}_n(t) dx \bigg|_m dt  \nonumber \\
    \leq& CT^\frac 12h^\frac 12  \|e^{-W}Y_n\|_{\calv^2(0,T)} (\|e^{-W}VX_n\|_{\calv^2(0,T)}+ \| e^{-W}X_n\|_{\calu^2(0,T)}).
\end{align*}
By virtue of uniform estimates \eqref{E-Xu-U2-Integ}, \eqref{E-VXu-V2} and \eqref{E-Yu-V2},
we obtain
\begin{align*}
           \bbe \sup\limits_{0\leq h\leq \delta}
        \int_0^{T-h} \bigg|\int V X_n(t+h)\ol{Y}_n(t+h) dx - \int V X_n(t)\ol{Y}_n(t) dx \bigg|_m dt
    \leq  CT^\frac 12 \delta^\frac 12,
\end{align*}
where $C$ is independent of $n$.
This along with \eqref{esti-un-h.1} and \eqref{esti-etan}
implies  \eqref{L1-h-un} with $b= 1/2$,
thereby yielding the tightness of $\{\bbp\circ u_n^{-1}\}$
on $L^1(0,T;\bbr^m)$.

Therefore,
the proof of Lemma \ref{Lem-tight-contr} is complete.
\hfill $\square$ \\

We are now ready to prove Theorem \ref{Thm-Control}.

{\bf Proof of Theorem \ref{Thm-Control}.}
We use similar arguments as in Section 3 of \cite{BRZ18}.
Set   $\calx_n:=(\bbx_T, \bbx, \beta, u_n)$ in the space
$\mathcal{Y}:=   L^2(  \bbr^d) \times L^2( (0,T)\times \bbr^d) \times C([0,T];\bbr^N) \times L^1(0,T;\bbr^m)$,
where $\beta=(\beta_1,\cdots,\beta_N)$ .
Lemma \ref{Lem-tight-contr} yields the tightness of the induced probability measures of $\{\calx_n\}_{n\geq 1}$ on $\mathcal{Y}$.
Tihs,
via Prohorov's theorem and Skorohod's representation theorem,
implies that there exist a probability space
$(\Omega^*, \mathcal{F}^*, \bbp^*)$
and squences $\calx^*_n := ((\bbx^*_{T})_n, \bbx^*_{n},\beta^*_n, u^*_n)$,
$\calx^* :=(\bbx^*_T,\bbx^*, \beta^*, u^*)$ in $\mathcal{Y}$, $n\in \mathbb{N}$,
such that the joint distributions of $\calx_n^*$ and $\calx_n$
coincide,
$u^*_n, u^* \in K$, $d\bbp^* \times dt$-a.e.,
and
\begin{align} \label{conv-calxn-calx}
  \calx^*_n \to \calx^*,\ \ in\ \caly,\ as\ n\to \9,\ \bbp^*-a.s..
\end{align}
In particular,
by the bounded dominated convergence theorem,
\begin{align}
   u^*_n \to u^*,\ \ in\ L^2(0,T; \bbr^m),\ as\ n\to \9,\ \bbp^*-a.s..
\end{align}
We also note that
$((\bbx_T^*)_n, \bbx^*_{n}, \beta^*_n)$ and $(\bbx^*_T, \bbx^*, \beta^*)$
have the same distributions as $(\bbx_T, \bbx, \beta)$.

Then,
for each $n\geq 1$,
define $\calf^*_{t,n}:=\sigma(\calx^*_n(s), s\leq t)$.
It follows that
$\bbx^*_n(T) \in L^2(\Omega, \mathcal{F}^*_{T,n}, \bbp^*; L^2)$,
$\bbx^*_{n} \in L^2_{ad^*}(0,T; L^2(\Omega; L^2))$,
$u^*_n \in \calu_{ad^*}$,
and that
$(\beta^*_n(t), \calf^*_{t,n})$, $t\in [0,T]$, is a Wiener process.
Moreover, Theorem \ref{Thm-Equa-X} implies
a unique strong $L^2$-solution $X^*_n$ to \eqref{equa-x}
related to $(\beta^*_n, u^*_n)$.
Hence, we deduce that
$(\Omega^*, \mathcal{F}^*, \bbp^*, (\mathcal{F}^*_t)_{t\geq 0}, (\bbx^*_T)_n, \bbx_n^*, \beta_n^*, u_n^*, X_n^*)$
is an admissible system
in the sense of Definition \ref{Def-Contr}.

Similarly, let $\calf^*_t:= \sigma(\calx^*(s), s\leq t)$
and $X^*$ be the unique strong $L^2$-solution to \eqref{equa-x} corresponding to
$(\beta^*, u^*)$. Then,
$(\Omega^*, \mathcal{F}^*, \bbp^*, (\mathcal{F}^*_t)_{t\geq 0}, \bbx^*_T, \bbx^*,$
$\beta^*, u^*, X^*)$
is also an admissible system.

We also note that,
since the controlled solution to \eqref{equa-x} is a measurable map of
Wiener processes and controls,
the distributions of $((\bbx^*_T)_n, \bbx^*_{1,n}, \beta_n^*, u_n^*,$ $X_n^*)$
and
$(\bbx_T, \bbx_{1}, \beta, u_n, X^{u_n})$ coincide.
In particular,
if $\Phi^*_n$ is defined similarly as in \eqref{Def-Phin}
on the new filtrated probability space,
we have $\Phi_n^*(u^*_n) = \Phi(u_n)$.

Thus,
in view of  Proposition \ref{Prop-Phi-conti} and \eqref{inf-phi-n},
we obtain that, if $\g_3=0$,
\begin{align*}
    \Phi^*(u^*) = \lim\limits_{n\to \9} \Phi_n^*(u^*_n)
    =  \lim\limits_{n\to \9} \Phi(u_n)
    =  \inf\{ \Phi(u); u\in \calu_{ad}\},
\end{align*}
which yields \eqref{e2.10}.

Finally, the characterization formula \eqref{e5.1}
of optimal control $u^*$
can be proved similarly as in the proof of Lemma \ref{Lem-chara-contr}.

Therefore, the proof of Theorem \ref{Thm-Control} is complete.
\hfill $\square$

\section{Appendix} \label{APPDIX}

{\bf Proof of Proposition \ref{Prop-V-l2}.}
First,
since $\Re \Phi=0$, $\ol{e^{-\Phi}} = e^{\Phi}$
and $V(t,0)V(0,t) = Id$,
we see that
\eqref{V*-V-L2} follows from \eqref{V*-V}.

Below we prove \eqref{V*-V}.
It is equivalent to prove it for $v\in H^2$,
since $H^2$ is dense in $L^2$.
Let $A(t)(\cdot):= -ie^{-\Phi(t)}\Delta (e^{\Phi(t)} \cdot)$
and $A^*(t)$ be the dual operator of $A(t)$, $t\in \bbr^+$.
Note that, $A^*(t)=-A(t)$.

On the one hand,
since
$\p_t V(t,0)v = A(t) V(t,0)v$
and $(A(t))^* = - A(t)$,
we have for any $w\in H^2$,
\begin{align*}
   \p_t\<V^*(t,0)v, w\>_2
   =  \<v, \p_t V(t,0)w\>_2
   =\<u, A(t)V(t,0)w\>_2
   =\<-V^*(t,0) A(t)v, w\>_2,
\end{align*}
which implies that
\begin{align} \label{equa-U*}
   \p_t V^*(t,0)v = -V^*(t,0) A(t)v, \ \ V^*(0,0)v= v.
\end{align}
On the other hand,
since $V(0,t)V(t,0) = Id$,
we have that 
%for any $w\in H^2$,
%\begin{align*}
%   0=&(\p_t V(0,t)) V(t,0)w + V(0,t) \p_t V(t,0)w \\
%    =& (\p_t V(0,t))  V(t,0)w + V(0,t) A(t)V(t,0)w.
%\end{align*}
%Then, taking $v :=V(0,t)w$ we obtain
\begin{align} \label{equa-U}
   \p_t V(0,t)v= - V(0,t) A(t)v, \ \ V(0,0)v = v.
\end{align}
Thus,
we infer from \eqref{equa-U*} and \eqref{equa-U} that
$V^*(t,0)$ and  $V(0,t)$
satisfy the same equation with the same initial data.
In view of the uniqueness of solutions
(see e.g. \cite{D96}),
we obtain \eqref{V*-V}.

Regarding \eqref{V*-V.2},
since $Id=V(0,t)V(t,0)$,
we deduce from \eqref{V*-V} that
\begin{align*}
  Id= V^*(t,0)V^*(0,t) = V(0,t)V^*(0,t).
\end{align*}
Then,
applying $V(t,0)$ to both sides above
we prove \eqref{V*-V.2}.
\hfill $\square$

{\bf Proof of Lemma \ref{Lem-Embed-UpVp*}.}
Let $I_j=[jh,(j+1)h]$, $0\leq j\leq [\frac{T}{h}]-1=:L^*$.
Since $v\in C([0,T]; H)$,
there exists $t_j\in I_j$ such that
$\sup_{t\in I_j} \|v(t+h)-v(t)\|_H = \|v(t_j+h)-v(t_j)\|_H$,
 $0\leq j\leq L^*$.
Then,
\begin{align*}
  \int_0^{T-h} \|v(t+h)-v(t)\|_H^p dt
  \leq& h \sum\limits_{j=0}^{L^*} \|v(t_j+h)-v(t_j)\|_H^p.
\end{align*}
Moreover, letting $\{t_j'\}_{j=0}^{L} = \{t_j+h, t_j\}_{j=0}^{L^*}$
we have that
\begin{align*}
   \sum\limits_{j=0}^{L^*} \|v(t_j+h)-v(t_j)\|_H^p
   \leq 2^{p+1} \sum\limits_{j=0}^{L-1} \|v(t'_{j+1})-v(t'_j)\|_H^p
   \leq 2^{p+1} \|v\|_{V^p(0,T)}^p.
\end{align*}
Thus, combining the estimates above we prove \eqref{embed-Vp-Bp*}.
\hfill $\square$

{\bf Proof of Proposition \ref{Pro-LpLq-UqVq}.}
First we prove \eqref{Hom-LpLq-Uq}.
It suffices to prove the case that $V(t_0, \cdot) u$ is a $U^q$-atom of the form
$V(t_0, t) u(t) = \sum_{j=1}^n u_j {\calx}_{[t_j, t_{j+1})}(t)$, $t\in I$,
such that $\{t_j\}_{j=1}^n \subseteq I$, $\sum_{j=1}^n |u_j|_{L^2}^q = 1$.
This yields that
$u(t) =  \sum_{j=1}^n V(t, t_0) u_j {\calx}_{[t_j, t_{j+1})}(t)$,
and so
\begin{align}  \label{esti-u-LqLp-Uq.0}
   \|u\|^q_{L^q(I; L^p)}
   = \sum\limits_{j=1}^n \|V(\cdot, t_0)u_i\|^q_{L^q(t_j, t_{j+1}; L^p)}.
\end{align}
Using Theorem \ref{Thm-Stri}
and $\|V(t_j, t_0)\|_{\mathcal{L}(L^2, L^2)} \leq C(T) \in L^\rho(\Omega)$
we have
\begin{align} \label{esti-u-LqLp-Uq.1}
   \|V(\cdot, t_0)u_i\|^q_{L^q(t_j, t_{j+1}; L^p)}
   =& \|V(\cdot, t_j) V(t_j, t_0)u_j\|^q_{L^q(t_j, t_{j+1}; L^p)} \nonumber \\
   \leq& C_T^q |V(t_j, t_0) u_j |^q_{L^2}
   \leq (C'(T))^q |u_j|^q_{L^2},
\end{align}
where
$C'(T)$ is independent of $(p,q)$ and
$C'(T) \in L^\rho(\Omega)$ for any $1\leq \rho <\9$.
Plugging \eqref{esti-u-LqLp-Uq.1} into \eqref{esti-u-LqLp-Uq.0}
yields
\begin{align} \label{esti-u-LqLp-Uq}
   \|u\|_{L^q(I; L^p)}
   \leq C'(T) (\sum\limits_{j=1}^n|u_j|^q_{L^2})^{\frac 1q}
    \leq C'(T),
\end{align}
which implies \eqref{Hom-LpLq-Uq}.

In order to prove \eqref{Inhom-Vq-LpLq},
we see that
for any partition $\{t_j\}_{j=0}^m \in \calz$,
since $\|V(t_0,t_j)\|_{\mathcal{L}(L^2, L^2)} \leq C(T) \in L^\rho(\Omega)$,
if $\wt{f}:= \calx_I f$,
\begin{align} \label{esti-inhom.1}
    \bigg|\int_{t_{j-1}}^{t_{j}} V(t_0, s) \wt{f}(s) ds \bigg|_{L^2}
   \leq& C(T)  \bigg|\int_{t_{j-1}}^{t_{j}} V(t_j, s) \wt{f}(s) ds \bigg|_{L^2}  \nonumber \\
   =& C(T) \sup\limits_{|z|_{L^2}\leq 1} \<\int_{t_{j-1}}^{t_{j}} V(t_{j-1}, s) \wt{f}(s) ds, z\>_2 \\
   \leq& C(T) \|\wt{f}\|_{L^{q'}(t_{j-1}, t_{j}; L^{p'})} \sup\limits_{|z|_{L^2}\leq 1}
          \|V^*(t_{j-1}, \cdot) \calx_I z\|_{L^q(t_{j-1}, t_{j}; L^p)}.  \nonumber
\end{align}
Estimating as in \cite[Lemma 5.3]{Z17},
we have
\begin{align}  \label{Homo-Stri-LS-V*}
    \|V^*(t_{j-1}, \cdot) \calx_I z\|_{L^q(t_{j-1}, t_{j}; L^p) \cap L^2(t_{j-1}, t_{j}; H^\frac 12_{-1})}
    \leq C(T) |z|_{L^2}
    \leq C(T),
\end{align}
where $C(T) \in L^\rho(\Omega)$
for any $1\leq \rho <\9$.
Then, we get
\begin{align} \label{esti-inhom.2}
   \sum\limits_{j=1}^m
    \bigg|\int_{t_{j-1}}^{t_{j}} V(t_0, s) \wt{f}(s) ds \bigg|^{q'}_{L^2}
  \leq& (C''(T))^{q'} \sum\limits_{j=1}^m  \|\wt{f}\|^{q'}_{L^{q'}(t_{j-1}, t_{j}; L^{p'})} \nonumber  \\
  =&  (C''(T))^{q'} \|f\|^{q'}_{L^{q'}(I; L^{p'})},
\end{align}
where $C''(T)$ is independent of $(p,q)$
and $C''(T)\in L^\rho(\Omega)$,
$\forall 1\leq \rho<\9$.
This implies \eqref{Inhom-Vq-LpLq}.

Regarding \eqref{Inhom-V2-LpLq-LS},
let $f = f_1+f_2$
with $f_1\in L^{q'}(I;L^{p'})$,
$f_2 \in L^2(I; H^{-\frac 12}_1)$,
and set $\wt{f}_j =  \calx_I f_j$, $j=1,2$.
We can take a finer partition $\{t_j\}_{j=0}^m$
such that $\|f\|_{L^{q'}(t_{j-1}, t_{j}; L^{p'})} \leq 1$,
$1\leq j\leq m$.

Then,
estimating as in \eqref{esti-inhom.2},
since $q>2$, $q'<2$, we have
\begin{align*}
   \sum\limits_{j=1}^m
    \bigg|\int_{t_{j-1}}^{t_{j}} V(t_0, s) \wt{f}_1(s) ds \bigg|^{2}_{L^2}
  \leq& (C''(T))^{2} \sum\limits_{j=1}^m  \|\wt{f}_1\|^{2}_{L^{q'}(t_{j-1}, t_{j}; L^{p'})} \\
  \leq& (C''(T))^{2} \|f_1\|^{q'}_{L^{q'}(I; L^{p'})},
\end{align*}
which implies that
\begin{align} \label{esti-inhom-V2.1}
   \bigg\|\int_{t_0}^\cdot V(t_0, s) f_1(s) ds \bigg\|_{V^2(I)}
   \leq& C''(T) \|f_1\|^{\frac{q'}{2}}_{L^{q'}(I; L^{p'})}  \nonumber  \\
   \leq& C''(T) (1+ \|f_1\|_{L^{q'}(I; L^{p'})}).
\end{align}

Moreover,
arguing as in the proof of \eqref{esti-inhom.1}
and using \eqref{Homo-Stri-LS-V*}
we have
\begin{align*}
     \sum\limits_{j=1}^m
    \bigg|\int_{t_{j-1}}^{t_{j}} V(t_0, s) \wt{f}_2(s) ds \bigg|^2_{L^2}
    \leq  (C''(T))^2 \|f_2\|^2_{L^{2}(I; H^{-\frac 12}_1)}.
\end{align*}
This yields that
\begin{align}  \label{esti-inhom-V2.2}
    \bigg\|\int_{t_0}^\cdot  V(t_0, s) f_2(s) ds \bigg\|_{V^2(I)}
   \leq  C''(T) \|f_2\|_{L^{2}(I; H^{-\frac 12}_1)}.
\end{align}

Therefore, combining \eqref{esti-inhom-V2.1} and \eqref{esti-inhom-V2.2} together
we prove  \eqref{Inhom-V2-LpLq-LS}.
\hfill $\square$

In order to prove Theorem \ref{Thm-Sta-L2},
we first prove the short-time perturbation result below as in \cite{Z18}.
\begin{proposition}  \label{Pro-ShortP-L2} ({\it Mass-Critical Short-time Perturbation}).
Consider the situations in Theorem \ref{Thm-Sta-L2}.
Assume also  the smallness conditions
\begin{align}
    & |I|+ \|\wt{v}\|_{L^{q}(I; L^p)} \leq \delta,   \label{Short-L2-ve.1}\\
    & \| V(\cdot, t_0)(v(t_0)-\wt{v}(t_0) -R(t_0)) \|_{L^{q}(I; L^p)} \leq \ve, \
    \|R\|_{L^{q}(I; L^p) \cap L^1(I; L^2)} \leq \ve,\ \label{Short-L2-ve.2} \\
    & \| e \|_{{N}^0(I) + L^2(I; H^{- \frac 12}_{1})} \leq \ve   \label{Short-L2-ve.3}
\end{align}
for some $0<\ve\leq \delta$
where $\delta = \delta(C_T)>0$ is a small constant,
and $C_T$ is  as in Theorem \ref{Thm-Sta-L2}.
Then, we have
\begin{align}
   & \|v-\wt{v} -R\|_{L^{q}(I; L^p) \cap C(I; L^2)} \leq C(C_T) \ve, \label{Short-L2.1} \\
    & \|(F(v) - F(\wt{v})) + G (v-\wt{v} -R)) \|_{{N}^0(I)}
     \leq C(C_T) \ve,   \label{Short-L2.4}
\end{align}
where $(\delta(C_T))^{-1}$, $C(C_T)$
can be taken to nondecreasing with respect to $C_T$.
\end{proposition}

{\bf Proof.}
We mainly prove Proposition \ref{Pro-ShortP-L2}
for the case
where $p$ satisfies that
 $\frac 1p \in (\max\{\frac{1}{2\a}, \frac 12 - \frac{1}{2d}\}, \frac{1}{\a}(\frac 12 + \frac 1d))$
with  $1\leq d\leq 3$, $\a=1+\frac 4d$.
The case $p =2+\frac 4d$ with $d\geq 1$ can be proved similarly.

As mentioned below Hypothesis $(H0)^*$,
there exists another Strichartz pair $(\wt{p}, \wt{q})$
such that $(\frac{1}{\wt{p}'}, \frac{1}{\wt{q}'}) = (\frac{\a}{p}, \frac{\a}{q})$,
where
$q\in (2,\9)$ is such that $(p,q)$ is a Strichartz pair,
and $\wt{p}', \wt{q}'$ are the conjugate numbers of $\wt{p}, \wt{q}$ respectively.

Let $z:= v -\wt{v} - R$,
$F(\wt{v}):=|\wt{v}|^{\frac 4d} \wt{v}$
and $F(z +R+ \wt{v})$ be defined similarly.
By equations \eqref{equa-v-p} and \eqref{equa-wtv-p},
\begin{align} \label{equa-z-p}
   z(t) =& V(t,t_0)z(t_0)
          + \int_{t_0}^t V(t,s) \big(i (F(z + R+\wt{v}) - F(\wt{v})) + Gz  + GR + e \big) ds.
\end{align}

We set $S(I):= \|i  (F(z + R+\wt{v}) - F(\wt{v})) +Gz\|_{ N^0(I)}$.
By H\"older's inequality and \eqref{Short-L2-ve.1}-\eqref{Short-L2-ve.3},
\begin{align} \label{esti-S-L2}
   S(I) \leq& \|(F(z+R+\wt{v}) - F(\wt{v}))\|_{L^{\wt{q}'}(0,T; L^{\wt{p}'})}
              + \|Gz\|_{L^1(0,T;L^2)} \nonumber \\
        \leq& C(\|\wt{v}\|_{L^q(I; L^p)}^{\frac 4d} \|z + R\|_{L^q(I; L^p)} + \|z + R\|_{L^q(I; L^p)}^{1+\frac 4d})
              + \|Gz\|_{L^1(I; L^2)}   \nonumber \\
        \leq& C_1(\delta^{\frac 4d}\ve + \delta^{\frac 4d} \|z\|_{L^q(I; L^p)}  + \delta \|z\|_{C(I; L^2)}  + \|z\|_{L^q(I; L^p)}^{1+\frac 4d}),
\end{align}
where $C_1(\geq 1)$ depends on $d$ and $|G|_{L^\9}(I\times \bbr^d)$.
Moreover,
applying Theorem \ref{Thm-Stri} to \eqref{equa-z-p}
and using \eqref{Short-L2-ve.2} and \eqref{Short-L2-ve.3}  we have
\begin{align} \label{esti-S-L2*}
     &\|z\|_{L^q(I; L^p) \cap C(I; L^2)} \nonumber \\
     \leq& C_T (\|V(\cdot, t_0)z(t_0)\|_{L^q(I; L^p)} + S(I)
            + \|GR\|_{L^1(I; L^2)}
            + \|e\|_{N^0(I) + L^2(I; H^{- \frac 12}_{1})}) \nonumber \\
     \leq& C_2C_T (\ve + S(I)),
\end{align}
where $C_2$ depends on $|G|_{L^\9}(I\times \bbr^d)$.
Then, combining \eqref{esti-S-L2}, \eqref{esti-S-L2*} we get
\begin{align*}
      \|z\|_{L^{q}(I; L^p)\cap C(I; L^2)}
      \leq& C_1C_2C_T \big(2\ve +  (\delta^{\frac 4d}
                +  \delta) \|z\|_{L^q(I; L^p)\cap C(I; L^2)}   +   \|z\|_{L^q(I; L^p)}^{1+\frac 4d}\big).
\end{align*}
Thus, in view of Lemma $6.1$ in \cite{BRZ18},
taking $\delta = \delta(C_T)$ small enough such that
$C_1C_2C_T(\delta^{\frac 4d} + \delta) \leq \frac 12$
and  $4C_1C_2C_T \delta < (1-\frac 1\a)(2\a C_1C_2C_T)^{-\frac{1}{\a-1}}$
we obtain \eqref{Short-L2.1}.
Moreover,
plugging \eqref{Short-L2.1} into
\eqref{esti-S-L2} we obtain
\eqref{Short-L2.4}.

Therefore,
the proof is complete. \hfill $\square$ 

{\bf Proof of Theorem \ref{Thm-Sta-L2}.}
First fix $\delta = \delta(C_T)$ as in Proposition \ref{Pro-ShortP-L2}.
We divide $[t_0,T]$ into finitely many
subintervals
$[t'_j,t'_{j+1}]$,  $0\leq j\leq l'$,
such that
$t'_{j+1} = \inf\{t>t'_j; \|\wt{v}\|_{L^{q}(t'_j,t; L^p)} =\frac \delta 2\} \wedge T$.
Then,
$l'\leq (2L/\delta)^{2+\frac 4d}<\9$.

Take a new partition $\{I_j\}_{j=0}^l:= \{[t_j, t_{j+1}]\}_{j=0}^l$ of $[0,T]$,
such that
$\{t_j; 0\leq j\leq l+1\} = \{t'_j; 0\leq j\leq l'+1\}
      \cup \{t_0 + \frac{j\delta}{2}; 0\leq j\leq [\frac{2(T-t_0)}{\delta}]\}$.
Then,
\begin{align*}
    l\leq (2L/\delta)^{2+\frac 4d} + [2(T-t_0)/\delta], \ \
    |I_j| + \|\wt{v}\|_{L^{q}(I_j; L^p)} \leq \delta,\ \ 0\leq j\leq l.
\end{align*}

Let
$C(0) =  C(C_T)$,
$C(j+1) =  C(0)C_T^2 (\sum_{k=0}^j C(k)+ |G|_{L^\9}$ $+2)$,
$0\leq j\leq l-1$,
where $C(C_T)$ and $C_T$ are the constants
in  Proposition \ref{Pro-ShortP-L2} and Theorem \ref{Thm-Stri}, respectively.
Choose $\ve_*= \ve_*(C_T, L)$ sufficiently small such that
\begin{align} \label{ve*-L2}
  C_T^2 (\sum\limits_{k=0}^l C(k) +|G|_{L^\9} + 2) \ve_* \leq \delta.
\end{align}

We claim that
\eqref{Short-L2.1} and \eqref{Short-L2.4} hold on $I_j$
with $C(j)$ replacing $C(C_T)$
for every $0\leq j\leq l$.

To this end,
we first see that Proposition \ref{Pro-ShortP-L2} yields
the claim   for $j=0$.
Suppose that  the claim is also valid for each $0\leq k\leq j<l$.
We shall use  Proposition \ref{Pro-ShortP-L2} to show that
it also holds on $I_{j+1}$.

For this purpose,
applying Theorem \ref{Thm-Stri} to \eqref{equa-z-p} again
and using \eqref{Sta-L2-ve.2}, \eqref{ve*-L2}
and the inductive assumption we have
\begin{align*}
  &\|V(\cdot, t_{j+1})(v(t_{j+1}) - \wt{v}(t_{j+1}) -R(t_{j+1}))\|_{L^{q}(I_{j+1}; L^p)}  \\
  \leq&  \|V(\cdot, t_0)(v(t_0)-\wt{v}(t_0))\|_{L^{q}(I; L^p)}
         \\
      & + C^2_T  \big( S(t_0, t_{j+1})    +   \|GR\|_{L^1(t_0,t_{j+1};L^2)} +   \|e\|_{N^0(t_0, t_{j+1}) + L^2(t_0,t_{j+1}; H^{-\frac 12}_1)} \big) \\
  \leq&  \ve + C_T^2 (\sum\limits_{k=0}^j C(k) \ve + |G|_{L^\9(I\times \bbr^d)}\ve+\ve)
  \leq \delta.
\end{align*}
Thus, the conditions \eqref{Short-L2-ve.1}-\eqref{Short-L2-ve.3}  of Proposition \ref{Pro-ShortP-L2} are satisfied
with
$C_T^2$ $(\sum_{k=0}^j C(k)+|G|_{L^\9}+2)\ve$
replacing  $\ve$.
Applying Proposition \ref{Pro-ShortP-L2}
we prove the claim on $I_{j+1}$.

Therefore, using inductive arguments
we prove the claim on $I_j$   for every $0\leq j\leq l$,
thereby proving Theorem \ref{Thm-Sta-L2}.
The proof is complete.
\hfill $\square$

{\bf Proof of \eqref{Integ-M12}.}
The proof is similar to that in \cite{FX18.1}.
Without loss of generality,
we may assume $\rho \geq q$.
By Minkowski's inequality
and Burkholder's inequality,
\begin{align} \label{esti-M1}
  \|M_1^*\|_{L^\rho(\Omega; L^1(0,T))}
  \leq& C \|M_1^*\|_{L^1(0,T; L^\rho(\Omega))}  \nonumber \\
  \leq& C \bigg\| \bigg|\int_0^\cdot \|e^{-i(\cdot -s)\Delta} X(s)\Phi\|^2_{HS(\bbr^N; L^2)} ds\bigg|^{\frac 12} \bigg \|_{L^1(0,T)},
\end{align}
where the operator $\Phi: \bbr^N \mapsto L^\9$ is defined by
$\Phi(x) = \sum_{j=1}^N \mu_j e_j x_j$
for $x= (x_1,\cdots, x_N)$,
and $\| \cdot \|_{HS(\bbr^N; L^2)}$ denotes the Hilbert-Schmidt norm
(see, e.g., \cite{LR15}).
Note that,
\begin{align*}
   \|e^{-i(t -s)\Delta} X(s)\Phi\|^2_{HS(\bbr^N; L^2)}
   \leq& C \sum\limits_{j=1}^N |\mu_j|^2|e_j|^2_{L^\9} |X_0|^2_{L^2}.
\end{align*}
Plugging this into \eqref{esti-M1}
yields immediately that $\|M_1^*\|_{L^\rho(\Omega; L^1(a,b))}  <\9$.

Similarly,
again by Minkowski's inequality
and Burkholder's inequality,
\begin{align} \label{esti-M2}
  \|M_2^*\|_{L^\rho(\Omega; L^q(0,T))}
  \leq& \|M_2^*\|_{L^q(0,T; L^\rho(\Omega))}  \nonumber \\
  \leq& C \bigg\| \bigg|\int_0^\cdot \|e^{-i(\cdot -s)\Delta} X(s)\Phi\|^2_{\mathcal{R}(\bbr^N; L^\rho)} ds\bigg|^{\frac 12} \bigg \|_{L^q(0,T)},
\end{align}
where $\| \cdot \|_{\mathcal{R}(\bbr^N; L^\rho)}$ denotes the $\g$-radonifying norm
(see, e.g., \cite{BD99,BD03,FX18.1}).

Note that,
by the dispersive inequality
$\|e^{-i(t-s)\Delta}\|_{\mathcal{L}(L^{p'}, L^p)} \leq C (t-s)^{-\frac{d}{2}(1-\frac 2p)}$,
\begin{align*}
   \|e^{-i(t -s)\Delta} X(s)\Phi\|_{\mathcal{R}(\bbr^N; L^\rho)}
   \leq& \|\Phi\|_{\mathcal{R}(\bbr^N; L^{\frac{2p}{p-2}})}
         \|X(s)\|_{\mathcal{L}(L^{\frac{2p}{p-2}}, L^{p'})}
         \|e^{-i(t-s)\Delta}\|_{\mathcal{L}(L^{p'}; L^p)}   \\
   \leq& C(N)|X_0|_{L^2} (\sum\limits_{j=1}^N |\mu_j|^2 |e_j|^2_{L^{\frac{2p}{p-2}}})^\frac 12
          (t-s)^{-\frac{d}{2}(1-\frac 2p)}.
\end{align*}
Plugging this into \eqref{esti-M2} yields that
$ \|M_2^*\|_{L^\rho(\Omega; L^q(0,T))} <\9$ if $p<\frac{2d}{d-1}$,
thereby finishing the proof.
\hfill $\square$

{\it \bf Acknowledgements.}
The author would like to thank Yiming Su for helpful discussions on tightness.
The author also thanks Chenjie Fan for conversations on integrability of
stochastic controlled solutions in the mass-critical case.
Financial support by the NSFC (No. 11871337) is gratefully acknowledged.

\end{document}